\documentclass[hidelinks,onefignum,onetabnum]{siamart250106}

\usepackage{graphicx} \usepackage{booktabs}
\usepackage{amsfonts}
\usepackage{amsmath}
\usepackage{amssymb}
 \usepackage{array}
\usepackage{mathtools}
\usepackage{colortbl} \usepackage{multicol}
\usepackage{enumitem}
\usepackage{bm}
\usepackage{comment}
\usepackage{caption}
\usepackage{pgfplots}
\definecolor{CeruleanRef}{RGB}{12,127,172}
\hypersetup{colorlinks=true,allcolors=CeruleanRef}
\hypersetup{pdfstartview=FitB,pdfpagemode=UseNone}
\usepackage[nameinlink]{cleveref}
\crefname{equation}{}{}
\usepackage{rotating}
      \usepackage{standalone}
\usetikzlibrary{trees}
\usepackage{subcaption}
\usepackage{hhline}
\usepackage{multirow}
 \usepackage{wrapfig}

\newcommand{\dd}{\mathop{}\!\mathrm{d}}

\crefname{equation}{}{}

\usepackage{todonotes}

\DeclareMathOperator*{\argmin}{arg\, min}
\DeclareMathOperator*{\C}{\mathsf{C}}
\DeclareMathOperator*{\tr}{tr}
\DeclareMathOperator*{\tanhm}{tanhm}

\makeatletter
\newcommand{\fcolon}{  \mathrel{\mathpalette\fcolon@\relax}}
\newcommand{\fcolon@}[2]{  \sbox\z@{$\m@th#1:$}  \vbox to\ht\z@{    \hbox{$\m@th#1.$}    \vss
    \hbox{$\m@th#1.$}    \vss
    \hbox{$\m@th#1.$}  }}
\makeatother

\title{The latent variable proximal point algorithm for variational problems with inequality constraints}
\author{J{\O}rgen S.\ Dokken\thanks{Department of Scientific Computing and Numerical Analysis, Simula Research Laboratory, Oslo, Norway} (\email{dokken@simula.no}, \email{thomasms@simula.no})
\and 
Patrick E.\ Farrell\thanks{Mathematical Institute, University of Oxford, Oxford, United Kingdom
and
Mathematical Institute, Faculty of Mathematics and Physics, Charles University, Prague, Czechia}
(\email{patrick.farrell@maths.ox.ac.uk})
\and
Brendan Keith\thanks{Division of Applied Mathematics, Brown University, Providence, RI, United States of America}
(\email{brendan\_keith@brown.edu})
\and
Ioannis P.\ A.\ Papadopoulos\thanks{Weierstrass Institute, Berlin, Germany}
(\email{papadopoulos@wias-berlin.de})
\and
Thomas M.\ Surowiec\footnotemark[1]
}

\begin{document}

\maketitle

\begin{abstract}
The latent variable proximal point (LVPP) algorithm is a framework for solving infinite-dimensional variational problems with pointwise inequality constraints. The algorithm is a saddle point reformulation of the Bregman proximal point algorithm. At the continuous level, the two formulations are equivalent, but the saddle point formulation is more amenable to discretization because it introduces a structure-preserving transformation between a latent function space and the feasible set. Working in this latent space is much more convenient for enforcing inequality constraints than the feasible set, as discretizations can employ general linear combinations of suitable basis functions, and nonlinear solvers can involve general additive updates. LVPP yields numerical methods with observed mesh-independence for obstacle problems, contact, fracture, plasticity, and others besides; in many cases, for the first time. The framework also extends to more complex constraints, providing means to enforce convexity in the Monge--Amp\`ere equation and handling quasi-variational inequalities, where the underlying constraint depends implicitly on the unknown solution. In this paper, we describe the LVPP algorithm in a general form and apply it to ten problems from across mathematics.
\end{abstract}

\section{Introduction} \label{sec:intro}

Many problems in science, engineering, finance, and mathematics involve solving for a function subject to inequality constraints. Prominent examples include contact and damage in solid mechanics \cite{NKikuchi_JTOden_1988,Bourdin2000}, non-negativity of probability densities \cite{Risken1996-wo}, invariant domain properties of hydrodynamic flow models \cite{wu2023geometric}, the pricing of American options \cite{jaillet1990}, and the convexity condition on solutions to the Monge--Amp\`ere equation \cite{santambrogio2015}.
While there are many popular methods for such problems, in general, it is quite difficult to achieve \emph{mesh-independence}---that the number of iterations required for convergence of an optimization solver or Newton-type method does not grow unboundedly as the discretization is refined \cite{allgower1986mesh,schwedes2017mesh,weiser2005asymptotic}---or \emph{high-order accuracy}. 

Where mesh-independent algorithms exist, they often return infeasible iterates \cite{Glowinski2014-ii,NKikuchi_JTOden_1988} or are tightly tied to lowest-order discretizations (e.g.,~piecewise linear finite elements), as they rely on monotonicity properties of the basis functions employed \cite{bueler2023}.
Typically, they also require parameters to pass to infinity for convergence \cite{hintermuller2006feasible,hintermueller2006b,adam2019semismooth}.
While some methods have been proposed using high-degree basis approximations to construct solutions (e.g., Bernstein polynomial-based discretizations \cite{Kirby2010} or biorthogonal bases \cite{banz2015}), these approaches are limited to pointwise bound constraints, and are generally mesh-dependent.
The first main challenge is the general lack of smoothness of many variational problems with inequality constraints.
For example, the optimality conditions for the classical obstacle problem \cite{DKinderlehrer_GStampacchia_2000} are, in general, insufficiently regular to define a semismooth Newton method \cite{hintermueller2002,ulbrich2003,ulbrich2011,hintermuller2006feasible} at the continuous level. This can seemingly be avoided on the finite-dimensional level by first discretizing the problem, where semismooth Newton methods are directly related to certain active set solvers for complementarity problems \cite{munson2001,hintermueller2002}. However,
the algorithm's lack of validity at the continuous level manifests as mesh dependence \cite[\S4]{farrell2018a} as demonstrated in \Cref{sub:ObstacleProblem} of this paper.

The second main challenge is that the feasible set $K$, implied by the underlying inequality constraint(s), is not a vector space; therefore, one cannot take arbitrary linear combinations of functions in $K$ and remain in the set.
Most successful discretization approaches approximate solutions using linear combinations of basis functions, but this strategy conflicts with the geometry imposed by the constraints. Moreover, most practical algorithms for solving optimization problems or nonlinear equations involve additive updates, which again conflicts.

In this paper, we introduce a general framework that adaptively regularizes variational problems with inequality constraints to overcome these challenges.
The employed regularization permits $K$ to be identified with a latent function space via a flexible class of smooth coordinate transformations arising from the set's convex geometry.
In turn, the framework enables the natural use of linear combinations of basis functions for approximation and additive update formulae.
Furthermore, it results in regularized subproblems that are smooth enough to be solved with well-established mesh-independent techniques.
We call this framework the \emph{latent variable proximal point algorithm} (LVPP) and demonstrate some of its many potential applications by solving ten challenging problems from across mathematics:
\vspace*{-0.9em}
\begin{multicols}{2}
\begin{itemize}[leftmargin=*]
    \item The obstacle problem
    \item The Signorini problem
    \item Variational fracture
    \item Multi-phase gradient flow
    \item Quasi-variational inequalities
    \item Gradient constraints
    \item Eigenvalue constraints
    \item Intersections of constraints
                \item The eikonal equation
    \item The Monge--Amp\`ere equation
\end{itemize}
\end{multicols}
\vspace*{-0.9em}
The LVPP algorithm, which was first proposed by Keith and Surowiec in \cite{keith2023proximal}, is highly adaptable to various problem types, and possesses the following key properties, among others:
\begin{enumerate}[leftmargin=*,label=(\roman*)]
\item an infinite-dimensional formulation;
\item observed mesh and order independence in both the mesh size and discretization degree;
\item a simple mechanism for enforcing pointwise constraints on the discrete level without the need for a projection;
\item ease of implementation --- the algorithm reduces to repeatedly solving a system of smooth, nonlinear PDEs without requiring specialized discretizations, making it compatible with high-order methods and standard finite element libraries;
\item robust numerical performance, as convergence does not rely on any parameter passing to infinity.
\end{enumerate}
No other framework satisfies all these properties simultaneously. For instance, penalty methods \cite{hintermuller2006feasible,hintermueller2006b,adam2019semismooth} satisfy (i), (ii), and (iv) but, without special care, often violate (iii) and (v), deliver infeasible solutions, and are known to be sub-optimal for high-order discretizations \cite[\S3.1]{gustafsson2017on}. Similarly, although multigrid methods \cite{hackbusch1983, hoppe1987,kornhuber1994} satisfy (iii) and (v), they are mildly mesh dependent, tied to low-order discretizations, and the implementation may be delicate for general problems.

\begin{figure}[ht]
\centering
 \begin{tikzpicture}[%
every node/.style={draw=black, thick, anchor=west},
grow via three points={one child at (0.0,-0.7) and
    two children at (0.0,-0.7) and (0.0,-1.4)},
edge from parent path={(\tikzparentnode.210) |- (\tikzchildnode.west)}]
\node {Variational problem with inequality constraints}
child { node {Apply LVPP: sequence of nonlinear systems of PDEs \cref{eq:intro:lvpp}}
    child { node {Discretize \cref{eq:intro:lvpp}: sequence of nonlinear systems of equations}
        child { node{Newton solver} 
        }
    }
};
\end{tikzpicture}
\caption{LVPP approximates a solution of a variational problem with inequality constraints by solving a sequence of nonlinear systems of PDEs. In practice these systems are discretized and the resulting nonlinear systems of equations are solved via Newton's method.}\label{fig:solver-diagram}
\end{figure}
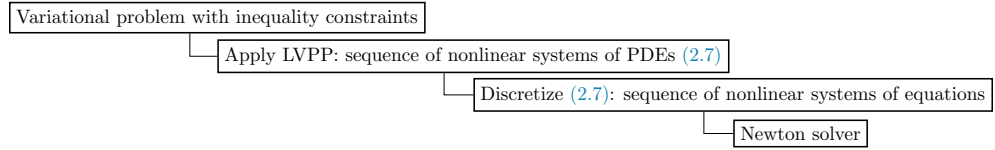

In \Cref{sec:LVPPalgorithm}, we describe the theoretical foundations of the LVPP algorithm.
\noindent\noindent\Cref{sec:BoundConstraints} is dedicated to applying LVPP to pointwise bound constraints and includes example applications for the obstacle problem, contact mechanics, variational fracture, multiphase free boundary problems, and obstacle-type quasi-variational inequalities.
\Cref{sec:Non-Polyhedral} focuses on more complicated convex constraints, such as norm and eigenvalue constraints, as well as how to handle multiple simultaneous inequality constraints.
Finally, in \Cref{sec:NonlinearPDEs}, we describe applications to solving fully-nonlinear first- and second-order PDEs, using the eikonal and Monge--Amp\`ere equations as examples.
Our findings are summarized in \Cref{sec:Conclusion}. In all the examples, we derive our numerical methods via the solver diagram blueprint found in \cref{fig:solver-diagram}.

Supplementing the ten inequality constraint applications in this work, \Cref{app:LinearEqualityConstraints} considers problems with general linear \emph{equality} constraints.
Here, we show that LVPP reduces, in the appropriate limit, to the classical, Lagrange multiplier saddle-point formulation, as established and studied thoroughly in, e.g.,~\cite{boffi2013mixed}.

\section{Theoretical foundations}\label{sec:LVPPalgorithm}

We formally derive the general LVPP algorithm, drawing together classical notions from convex analysis and optimization, such as Legendre functions and the Bregman proximal point method.
The section closes with a brief discussion of proximal numerical methods derivable from LVPP, including the proximal Galerkin method \cite{keith2023proximal}. {\color{black} It is important to note that while the majority of statements made in this section are rigorous, several key steps leading to the \textit{general} form of LVPP have only been proven to date for obstacle-type constraints, see \cite[Appendix A]{keith2023proximal}. Beyond this setting, and to be completely rigorous at the continuous level, every new application of LVPP introduced below requires a proof of uniform bounds on the LVPP iterates that ensure they are strictly feasible in a quantifiable way. Such results are crucial for the well-posedness of the LVPP subproblems and, ultimately, for linking the saddle point problems back to the original proximal point method. The bounds themselves allow us to differentiate the nonlinear functionals and operators in the subproblems. Though these proofs go beyond the scope and purpose of this article, we point the interested reader to Appendices A1, A2, and A3 in \cite{keith2023proximal} for a ``roadmap'' detailing what is necessary.
}

We choose to focus on variational problems posed on Lipschitz domains $\Omega \subset \mathbb{R}^n$ with feasible sets $K$ of the following form:
\begin{equation} \label{eq:intro:feasible_set}
K = \{ v \in V \mid Bv(x) \in C(x) \text{ for almost every } x \in \Omega_d \subset \overline{\Omega} \}.
\end{equation}
Here, $\Omega_d$ is a subset of the closure of $\Omega$, $V \subset L^1(\Omega)$ is a Banach space on $\Omega$, and $B \colon V \to L^1(\Omega_d)$ is a bounded linear operator on $V$ whose arguments in~\cref{eq:intro:feasible_set} are constrained to map pointwise a.e.\ in $\Omega_d$ into a convex set $C(x) \subset \mathbb{R}^m$.
Generally speaking, $\Omega_d$ should be Hausdorff-measurable, with Hausdorff dimension $d$ not necessarily equal to $n$. For example, in contact problems, we may wish to take $\Omega_d \subset \partial \Omega$, where $d = n - 1$.
Moreover, $C$ may vary from point to point, as in an obstacle problem, where $B = \operatorname{id}$ is the identity operator and $\Omega_d = \Omega$, specified functions give the upper and lower bounds on every $v \in K$. We refer to $C = C(x)$ as the (pointwise) \emph{feasible image} of $B$.

\subsection{Legendre functions}

The geometry of a closed convex set $C \subset \mathbb{R}^m$ with a non-empty interior, $\operatorname{int}C \neq \emptyset$, can be encoded into a \emph{Legendre function} $R: \mathbb{R}^m \to \mathbb{R}\cup\{+\infty\}$. We recall the definition here.
\begin{definition}[Legendre function]
\label{def:LegendreFunction}
Let the essential domain of a function $R$ be defined as $\operatorname{dom} R \coloneqq \{ a \in \mathbb{R}^m \mid R(a) < \infty \}$. We call a proper convex function a \emph{Legendre function} $R: \mathbb{R}^m \to \mathbb{R}\cup\{+\infty\}$ if
\begin{itemize}[leftmargin=*,itemsep=1pt,topsep=2pt]
\item $\operatorname{int}(\operatorname{dom} R) \neq \varnothing$;
\item R is differentiable on $\operatorname{int}(\operatorname{dom} R)$;
\item $\lim_{t \to 0^+} \langle \nabla R(a+t(b-a)), b-a\rangle = -\infty$ for all $a \in \partial( \operatorname{dom} R)$ and $b \in \operatorname{int}(\operatorname{dom} R)$;
\item R is strictly convex on $\operatorname{int}(\operatorname{dom} R)$.
\end{itemize}
\end{definition}
Introduced by Rockafellar in 1967~\cite{Rockafellar1967Conjugates}, see also \cite[Chap.~26]{RTRockafellar_1970}, Legendre functions constitute a special class of proper convex functions whose gradients $\nabla R$ become singular on the boundary of their essential domains.
We denote by $R^*(a^\ast) \coloneqq \sup \{ a \cdot a^\ast - R(a) \mid a \in \mathbb{R}^m \}$ the convex conjugate of $R$.

\begin{theorem}[Rockafellar \text{\cite[Thm. 1]{Rockafellar1967Conjugates}}]
\label{thm:Rockafellar}
A proper convex function $R$ is a Legendre function if and only if its convex conjugate $R^\ast$ is also a Legendre function.
Moreover, $\nabla R \colon \operatorname{int}(\operatorname{dom} R) \to \operatorname{int}(\operatorname{dom} R^*)$ is a topological isomorphism with $(\nabla R)^{-1} = \nabla R^\ast$.
\end{theorem}

\begin{figure}
    \centering
    \includegraphics[width=0.55\linewidth]{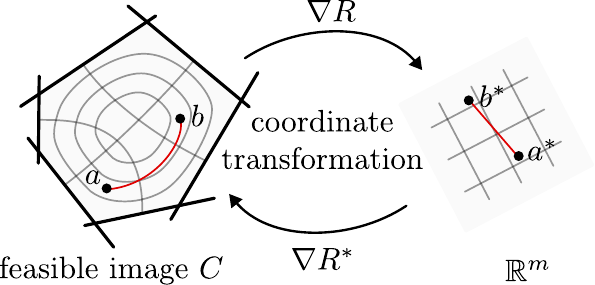}
    \caption{The gradient of the Legendre function $R$ is an isomorphism between the interior of the feasible image of an operator $B$, denoted $\operatorname{int} C$, and $\mathbb{R}^m$.
    Here, curved geodesics in $C$ are transformed into straight lines in $\mathbb{R}^m$ \cite{Amari2016}, with the Riemannian metric on $\operatorname{int} C$ given by the Hessian $\nabla^2 R$ \cite{alvarez2004hessian}.
    The grey lines in the illustration depict the different coordinate systems on the two manifolds, while the red lines represent geodesics between points $a$ and $b$ in $C$  and dual geodesics between $a^\ast \coloneqq \nabla R(a)$ and $b^\ast \coloneqq \nabla R(b)$ in $\mathbb{R}^m$.
    By \Cref{thm:Rockafellar}, the gradient of the convex conjugate $R^\ast$ is the inverse mapping $\nabla R^\ast = (\nabla R)^{-1}$.
    }
    \label{fig:Geodesics}
\end{figure}

From now on, we assume that ${R(a)}/{|a|} \to +\infty$ as $|a| \to \infty$ and $\operatorname{dom} R = C$.\footnote{In the more general setting, where $C = C(x)$ varies with $x \in \Omega$, it is necessary for $R = R(x,a)$ to be a Carath\'eodory function. In this case, $\nabla$ denotes the gradient with respect to the second argument and $R^\ast(x,a^\ast) = \sup \{ a \cdot a^\ast - R(x,a) \mid a \in \mathbb{R}^m \}$. For further details, see \cite[Section~3]{fu2024locallyconservativeproximalgalerkinmethod}.}
The first property is equivalent to $\operatorname{dom} R^* = \mathbb{R}^m$ (see, e.g., \cite[Corollary 13.3.1]{RTRockafellar_1970} and \cite[Proposition 2.16]{bauschke1997legendre}), and so \Cref{thm:Rockafellar} ensures that the gradient $\nabla R$ is an isomorphism between $\operatorname{int} C$ and $\mathbb{R}^m$.
Geometrically, this map defines geodesics in $C$ that appear as straight lines in $\mathbb{R}^m$ \cite{Amari2016}; see \Cref{fig:Geodesics} for a depiction.
As explained below, the inverse of the gradient of $R$, namely, $\nabla R^\ast \colon \mathbb{R}^m \to \operatorname{int} C$, exposes a latent vector space enabling simple, high-order discretizations of a wide variety of variational problems with inequality constraints.

The Legendre function $R$ induces a non-symmetric notion of (squared) distance within the feasible image $C$ via the concept of a Bregman divergence~\cite{BREGMAN1967200}:
\begin{equation}
\label{eq:BregmanDivergence}
D_R(a, b) \coloneqq R(a) - R(b) - \nabla R(b) \cdot (a-b)\,, \quad a \in C\,,~b \in \operatorname{int} C\,,
\end{equation}
which measures the error in the first-order Taylor expansion of $R$ over $C$. In turn,~\cref{eq:BregmanDivergence} induces a notion of distance in the feasible set $K$, given by $\int_{\Omega_d} D_R(Bu, Bv) \dd \mathcal{H}_d$, where $\dd \mathcal{H}_d$ denotes the $d$-dimensional Hausdorff measure on $\Omega_d$.
In this work, we only consider the settings $\Omega_d = \Omega$ and $\Omega_d \subset \partial\Omega$, in which case $\dd \mathcal{H}_d = \dd \mathcal{H}_{n}$ and $\dd \mathcal{H}_d = \dd \mathcal{H}_{n-1}$ coincide with the standard Lebesgue volume and surface measures, respectively.

\subsection{The latent variable proximal point algorithm}

The Bregman divergence induced by a Legendre function has a regularizing effect on constrained optimization problems of the form
\begin{equation} \label{eq:intro:constrained}
\min_{u \in K} J(u)
\,,
\end{equation}
that is exploited in the \emph{Bregman proximal point algorithm}~\cite{MTeboulle_2018}: for a suitable initial guess $u^0 \in V$ and given a positive sequence of proximity parameters $\{\alpha_k\}$, compute
\begin{equation} \label{eq:intro:primalbregman}
u^{k} \in \argmin_{u \in K} J(u) + \alpha_{k}^{-1} \int_{\Omega_d} D_R(Bu, Bu^{k-1}) \dd \mathcal{H}_d, \quad k = 1, 2, \dots .
\end{equation}
{\color{black}
This algorithm adaptively regularizes~\cref{eq:intro:constrained}, with the striking property that it converges arbitrarily quickly for convex $J$:
\begin{equation} \label{eq:intro:striking_convergence}
J(u^k) - J(u^\star) \le \bigg(\sum_{i=1}^k \alpha_i\bigg)^{-1}\int_{\Omega_d} D_R(Bu^\star, Bu^0) \dd \mathcal{H}_d,
\end{equation}
where $u^\star$ is any minimizer of $J$ in $K$ \cite{chen1993convergence}.
If $J$ is strongly convex, then $u^\star$ is unique, and we obtain strong convergence of $u^k \to u^\star$ in $V$. Since $\alpha_k$ does not need to go to infinity for convergence, the algorithm can be very robust; the algorithm can also be very fast if the sequence grows rapidly. In addition to convex problems, we also consider non-convex functionals (\Cref{sec:fracture}) and quasi-variational inequalities (\Cref{sec:qvi}) in this work.
}

{\color{black}Solving \cref{eq:intro:primalbregman} generally requires gradient-based methods and thus, differentiability of the Bregman divergence.}
To this end, the singularity of $\nabla R$ on $\partial C$ leads us to a condition familiar to interior point methods: $B u^k(x) \in \operatorname{int} C(x)$ for a.e.\ point $x \in \Omega_d$.
{\color{black}For bound constraints, this inclusion is even uniformly satisfied, i.e., the distance of $B u^k$ to the bound of $C$ is uniformly greater than a positive constant dependent on $\alpha$, see \cite{keith2023proximal}.}
In turn, each perturbed subproblem~\cref{eq:intro:primalbregman} {\color{black}will have} first-order optimality conditions that are smooth \emph{equations}. In contrast, the first-order optimality conditions for traditional quadratic penalty or augmented Lagrangian methods naturally include nonsmooth operators. Thus, the first order optimality conditions for  \cref{eq:intro:primalbregman}, amount to
\begin{equation} \label{eq:intro:fooc}
    \text{find } u^k \in K: \quad \alpha_k J^\prime(u^k) + B^* \nabla R(Bu^k) - B^* \nabla R(Bu^{k-1}) = 0
    ~~\text{ in } V^\prime
    ,
\end{equation}
where $J^\prime$ is the Fr\'echet derivative of $J$ and $B^*$ is the dual (conjugate) operator of $B$; {\color{black}cf.\ the derivation of the \emph{entropic Poisson equation} in \cite{keith2023proximal}}.

While the Bregman proximal point algorithm is theoretically attractive, its practical implementation \cref{eq:intro:fooc} requires discretizing functions $u^k$ in the feasible set $K$. The key novelty of the \textit{latent variable proximal point} (LVPP) \textit{algorithm} is to rewrite the Bregman proximal point algorithm \cref{eq:intro:primalbregman} as a sequence of saddle point problems: starting from some $\psi^0 \in W$, find $(u^{k}, \psi^{k}) \in V \times W$ satisfying
\begin{subequations} \label{eq:intro:lvpp}
\begin{align}
\label{eq:intro:lvpp:a}
\alpha_{k} J'(u^{k}) + B^*\psi^{k} &= B^*\psi^{k-1}, \\
\label{eq:intro:lvpp:b}
Bu^{k} - \nabla R^*(\psi^{k}) &= 0,
\end{align}
\end{subequations}
for $k = 1,2, \dots$, where $\psi^k \coloneqq \nabla R(Bu^k)$ is a \emph{latent variable} lying in a suitable Banach space $W$ over which $\nabla R^\ast \colon W \to B(K)$ is well-defined and continuously invertible onto its range.
In this paper, we adopt the simplifying conventions $\psi^0 = 0$ and $W = L^\infty(\Omega_d)$. 
However, the identity $W = L^\infty(\Omega_d)$ can only be guaranteed by technical regularity conditions that are outside the scope of this work; cf.\ \cite{keith2023proximal}.
In general, problems \cref{eq:intro:primalbregman} and \cref{eq:intro:lvpp} are equivalent in that $u^{k}$ coincide for each iteration $k$.
However, the saddle point formulation~\cref{eq:intro:lvpp} yields two approximations to the pointwise observable $Bu$; namely, $Bu^{k}$ and $\nabla R^*(\psi^{k})$, which differ from each other after discretization.
Since convergence of the algorithm is proven at the infinite-dimensional level, one can generically hope for mesh-independent convergence of its discretizations.

The saddle point reformulation~\cref{eq:intro:lvpp} has at least four major advantages.
First, \emph{this formulation does not require discretizing the feasible set $K$}.
Instead, it only requires discretizing the vector spaces $V \supset K$ and $W$.
Thus, LVPP removes any need for special discretizations, e.g., positive basis functions, enabling the use of convenient and familiar discretizations with well-understood approximation properties, e.g., $hp$-FEM \cite{papadopoulos2024hierarchicalproximalgalerkinfast} and high-degree spectral approximations for the obstacle problem (cf.\ \Cref{sub:ObstacleProblem}). Second, LVPP is guaranteed to return discretized iterates in $K$.
Indeed, although the approximate observable $Bu^{k}_h$ may not belong to $K$ after discretizing \cref{eq:intro:lvpp}, the second approximation $\nabla R^*(\psi_h^{k})$ always does by construction (recall $\nabla R^\ast \colon \mathbb{R}^m \to \operatorname{int} C$).
Third, the subproblems~\cref{eq:intro:lvpp} are often smoother than those arising in other algorithms and can be solved by standard Newton methods.
Finally, because of this smoothness, the number of required iterations is typically bounded as the resolution of the discretization increases.

LVPP was first introduced in the context of pointwise bound constraints ($B = \operatorname{id}$, $\Omega_d = \Omega$) \cite{keith2023proximal}.
In this work, we make several major extensions, allowing us to systematically tackle a wide variety of variational problems with gradient ($B = \nabla$, $\Omega_d = \Omega$), Hessian ($B = \nabla^2$, $\Omega_d = \Omega$), and trace constraints ($B = \gamma$, $\Omega_d = \partial\Omega$), where $\gamma$ is the trace operator.
See \Cref{tab:ConstraintTable} for a partial list.
As a general principle, we choose $R$ so that $\nabla R^\ast$ is Fr\'echet differentiable.
For example, in our treatment of the obstacle problem presented in \Cref{sub:ObstacleProblem}, we consider taking $R(a) = (a - \phi) \ln (a - \phi) - (a - \phi)$, leading to the exponential map $\nabla R^\ast(\psi^k) = \phi + \exp\psi^k$ in~\cref{eq:intro:lvpp}, which is infinitely Fr\'echet differentiable on $W = L^\infty(\Omega)$.
This particular setting was investigated in detail in \cite{keith2023proximal}.
However, we note that the choice of the Legendre function is flexible in LVPP, and often, many appealing choices can be found for solving the same problem.

\begin{table}
\centering
\small
\setlength{\tabcolsep}{5pt}
\renewcommand{\arraystretch}{1.4}
    \begin{tabular}{ c|c|c|c } 
     \toprule
      Feasible set $K$ & Legendre function $R$ & $B$ & $\nabla R^*(\psi)$ \\ 
     \midrule
     $\big\{ u \geq \phi \big\}$ & $(a - \phi) \ln (a - \phi) - (a - \phi)$ & $\operatorname{id}$ & $\phi + \exp\psi$ \\[2ex]
     $\big\{ \phi_1 \leq u \leq \phi_2 \big\}$ & $(a - \phi_1) \ln (a - \phi_1) + (\phi_2-a) \ln (\phi_2-a)$ & $\operatorname{id}$ & $\dfrac{\phi_1 + \phi_2\exp\psi}{1 + \exp\psi}$ \\[2ex]
          $\big\{\gamma u \ge \phi \big\}$ & $(a-\phi) \ln (a - \phi) - (a - \phi)$ & $\gamma$ &  $\phi + \exp\psi$ \\[2ex]
     $\big\{ (\gamma u)\cdot n \leq \phi \big\}$ & $(\phi-a) \ln (\phi - a) - (\phi - a)$ & $\gamma(\cdot) \cdot n$ & $\phi - \exp(-\psi)$ \\[2ex]
     $\big\{| \nabla u | \le \phi \big\}$ & $ -\sqrt{\phi^2 - | a |^2}$ & $\nabla$  & $\dfrac{\phi\psi}{\sqrt{1 + | \psi |^2}}$ \\[3ex] 
     $\big\{ u \ge 0,\; \sum_{i} u_i = 1 \big\}$ & $\sum_{i} a_i \ln(a_i)$ & $\operatorname{id}$ & $\dfrac{\exp\psi}{\sum_{i} \exp\psi_i}$\\[2ex]
     $\big\{ \det (\nabla^2 u) \geq 0 \big\}$ & $\operatorname{tr}(a\ln a - a)$ & $\nabla^2$ & $\exp \psi$ \\[2ex]
     \bottomrule
    \end{tabular}
                                                                                \caption{\label{tab:ConstraintTable} Some examples of the feasible set~\cref{eq:intro:feasible_set} with an associated Legendre function $R$. Note that the Legendre function choices are not unique, and many are often available for the same feasible set $K$. Here, $n$ is the unit outward normal to $\partial \Omega$, $\operatorname{id}$ is the identity operator, $\gamma$ is the trace operator, and $\nabla^2$ is the Hessian operator.
        The function $u$ can be scalar- or vector-valued, as appropriate.
    Likewise, $\exp$ denotes the scalar, component-wise, or matrix exponential function.
    }
\vspace{-1em}
\end{table}

\subsection{Proximal Galerkin and other proximal numerical methods}
\label{sub:PG}
To\linebreak date, the LVPP algorithm has mainly been used as a means of deriving proximal Galerkin finite element methods \cite{keith2023proximal,fu2024locallyconservativeproximalgalerkinmethod,papadopoulos2024hierarchicalproximalgalerkinfast}, which arise from Galerkin discretizations of the LVPP saddle point problems~\cref{eq:intro:lvpp} using conforming finite element spaces.
However, since LVPP is derived at the continuous level, it is agnostic to the method of discretization.
Therefore, \textit{LVPP should be seen as a framework for deriving numerical methods rather than a numerical method itself.}
We highlight this perspective in \Cref{sub:ObstacleProblem} below by introducing other proximal numerical methods in addition to proximal Galerkin.
Nevertheless, the majority of methods in this work are proximal Galerkin methods.

No exotic elements are needed for the finite element methods in this paper.
Although finite element discretizations of saddle-point problems can come in many forms, with different properties depending on chosen bases or subspaces, we choose to focus on elements discretizing the $L^2$ de Rham complex (continuous Lagrange, etc.)~\cite{arnold2018}.
These elements are widely available in free software, meaning that every derived method should be implementable by non-experts.

\subsection{Update strategies for the proximal parameter $\alpha_k$}

The sequence of proximal parameters $\alpha_k$ in \cref{eq:intro:lvpp:a} is chosen by the user. As was shown for the obstacle problem \cite[Theorem 4.13]{keith2023proximal}, LVPP converges to the solution of the obstacle problem even for a fixed constant sequence (e.g.,~$\alpha_k = 1$ for all $k \in \mathbb{N}$). We conjecture that the same convergence result is true for all the examples in this work. For a constant sequence, the convergence of the LVPP iterates is provably at least sublinear \cite[Corollary A.12]{keith2023proximal}, and observed to converge linearly to strict minimizers \cite[Remark 4.18]{keith2023proximal}. Faster convergence is obtained by growing the sequence sufficiently quickly. In particular, picking the sequence $\alpha_k = r^{1/(q-1)} \mu^{q^k} - \alpha_k$ for $k>2$ where $\mu, q, r >1$, and $\alpha_1 = r^{1/(q-1)} \mu$, yields a worst-case superlinear convergence with order $q$ and rate $r$ \cite[Corollary A.12]{keith2023proximal}.

Provided the initial proximal parameter is sufficiently small, then many choices for an initial guess will converge in a handful of Newton iterations to the solution of the first LVPP subproblem. We usually opt for the zero function for $u$ and $\psi$ with some adjustments if other variables are included in the formulation. For each subsequent LVPP subproblem the Newton method is initialized with the solution to the previous LVPP subproblem. The downside to an aggressive update strategy is that each subsequent nonlinear LVPP subproblem becomes harder to solve. This manifests as requiring more Newton iterations per subproblem. In fact, when $\alpha_k \gg 1$, the proximal parameter may cause numerical ill-conditioning and a mesh-dependent number of cumulative Newton iterations, akin to what is observed in penalty methods.

Hence, one must balance the rate of convergence of the LVPP iterates to the true solution with minimizing the number of Newton iterations to solve each nonlinear LVPP subproblem; this tradeoff also depends on the difficulty of the problem. Deriving an optimal update strategy is an open problem. However, an extensive study was conducted for various update strategies in the context of the obstacle problem \cite[Figures 7 \& 9]{keith2023proximal}. In this work, we use varying strategies across the different examples. We found that a conservative but effective strategy is $\alpha_k = \min(c \alpha_{k-1}, C)$ for some positive constants $c, C >0$. The linear growth tends to result in small numbers of Newton iterations at each subproblem, and the cap of $\alpha_k$ at $C$ mitigates ill-conditioning in the discretized Newton linear systems.
A notable exception is in the variational fracture example of \Cref{sec:fracture} where the difficult nonlinearity of the problem leads us to an adaptive choice of proximal parameter depending on the number of Newton iterations required to solve the previous LVPP subproblem.
As a rule of thumb, if the Newton solver struggles to converge at a particular LVPP subproblem then consider a smaller starting value for $\alpha_1$, a slower growing sequence $\alpha_k$, and a smaller maximal value $C$.

\section{Pointwise bound constraints}
\label{sec:BoundConstraints}

We begin our study of LVPP with pointwise bound constraints.
This is the simplest and best-understood class of problems, allowing us to compare the performance of the derived methods to other standard approaches in the literature.

\subsection{Example 1: The obstacle problem}
\label{sub:ObstacleProblem}
We first consider the (unilateral) obstacle problem: minimize the Dirichlet energy
\begin{equation}
\label{eq:dirichlet-energy}
J(u) = \frac{1}{2} \int_\Omega \nabla u \cdot \nabla u \dd x - \int_\Omega f u \dd x
,
\end{equation}
where $f \in L^2(\Omega)$ is a prescribed body force,
over the feasible set
\begin{equation} \label{eq:polyhedral:unilateral_K}
K = \{ v \in H^1_0(\Omega) \mid v \ge \phi_1 \text{ a.e.~in } \Omega\},
\end{equation}
which is obtained from the general feasible set  \cref{eq:intro:feasible_set} by choosing $V = H^1_0(\Omega)$, $B = \operatorname{id}$, $\Omega_d = \Omega$, and $C(x) = [\phi_1(x), \infty)$. LVPP for this problem was previously analyzed in \cite{keith2023proximal,fu2024locallyconservativeproximalgalerkinmethod}. We also consider the bilateral obstacle problem with the feasible set
\begin{equation} \label{eq:polyhedral:bilateral_K}
K = \{ v \in H^1_0(\Omega) \mid \phi_1 \le v \le \phi_2 \text{ a.e.~in } \Omega\},
\end{equation}
with $C(x) = [\phi_1(x), \phi_2(x)]$. Here, $\phi_1, \phi_2 \in H^1(\Omega) \cap L^\infty(\Omega)$ with nonpositive and nonnegative boundary values respectively.
We can choose as a Legendre function for the unilateral case the generalized Shannon entropy,
\begin{equation}
\label{eq:GeneralizedShannonEntropy}
R(a) = (a - \phi_1) \ln (a - \phi_1) - (a - \phi_1), \text{ with } \nabla R^*(a^\ast) = \phi_1 + \exp{a^\ast},
\end{equation}
and for the bilateral case, the generalized Fermi--Dirac entropy,
\begin{align}
\label{eq:GeneralizedFermiDiracEntropy}
\begin{split}
R(a) &= (a - \phi_1) \ln (a - \phi_1) + (\phi_2-a) \ln (\phi_2-a),\\
\text{ with } \nabla R^*(a^\ast) &= \dfrac{\phi_1 + \phi_2\exp a^\ast}{1 + \exp a^\ast}.
\end{split}
\end{align}
The resulting saddle point formulation (in weak form) is:
for $\psi^0 = 0$, find $(u^{k}, \psi^{k}) \in H^1_0(\Omega) \times L^\infty(\Omega)$ satisfying
\begin{subequations} \label{eq:polyhedral:obstacle}
\begin{align}
\alpha_{k} (\nabla u^k, \nabla v) + (\psi^k, v) &= \alpha_k (f, v) + (\psi^{k-1}, v), \\
(u^{k}, w) - (\nabla R^*(\psi^{k}), w) &= 0,
\end{align}
\end{subequations}
for all $(v, w) \in H^1_0(\Omega) \times L^\infty(\Omega)$. Here and throughout, $(\cdot, \cdot)$ denotes the $L^2(\Omega)$ inner product. Each nonlinear PDE in this sequence is solved with Newton's method. Discretizing this problem requires choosing approximation spaces for $H^1_0(\Omega)$ and $L^\infty(\Omega)$, which is a standard procedure that is much more straightforward than directly discretizing the feasible sets \cref{eq:polyhedral:unilateral_K} or \cref{eq:polyhedral:bilateral_K}. In the Galerkin finite element method, the pair of approximating spaces must satisfy a compatibility condition (cf.~\cite[\S 4.7]{keith2023proximal}), but many standard choices of finite elements work, including, among others, equal-order continuous Lagrange elements for both $u$ and $\psi$.
As in~\Cref{sub:PG}, we refer to these Galerkin finite element discretizations as proximal Galerkin methods.

Proximal Galerkin methods always deliver two solutions, $u_h$ and $\tilde{u}_h \coloneqq \nabla R^*(\psi_h)$, with the latter bound-preserving by construction. Moreover, they typically exhibit mesh-independent convergence and can converge superlinearly with the number of outer iterations $k$ if the step sizes ${\alpha_k}$ are chosen appropriately (cf.~\cref{eq:intro:striking_convergence}). The PDE subproblems \cref{eq:polyhedral:obstacle} are straightforward to implement in standard finite element software and can be discretized with arbitrary-degree polynomial bases, if appropriate.

\begin{figure}
\centering
\begin{subfigure}{\textwidth}
\centering
\footnotesize
\renewcommand{\arraystretch}{1.3}
\begin{tabular}{ |l| >{\centering}p{0.5cm}  >{\centering}p{0.5cm}  >{\centering}p{0.5cm} | >{\centering}p{0.5cm}  >{\centering}p{0.5cm}  >{\centering\arraybackslash}p{0.5cm} | }
   \hhline{~|------}
 \multicolumn{1}{c|}{} & \multicolumn{3}{c|}{\cellcolor{lightgray!02} Degree $p = 1$} & \multicolumn{3}{c|}{\cellcolor{lightgray!02} Degree $p = 2$}\\
 \hline
 \rowcolor{lightgray!10}
 Method & $h$ & $h/2$ & $h/4$ & $h$ & $h/2$ & $h/4$ \\
 \hline
\cellcolor{lightgray!02}   Proximal Galerkin & 15 & 13 & 12 & 15 & 16 & 12  \\
\hhline{|~|~|~|~|---|}
\cellcolor{lightgray!02}   Active Set \cite{benson2006,munson2001} & 11 & 16 & 25 & \multicolumn{3}{c|}{\cellcolor{lightgray!10}}\\
\cellcolor{lightgray!02}   Trust-Region \cite{gould2003}  & 6 & 12 & 19 & \multicolumn{3}{c|}{\cellcolor{lightgray!10}}\\
\cellcolor{lightgray!02}   Interior Point (IP) \cite{wachter2006}  & 9 & 9 & 8 & \multicolumn{3}{c|}{\cellcolor{lightgray!10}}\\
\cellcolor{lightgray!02}   IP without Hessian \cite{wachter2006}  & 90 & 260 & 500 & \multicolumn{3}{c|}{\cellcolor{lightgray!10}\multirow{-4}{*}{\parbox{2cm}{\centering Not bound preserving}}} \\
 \hline
\end{tabular}
\centering
\caption{Number of linear system solves for popular solvers using various mesh sizes $h$.}
\label{fig:ObstacleExperiment_MeshIndependence}
\end{subfigure} \\[1.2em]
\begin{subfigure}[c]{0.3\textwidth}
\centering
\raisebox{-0.5\height}{\includegraphics[clip, trim = 5cm 7cm 4cm 6cm, width=0.7\linewidth]{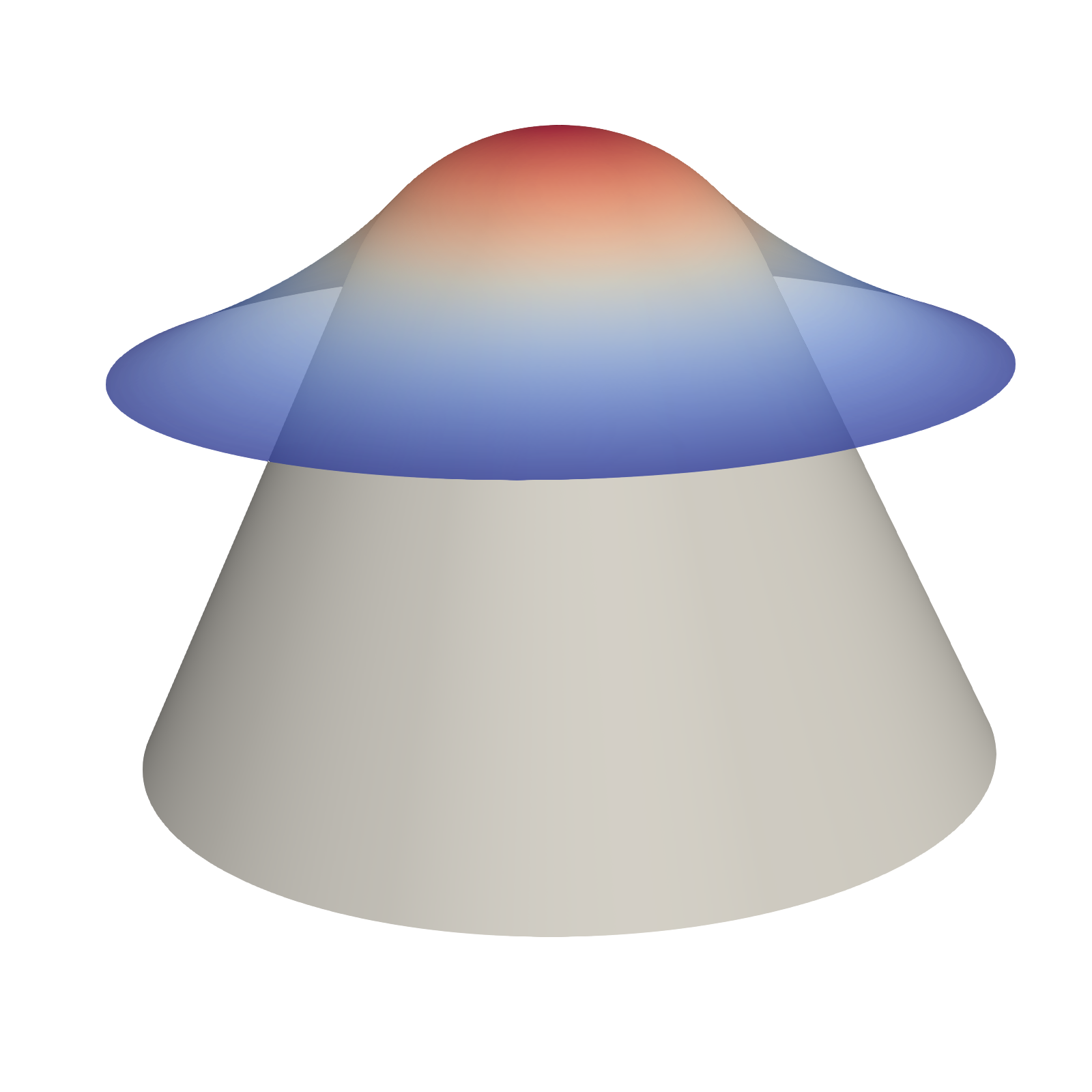}}
~
\raisebox{-0.5\height}{\includegraphics[width=0.09\linewidth]{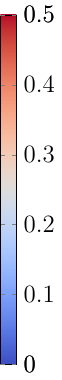}}
\caption{Obstacle $\phi$ (grey) and membrane $u$ (red/blue).}
\label{fig:ObstacleExperiment_Solution}
\end{subfigure}
~
\begin{subfigure}[c]{0.64\textwidth}
\centering
\footnotesize
\renewcommand{\arraystretch}{1.3}
\begin{tabular}{ |c| >{\centering}p{0.5cm}  >{\centering}p{0.5cm} >{\centering}p{0.5cm}  >{\centering}p{0.5cm}  >{\centering}p{0.5cm} >{\centering\arraybackslash}p{0.5cm}| }
 \hline
 \rowcolor{lightgray!10}
Mesh size $h$ & $2^{-1}$ & $2^{-2}$ & $2^{-3}$ & $2^{-4}$ & $2^{-5}$ & $2^{-6}$ \\
 \hline
\cellcolor{lightgray!02}  Finite Difference & 10 & 15 & 13 & 15 & 16 & 16  \\
\hline\hline
 \rowcolor{lightgray!10}
Degree $p$ & $8$ & $16$ & $24$ & $32$ & $40$ & $48$ \\
 \hline
\cellcolor{lightgray!02}  Spectral Method & 16 & 17 & 16 & 16 & 16 & 15 \\
 \hline
\end{tabular}
\caption{Number of linear system solves for the proximal finite difference and spectral methods.}
\label{fig:ObstacleExperiment_Discretizations}
\end{subfigure}
\caption{Example 1. The unilateral obstacle problem with the setup \cref{eq:ObstacleBenchmark}.}
\label{fig:ObstacleExperiment}
\end{figure}

We now compare proximal Galerkin with other popular solvers. We consider the benchmark unilateral obstacle problem from \cite[Experiment~4]{keith2023proximal}, where
\begin{align}
\begin{split}
\Omega &= \{ (x,y) \in \mathbb{R}^2 : 0 < r < 1\}, \quad r^2 = x^2+y^2, \quad f \equiv 0,\\
\phi(x,y) &= 
\begin{cases}
\sqrt{1/4-r^2} &  r \leq b, \\
d+b^2/d - b r/d & r > b,
\end{cases} \quad b=9/20, \quad d = \sqrt{1/4-b^2}.
\end{split} \label{eq:ObstacleBenchmark}
\end{align}
We plot the obstacle and solution in~\Cref{fig:ObstacleExperiment_Solution}. We use equal-order continuous Lagrange elements to discretize the solution $u$ and the latent variable $\psi$. The proximal parameter $\alpha_k$ is updated with a heuristic double-exponential rule; i.e.,~for all $ k \in \mathbb{N}$, we set
\begin{align}
\alpha_k = \min( \max(r^{q^k} - \alpha_{k-1},1), 10^2), \; r=q=3/2,\; \text{and}\; \alpha_{0} = 1.
\label{eq:alpha-double-exponential}
\end{align}
Note that we do not require the proximal parameter $\alpha_k$ to pass to infinity with $k$.
In \Cref{fig:ObstacleExperiment_MeshIndependence}, we report the number of linear system solves required by the proximal Galerkin method as well as an active set strategy \cite[Alg.~3.1]{benson2006} found in the popular scientific computing toolkit PETSc \cite{petsc-user-ref}\footnote{The active set method \cite[Alg.~3.1]{benson2006} can be interpreted as a semismooth Newton method for the discretized obstacle problem. However, the optimality conditions are, in general, insufficiently regular to define its infinite-dimensional formulation; see \Cref{sec:intro}.}, IPOPT \cite{wachter2006}, an interior point method, with and without Hessian access, and the bound-constrained trust-region method (TRB) as implemented in GALAHAD \cite{gould2003}.
The practical benefits of proximal Galerkin are immediately clear; the other methods, except IPOPT with Hessian access, exhibit mesh dependence.
While the interior point method performs the fewest linear solves, proximal Galerkin is the only solver that also preserves bounds for higher-order discretizations.
Note that by providing the Hessian matrix to IPOPT, we are dimensionally-rescaling the derivatives of the energy functional, which is known to often restore mesh independence in solvers for PDE-constrained optimization problems \cite{schwedes2017mesh}.

The LVPP saddle point subproblem \cref{eq:polyhedral:obstacle} can be discretized with many other techniques.
We also provide results where the subproblem \cref{eq:polyhedral:obstacle} is discretized with a coefficient-based Zernike sparse spectral method  \cite{burns2020,Papadopoulos2024,VasilDisk,Olver2013} and a five-point stencil finite difference method. For the finite difference scheme, we change the domain to the square $\Omega = (-1,1)^2$.  Here, we again use the double-exponential update rule \cref{eq:alpha-double-exponential} for $\alpha_k$. We terminate once $\| {\bf u}_k - {\bf u}_{k-1} \|_{\ell^2} < 10^{-9}$ where ${\bf u}_k$ is the discrete coefficient vector for $u$ at iteration $k$. The results are provided in \Cref{fig:ObstacleExperiment_Discretizations} where we observe $h$- and $p$-independent iteration counts for the proximal finite difference and spectral methods, respectively.
Further numerical experiments with the obstacle problem can be found in \cite{keith2023proximal,fu2024locallyconservativeproximalgalerkinmethod,papadopoulos2024hierarchicalproximalgalerkinfast}.

\subsection{Example 2: The Signorini problem}\label{ssec:sig}

We now consider the classical Signorini problem.
This problem demonstrates for the first time an extension of LVPP to pointwise bound constraints acting solely on the boundary of a computational domain $\Omega \subset \mathbb{R}^3$.
In this problem, we separate the boundary $\partial \Omega = \overline{ \Gamma_\mathrm{D} \cup \Gamma_\mathrm{T} }$ into disjoint measurable subsets for imposing displacement and traction boundary conditions.

The Signorini problem, posed by Signorini in 1959 \cite{signorini1959} and analyzed by Fichera in 1963 \cite{fichera1963}, is the essential first problem in contact mechanics. It models the deformation of a linear elastic body in the presence of a contact boundary constraint. The problem is posed on
\begin{equation}
\label{eq:Signorini_V}
    V = \big\{ u \in H^1(\Omega,\mathbb{R}^3)
        \mid u = g \text{ on } \Gamma_\mathrm{D\textbf{}}
    \big\}
    \,,
\end{equation}
and involves the minimization of the strain energy function
\begin{equation}
    J(u) =
    \frac{1}{2}\int_\Omega (\C \epsilon(u)) : \epsilon(u) \dd x
    -
    \int_\Omega f \cdot u \dd x
    \,,
\end{equation}
over the feasible set
\begin{equation}
    K = \big\{ u \in V
        \mid u \cdot \tilde{n} \leq \phi_1 \text{ on } \Gamma_\mathrm{T}
    \big\}
    \,.
\end{equation}
Here, $\epsilon : H^1(\Omega,\mathbb{R}^3) \to L^2(\Omega,\mathbb{R}^{3\times 3}_{\mathrm{sym}})$, $\epsilon \coloneqq (\nabla + \nabla^\top)/2$ denotes the symmetric gradient, $\C \colon \mathbb{R}^{3\times 3}_{\mathrm{sym}} \to \mathbb{R}^{3\times3}_{\mathrm{sym}}$ denotes the symmetric positive-definite elasticity tensor, $f \colon \Omega \to \mathbb{R}^3$ is an internal body force density, $\phi_1 \colon \Gamma_\mathrm{T} \to \mathbb{R}_+$ is a prescribed gap function, and $\tilde{n} \colon \Gamma_\mathrm{T} \to \mathbb{R}^3$ is a prescribed vector field. 
For simplicity of presentation, we assume that the displacement boundary conditions are homogeneous ($g = 0$) in the formulation below.

Notice that $K$ is obtained from the general feasible set  \cref{eq:intro:feasible_set} by choosing $V$ as in~\cref{eq:Signorini_V}, $B = \gamma(\cdot) \cdot \tilde{n}$, $\Omega_d = \Gamma_\mathrm{T}$, and $C(x) = (-\infty,\phi_1(x)]$.
Applying LVPP with the Legendre function~\cref{eq:GeneralizedShannonEntropy}, the resulting saddle-point formulation~\cref{eq:intro:lvpp} is:
for $\psi^0 = 0$, find $(u^{k}, \psi^{k}) \in V \times L^\infty(\Gamma_\mathrm{T})$ satisfying
\begin{subequations}
\label{eq:SignoriniVF}
\begin{align}
    ( \alpha_k\C \epsilon(u^k), \epsilon(v) ) - (\psi^k, v\cdot \tilde{n} )_{\Gamma_\mathrm{T}}
    &=
    (\alpha_k f, v) - ( \psi^{k-1}, v \cdot \tilde{n} )_{\Gamma_\mathrm{T}}
    \,,
    \\
    ( u^k\cdot \tilde{n} , w )_{\Gamma_\mathrm{T}} + ( \exp \psi^k, w )_{\Gamma_\mathrm{T}}
    &= ( \phi_1, w )_{\Gamma_\mathrm{T}}
    \,,
\end{align}
\end{subequations}
for all $(v, w) \in V \times L^\infty(\Gamma_\mathrm{T})$, where $(\cdot, \cdot  )_{\Gamma_\mathrm{T}}$ denotes the $L^2(\Gamma_\mathrm{T})$-inner product.

\begin{figure}
    \centering
    \begin{minipage}[c]{0.5\linewidth}
        \includegraphics[width=\linewidth]{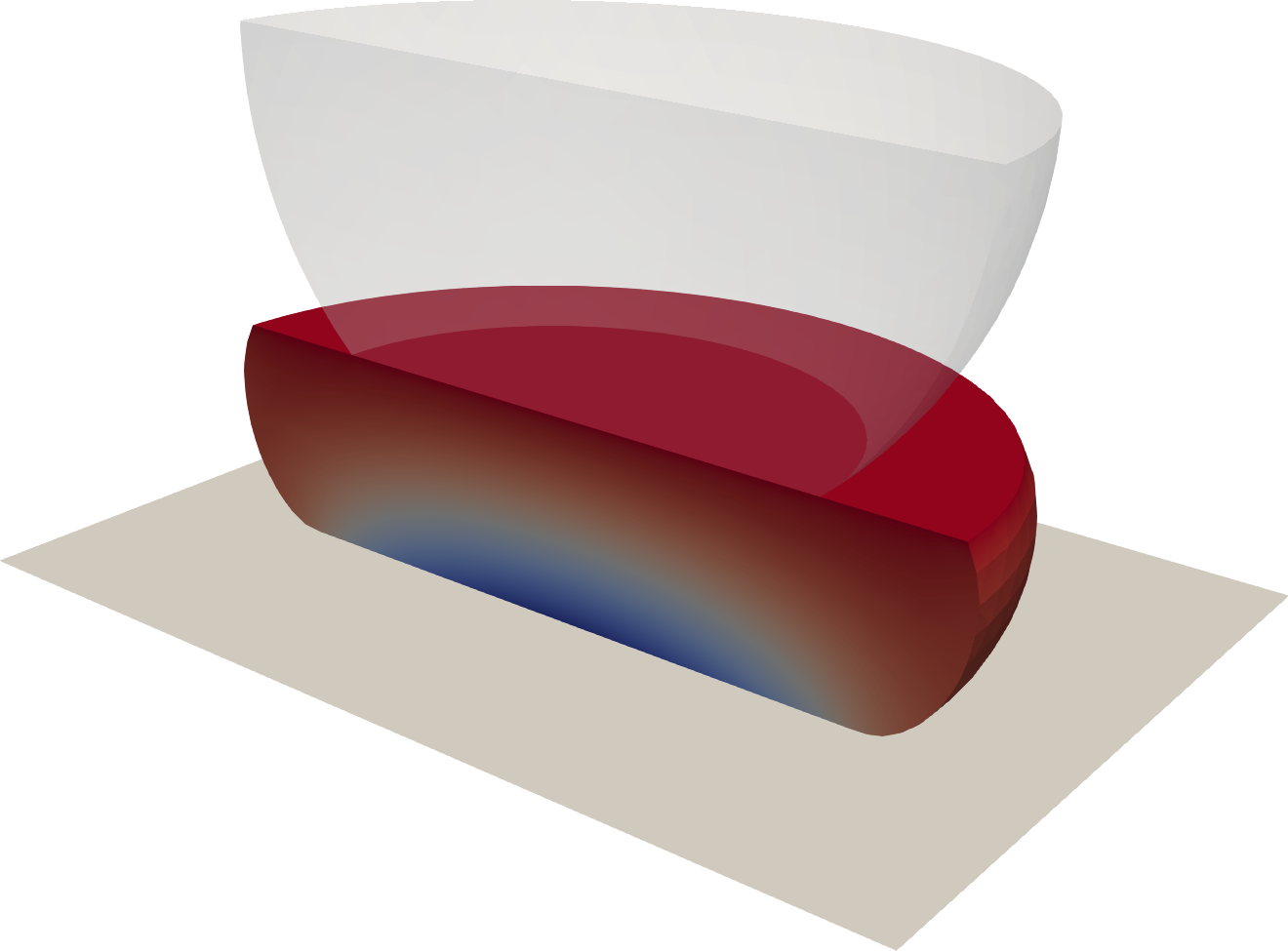}
    \end{minipage}
    ~
    \begin{minipage}[c]{0.05\linewidth}
        \includegraphics[width=\linewidth]{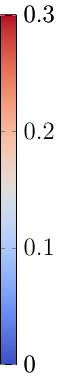}
    \end{minipage}
    \caption{Example 2.
    Final solution of the Signorini problem from \Cref{ssec:sig}. The undeformed configuration is shown in transparency. Both have been clipped vertically to emphasize the internal displacement (colored).}
    \label{fig:contact}
\end{figure}

As for the obstacle problem in \Cref{sub:ObstacleProblem}, we define a proximal Galerkin method by using equal-order continuous Lagrange spaces for the displacement and latent variable. Note that the spaces arising in~\cref{eq:SignoriniVF} are defined on manifolds of differing dimensions.
This is inherited in the discretization, and hence, the two discrete subspaces are not the same.
We use the mixed-dimensional assembly routines in DOLFINx \cite{baratta2023,dean2024} to solve the coupled problem.
The discrete problem is solved for half of a sphere centered at $(0,0,0.5)$ with radius $0.4$ coming into contact with the rigid $x,y$-plane. Here $\tilde{n}=(0,0,-1)^\top$ and $\phi_1$ measures the vertical distance between the $x,y$-plane and the boundary of the half sphere, i.e., $\phi_1(x,y,z) = z$. The displacement of the object and the rigid plane is visualized in~\cref{fig:contact}. Second-order Lagrange elements approximate the geometry, the displacement, and the latent variable.
The initial proximity parameter $\alpha_0 = 0.005$ is doubled at each subsequent proximal step (i.e., LVPP iteration). The tolerance for Newton’s method is set to $10^{-6}$. The LVPP algorithm requires a total of 13 proximal iterations, amounting to 17 linear solves. The maximum number of Newton steps was 5, the minimum 1. The number of degrees of freedom was $17,613$ for a mesh with $3,661$ cells.

\subsection{Example 3: Variational fracture} \label{sec:fracture}

The variational theory of brittle fracture of Francfort \& Marigo \cite{Francfort1998} is formulated in terms of a function $u$ over a domain $\Omega$ (representing the displacement) and an unknown subset $\Gamma \subset \Omega$ (representing the crack). The phase-field approach to fracture pioneered by Bourdin, Francfort \& Marigo \cite{Bourdin2000} approximates the crack set $\Gamma$ with an order parameter $c \colon \Omega \to [0,1]$, where $0$ represents the undamaged state and $1$ represents the complete presence of a crack, with the localization of the crack controlled by a length scale $\ell > 0$. 

The first major challenge with solving the resulting system is that it arises from a non-convex optimization problem. In fact, several fracture problems are known to possess multiple, physically-relevant solutions~\cite{Gerasimov2020}. However, the objective function is biconvex, in that the problem in $u$ or $c$ is convex if the other variable is fixed. This inspires the most popular algorithm for its solution, alternating minimization~\cite{Bourdin2000}, which iterates between minimizing only over displacement and damage, respectively. While convenient, alternating minimization can be difficult to apply to more complex physical models that lack this biconvex structure.

The second major challenge is that the numerical solution of phase-field fracture problems demands very high resolution. Experience shows that the mesh size $h$ for the discretization must be on the order of the crack localization length scale $\ell$, which is typically much smaller than the diameter of $\Omega$ for realistic simulations. These resolution requirements place a premium on mesh-independent algorithms.

We consider the anti-plane shear test formulated by Burke et al.~\cite{Burke2010}, which exhibits multiple solutions~\cite{Gerasimov2020}. This problem is formulated only in terms of the vertical displacement of the body, as opposed to the full displacement typically used for other problems, but this is not central. The resulting variational inequality is closely related to the Ambrosio--Tortorelli approximation of the Mumford--Shah functional \cite{Mumford1989,Ambrosio1990}.

The crack evolution is computed over a quasi-static incremental loading procedure. At each step, the load on the body is varied, and the body is allowed to equilibrate. Mathematically, at each loading step, we are given Dirichlet boundary data $g$ and the previous order parameter $c_{\rm prev}$. This defines our feasible set $K$:
\begin{equation}
\label{eq:FeasibleSetFracture}
K = \left\{ (u,c) \in H^1_g(\Omega) \times H^1(\Omega) \mid 0 \le c_{\rm prev} \le c \le 1 \; \text{a.e.~in} \; \Omega \right\},
\end{equation}
where
\begin{equation}
H^1_g(\Omega) = \left\{ u \in H^1(\Omega) \mid \gamma u = g \text{ on } \Gamma_D \right\},
\end{equation}
with $\Gamma_D \subset \partial \Omega$.
Given energy release rates $G, G_c > 0$, the length scale $\ell > 0$, and an artificial residual stiffness of a fully ruptured phase $\epsilon > 0$, we arrive at the (non-convex) objective functional:
\[
J(u,c) \coloneqq \frac{G}{2}\int_{\Omega}(\epsilon + (1-\epsilon)
(1-c)^2)|\nabla u|^2\,\mathrm{d}x 
+ 
\frac{G_c}{2}\int_{\Omega} \ell |\nabla c|^2 
+ 
\ell^{-1}|c|^2\,\mathrm{d}x.
\]

To obtain the feasible set $K$~\cref{eq:FeasibleSetFracture} from~\cref{eq:intro:feasible_set}, we select $V = H^1_g(\Omega) \times H^1(\Omega)$,\footnote{This is an affine Banach space, not a Banach space, but the arguments of \Cref{sec:intro} carry over, mutatis mutandis.} $B = (0,\operatorname{id})$, with $0$ denoting the zero operator, $\Omega_d = \Omega$, and $C(x) = \mathbb R \times [c_{\rm prev}(x),1]$.
Since there is no inequality constraint on the displacement variable $u$, we introduce only a single latent variable $\psi$ corresponding to the damage $c$.
We then use the Fermi--Dirac entropy~\cref{eq:GeneralizedFermiDiracEntropy} with $\phi_1 = c_{\rm prev}$ and $\phi_2 = 1$, giving us the following LVPP saddle point problem: for $\psi^0 = 0$, find $(u^k, c^k, \psi^k) \in H^1_g(\Omega) \times H^1(\Omega) \times L^\infty(\Omega)$ such that

\begin{subequations}
\label{eq:FractureVF}
\begin{align}
    \alpha_k G ((\epsilon + (1-\epsilon)(1-c^k)^2)\nabla u^k, \nabla v)
    &= 0
    \,,
    \\
    \begin{split}
    -\alpha_k G ((1-\epsilon)(1-c^k)|\nabla u^k|^2,d)\hspace{30mm}\\
    \indent +
    \alpha_k G_c ( \ell(\nabla c^k, \nabla d) + \ell^{-1}(c^k, d) )
    +(\psi^k,d)
    &=
    (\psi^{k-1},d)
    \,,
    \end{split}
    \\
    (c^k,w)
    -
    \bigg(\frac{c_\mathrm{prev} + \exp (\psi^k)}{\exp(\psi^k) + 1}, w\bigg)
    &=
    0
    \,,
\end{align}
\end{subequations}
for all $(v, d, w) \in H^1_0(\Omega) \times H^1(\Omega) \times L^\infty(\Omega)$.

\begin{figure}[t]
\centering
\begin{minipage}[c]{4cm}
\begin{tikzpicture}
    \node[anchor=south west,inner sep=0] (image) at (0,0) {\includegraphics[width=\textwidth]{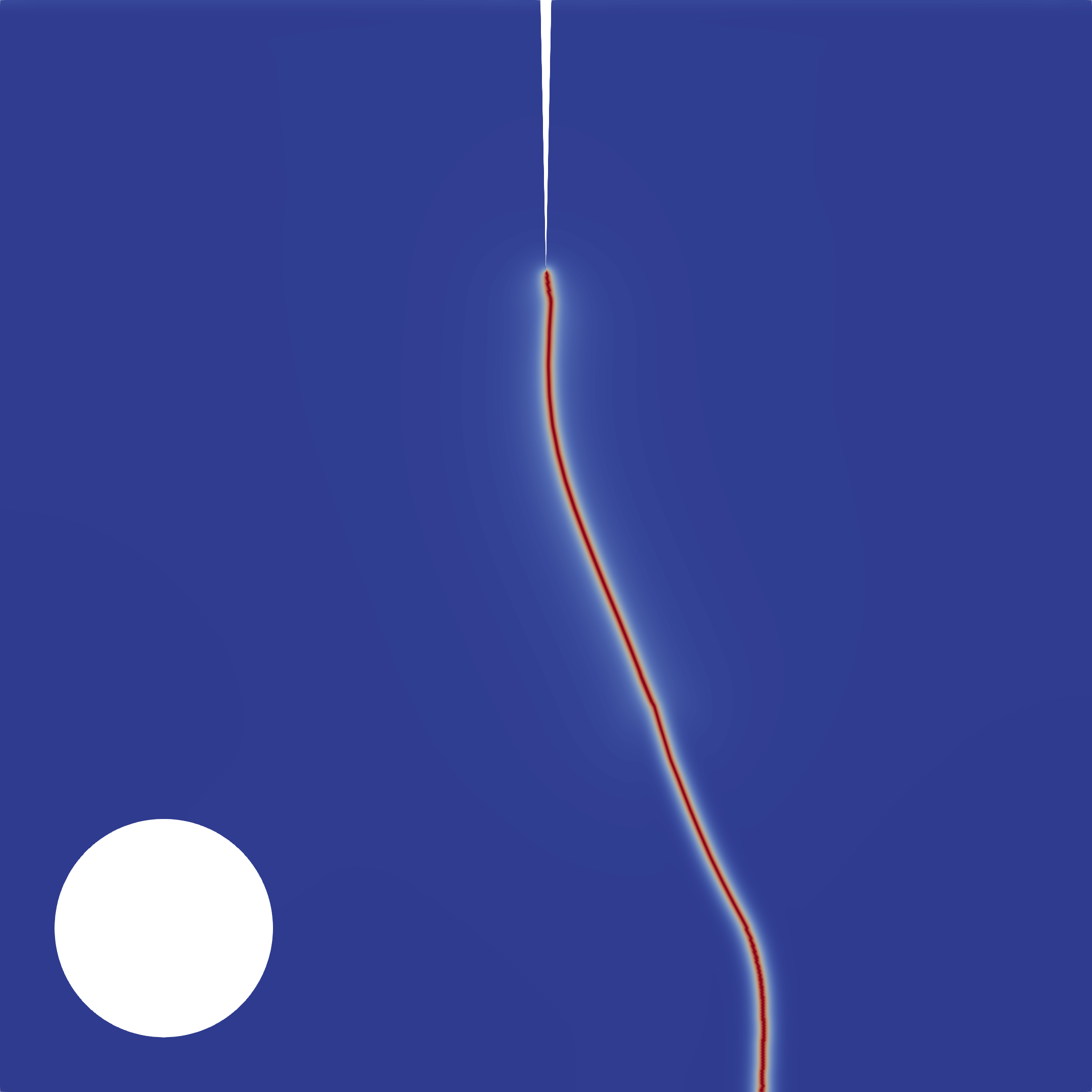}};
    \begin{scope}[x={(image.south east)},y={(image.north west)}]
        \node (bcleft) at (0.25, 1.05) {\small $u = -L$};
        \node (bcright) at (0.75, 1.05) {\small $u = +L$};
    \end{scope}
\end{tikzpicture}
    \end{minipage}
\hspace*{0.25cm}
 \begin{minipage}[c]{0.9cm}
     \includegraphics[width=\linewidth]{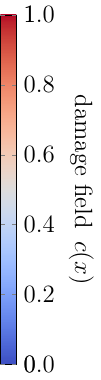}
 \end{minipage}
\caption{Example 3. The final damage field for the fracture problem described in \Cref{sec:fracture}.}
\label{fig:fracture}
\end{figure}

The anti-plane shear test imposes $u = +L$ on the top-right boundary (to the right of the initial notch) and $u = -L$ on the top-left (to the left of the notch), see \Cref{fig:fracture}. At each loading step, the imposed displacement $L$ is incremented by $0.005$, starting from zero. Natural boundary conditions are imposed on the remainder of the boundary.

We define a proximal Galerkin method using continuous piecewise linear finite elements for $u$, $c$, and $\psi$ and solve the resulting discrete systems with Newton's method. For robustness of the optimization procedure, we perturb the Jacobian of \eqref{eq:FractureVF} with a coercive modification, by adding
$
\iota ((u, v) + (c, d) - (\psi, w))
$
to the Jacobian with perturbation parameter $\iota = 10^{-3}$.
The proximity parameter $\alpha_k$ is updated using a simple heuristic based on the number of Newton iterations needed to solve the previous subproblem: if the number of Newton iterations is four or fewer, $\alpha_{k+1} = 2\alpha_k$; if the number of Newton iterations is ten or more, $\alpha_{k+1} = \alpha_k/2$; otherwise $\alpha_{k+1} = \alpha_k$. On average, we require 2.85 proximal steps (i.e.~LVPP iteration) per loading step, although with a large variance; the worst loading step took 147 proximal steps, with the second worst taking 10. Each proximal step required, on average, 5.44 Newton steps. The final damage field at $L = 2$ is depicted in \Cref{fig:fracture}.

\subsection{Example 4: Vector-valued multiphase problems}
\label{sec:cahn-hilliard}
Many challenging physical problems include multiple phases, materials, or components that are required to be balanced throughout the domain. For example, reservoir simulation in petroleum engineering  \cite[Chap.\ 3]{Fanchi2005-lp},  distillation columns in chemical engineering \cite[Chap.\ 11]{biegler2010nonlinear}, and optimal design in semiconductors \cite{Adam2019}, all explicitly require such constraints. From a mathematical perspective, the state variables, whether they be materials, components, or concentrations in the domain $\Omega$ are represented by vector-valued functions $u \colon \Omega \to \mathbb R^m$, where $m$ denotes the number of states.
The natural feasible set for such systems is a Gibbs simplex over a specified function space $V$; namely,
\begin{equation}
\label{eq:Simplex}
    K(V) := \left\{ 
    u \in V
    \left| \; 
    \sum_{i=1}^m u_i(x) = 1 \text{ and } 
    u_i(x) \ge 0\; 
    (i=1,\dots,m)
    \text{ f.a.e.\ } x\in\Omega
    \right.
    \right\}.
\end{equation}
This feasible set $K = K(V)$ fits into our framework by setting $B = \operatorname{id}$, $\Omega_d = \Omega$, and $C = \{a \in \mathbb{R}^m \mid \sum_{i=1}^m a_i = 1,\; a_i \geq 0\}$.
In this case, the LVPP subproblem has a relation $u = \nabla R^\ast(\psi)$ induced by the Legendre function
\begin{equation}
    R(a) = \sum_{i=1}^m a_i \ln(a_i)
    , ~\text{ with }~
    \frac{\partial R^*}{\partial a^\ast_i}(a^\ast) = \frac{\exp(a^\ast_i)}{\sum_{j=1}^m \exp(a^\ast_j)}
    ~\text{ for each } i=1,\ldots,m,
    \end{equation}
which couples the components of the primal variable $u : \Omega \to C \subset \mathbb{R}^m$ (the states) to the vector-valued latent variable $\psi : \Omega \to [-\infty,\infty]^m$.

For our numerical example, inspired by \cite{Garcke1998}, we construct a Cahn--Hilliard-type gradient flow toward a constrained minimizer of the Ginzburg--Landau-type energy functional
\[
    I(u)
    =
    \int_\Omega
        \frac{\varepsilon^2}{2} \sum_{i=1}^m|\nabla u_i|^2 + W(u)
    \dd x
    \,,
    \qquad
    W(u)
    =
    \sum_{i=1}^m
        u_i(1-u_i)
    +
    i_K(u)
    \,,
\]
where $i_K(u) = 0$ if $u \in K$ and $i_K(u) = +\infty$ otherwise denotes the indicator function for $K$.
Here, $\varepsilon > 0$ is a diffuse interface parameter that controls the width of the region containing mixed states.
The seminal works of Modica \cite{Modica1987a,Modica1987b} demonstrate that functionals of this type tend (in the sense of $\Gamma$-convergence \cite{DalMaso1993}) to a functional measuring the perimeter of the phase boundaries as $\varepsilon \to 0$. Therefore, the Ginzburg--Landau energy acts as a penalty on pathological free boundaries (phase transitions) and favors straight edges. This model problem can be easily modified with more complex energy functionals that, for instance, reflect anisotropic behavior amongst the various components represented by $u$ \cite{garcke1999multiphase}.

Following \cite[Equation (1)]{fife2000models}, we consider an $H^{-1}$-gradient flow over $\Omega$ described by the system of differential inclusions
\[
    \frac{\dd u}{\dd t} 
    \in
    -\varepsilon^2 \Delta^2 u + \Delta \partial W(u)
        \,,
            \]
where $\partial$ denotes the Clarke subdifferential of the nonsmooth functional $W$ \cite[Chap.\ 2]{Clarke1990}.
Using backward Euler to discretize this system in time is equivalent to solving a recursive sequence of minimization problems over the feasible set~\cref{eq:Simplex} with $V = H^2(\Omega)$.
In particular, at each time step, we minimize
\begin{equation}\label{eq:gl-euler}
    J(u) = \frac12\int_{\Omega}|u - u_{\rm prev}|^2 \, \mathrm{d}x
    +
    \tau \int_{\Omega}
        \frac{\varepsilon^2}{2} \sum_{i=1}^m (\Delta u_i)^2
        +
        \Delta W(u)
    \dd x
    \,,
\end{equation}
where $\tau > 0$ is the width of the time step, and $u_{\rm prev}$ is the solution at the previous point in time. As throughout the text, we assume high enough regularity so that $\Delta W(u) \in L^1(\Omega)$.
We can readily apply LVPP to these semi-discrete subproblems.

As is often done in practice, see, e.g., \cite{Elliott1989}, we avoid discretizing the natural function space for $u$, namely $H^2(\Omega;\mathbb R^m)$, by introducing slack variables $z_i = \varepsilon^2\Delta u_i - [\partial W(u)]_i$.
Given the previous solution $u_{\rm prev} \in H^1(\Omega; \mathbb R^m)$, the LVPP subproblems take the form: for $\psi^0 = 0$, find $(u^k, z^k, \psi^k) \in H^1(\Omega;\mathbb R^m) \times H^1(\Omega;\mathbb R^m) \times L^\infty(\Omega;\mathbb R^m)$ such that for, $i=1,\dots,m$, we have
\begin{subequations}
\label{eq:allen-cahn}
\begin{align}
    \alpha_k(z^k_i,y_i) +
    \varepsilon^2 \alpha_k(\nabla u^k_i,\nabla y_i) 
    - 2\alpha_k(u^k_i,y_i)
    + (\psi^k_i,y_i) 
    &= 
    (\psi^{k-1}_i,y_i) - \alpha_k(1, y_i),    \\
    (u^k_i,v_i) -
    \tau (\nabla z^k_i,\nabla v_i) 
    &= 
        (u_{\rm prev,i},v_i),        \\
    (u^k_i,w_i) - \left(\frac{\exp(\psi^k_i)}{\sum_{j=1}^m \exp(\psi^k_j)},w_i\right) &= 0,\label{eq:allen-cahn-b} \end{align}
\end{subequations}
for all $(v, y, w) \in H^1(\Omega;\mathbb R^m) \times H^1(\Omega;\mathbb R^m) \times L^\infty(\Omega;\mathbb R^m)$. Note that LVPP generates two semi-discrete flows using the primal variables $u(t_1),u(t_2),\dots$ and the latent solutions $\widetilde{u}(t_1),\widetilde{u}(t_2),\dots$ for the time steps $0 \le t_1 < t_2 < \cdots$, where each $\widetilde{u}_i := \frac{\exp(\psi^k_i)}{\sum_{j=1}^m \exp(\psi^k_j)}$ arises from~\cref{eq:allen-cahn-b}. The latent flow satisfies the constraints in $K$ even after discretizing~\cref{eq:allen-cahn} in space.

In our numerical experiments, we solve~\cref{eq:allen-cahn} on a square domain in $\mathbb R^2$ with Neumann/natural boundary conditions imposed on $u$ and $z$.
We increment the differential inclusion with a time step of size $\tau = 10^{-5}$ and start from an initial array of $m = 4$ phases that sum to one throughout the domain, inspired by \cite[Sec. 4.3]{wu2017multiphase}.

\begin{figure}
    \centering
    \subcaptionbox*{$t=0$} {\includegraphics[width=0.23\linewidth]{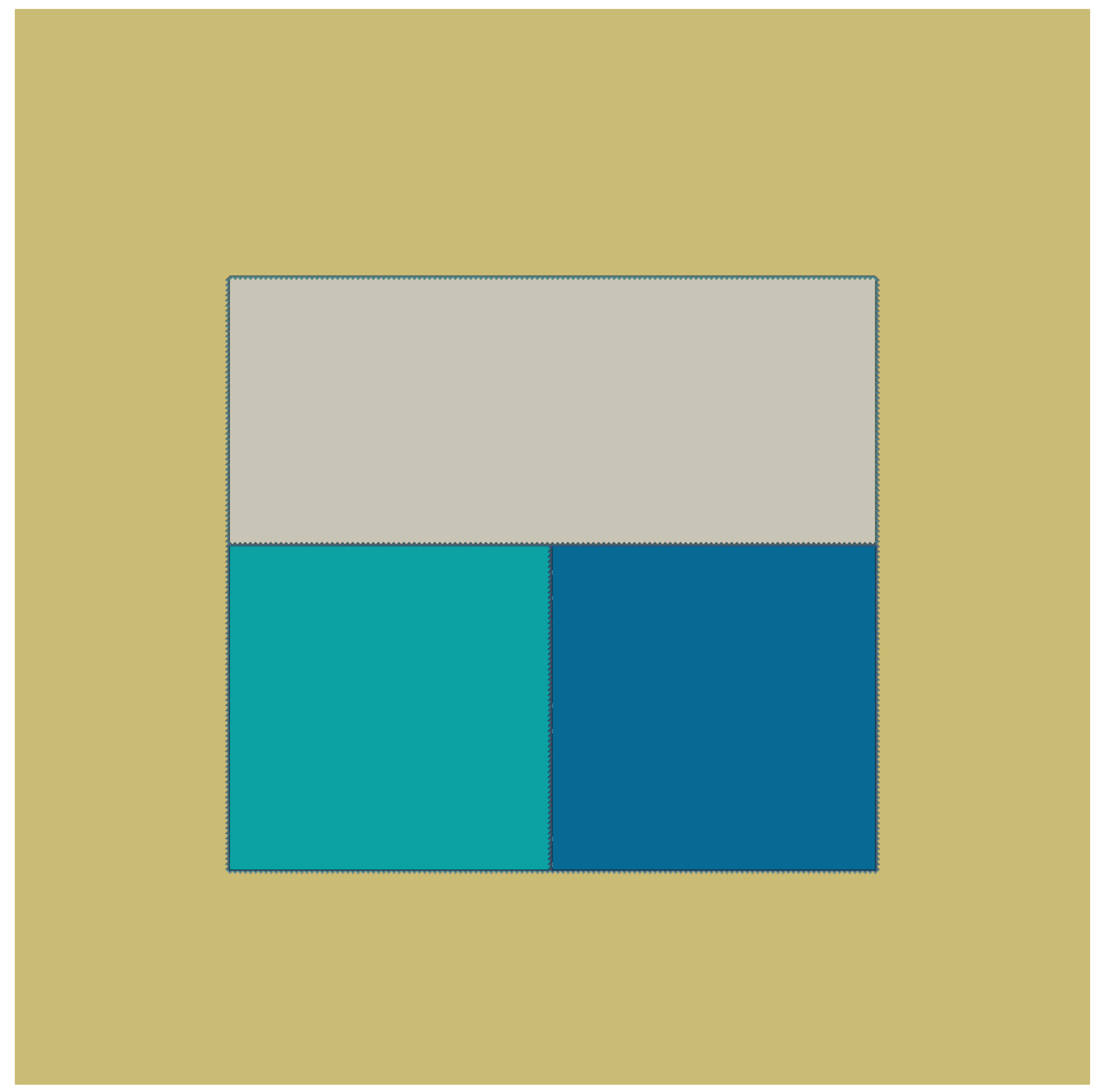}}
    \subcaptionbox*{$t=10^{-4}$} {\includegraphics[width=0.23\linewidth]{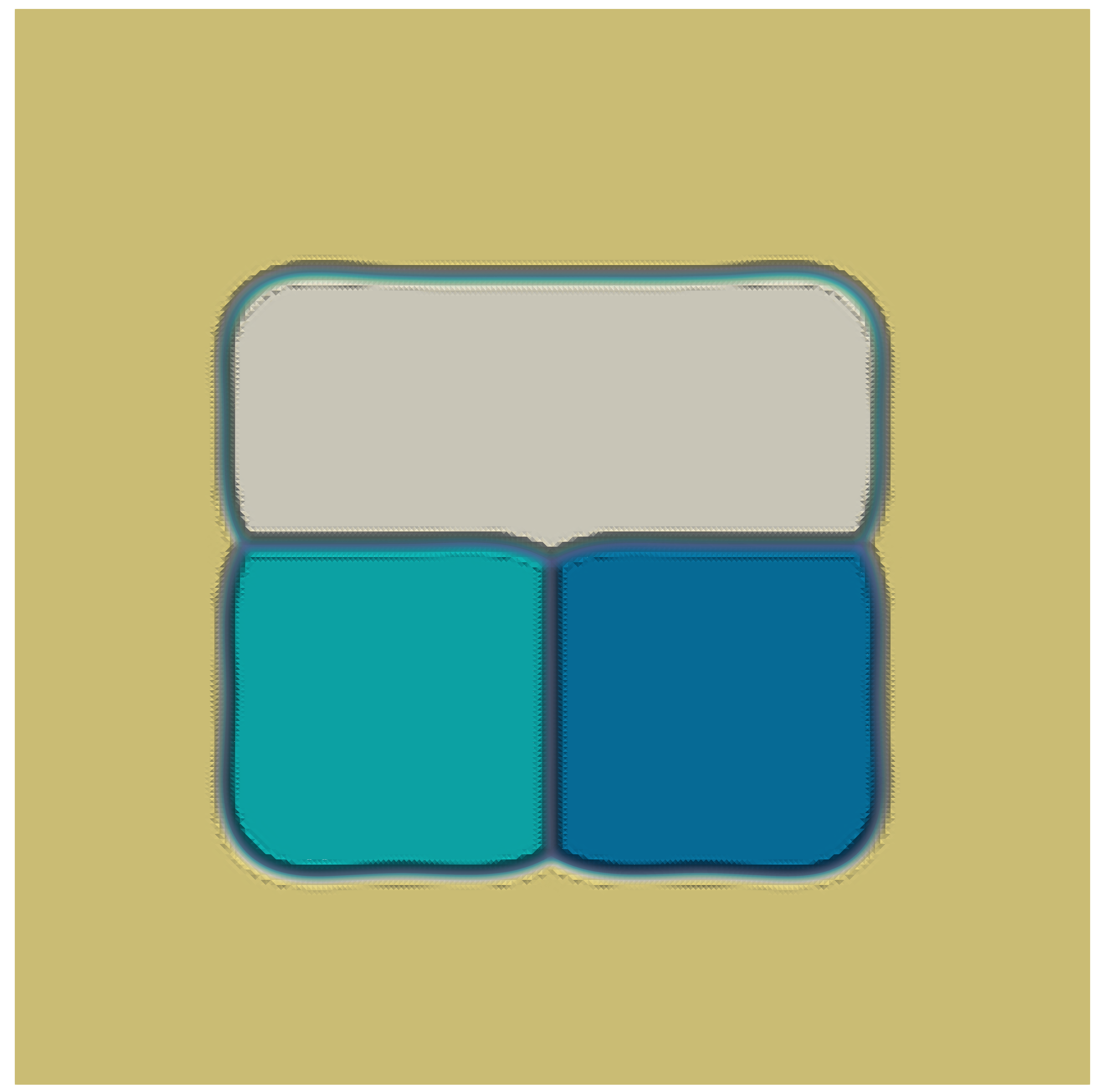}}
    \subcaptionbox*{$t=10^{-3}$} {\includegraphics[width=0.23\linewidth]{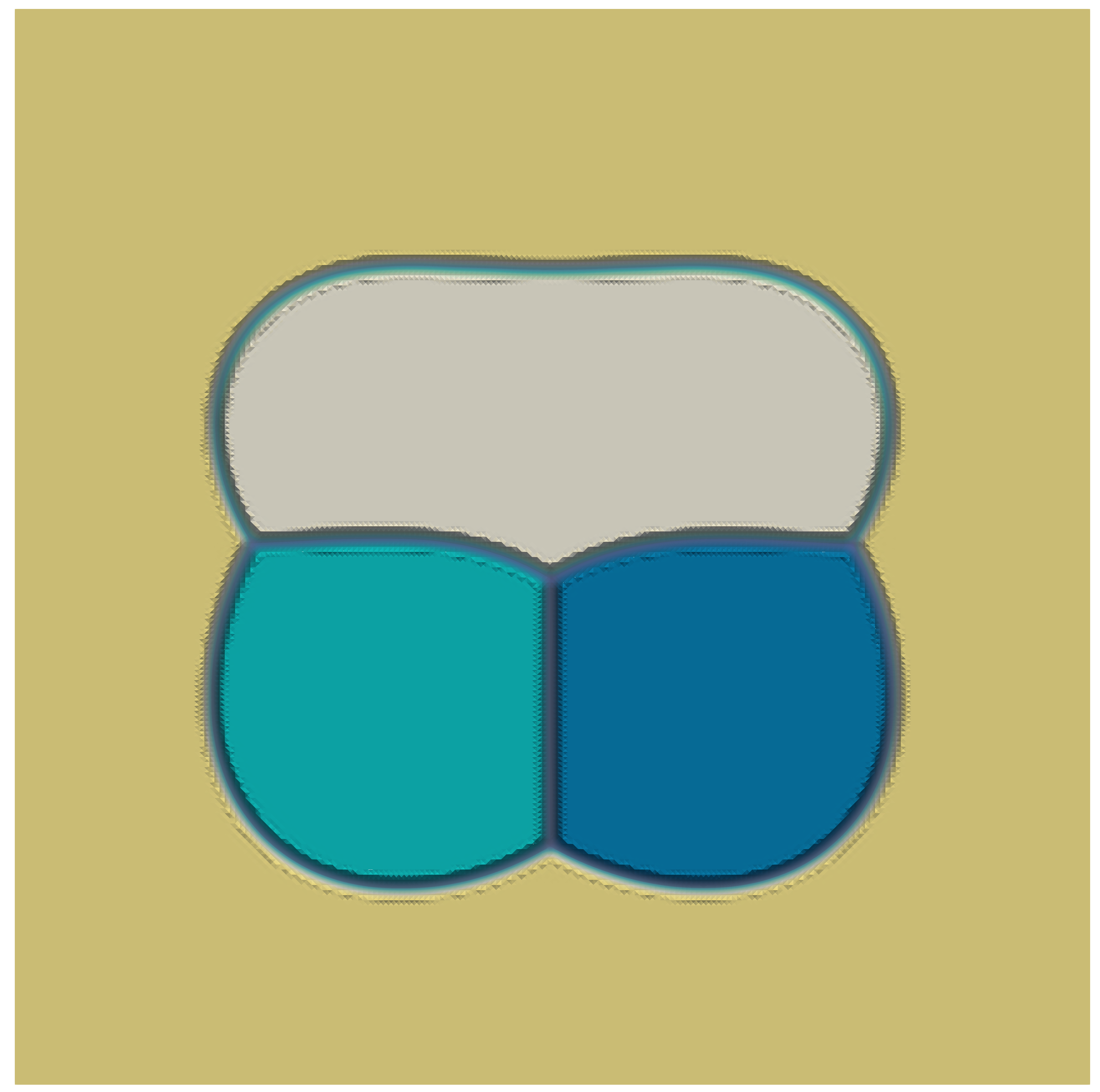}}
    \subcaptionbox*{$t=7\times 10^{-3}$} {\includegraphics[width=0.23\linewidth]{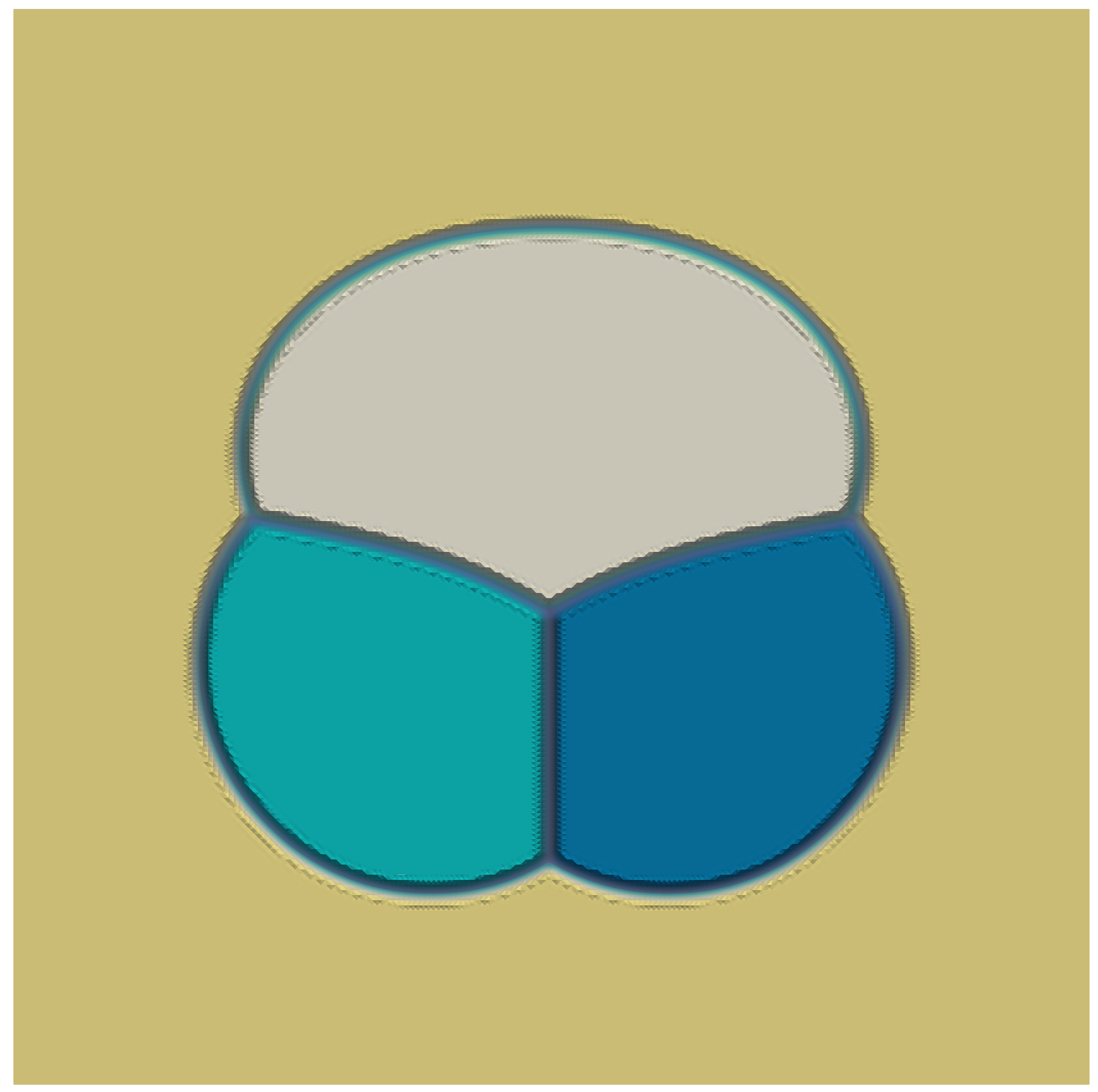}}
            \caption{Example 4. The evolution of the four phases $u_i$, $i=1,2,3$, of the Cahn--Hilliard problem described in \Cref{sec:cahn-hilliard} at times $t=0, 10^{-4}$, $10^{-3}$, and $7\times10^{-3}$.}
    \label{fig:cahn-hilliard}
\end{figure}

We discretize the underlying function spaces by piecewise linear continuous finite elements and solve the resulting discrete subproblems with Newton's method to a tolerance of $10^{-8}$. The proximity parameter $\alpha_k = 1$ is kept constant. On average, we require 8 proximal steps per time step, where the maximum number of proximal steps over all time intervals was 8 and the minimum 7. The average number of linear solves per Newton iteration was 2, resulting in roughly 16 linear solves per time step. We plot several snapshots of the evolution in
\Cref{fig:cahn-hilliard}.

\subsection{Example 5: Obstacle-type quasi-variational inequalities}
\label{sec:qvi}
Obstacle-type quasi-variational inequalities (QVIs) \cite{alphonse2019b} introduce an additional layer of complexity to the modeling framework. Here, the pointwise state constraints are no longer fixed but rather depend on the solution itself, which substantially contributes to the nonsmooth and nonlinear characteristics of the problem. The extremely versatile nature of QVIs makes them the model of choice for a number of important problems in physics and economics \cite{alphonse2019,aubin2007,baiocchi1984,bensoussan1975,bensoussan1984,lions1986,mosco1976}. However, this structure also means their numerical treatment poses a significant challenge, and the literature on solvers that tackle infinite-dimensional QVIs is limited. Fixed point methods are popular \cite{agarwal2003,alphonse2019b} but require that the QVI be reformulated as a fixed point problem with a contraction. This restricts its application to a small subclass of QVIs, and even if there is a contraction, the contraction coefficient may be close to one, requiring a computationally infeasible number of iterations for convergence. Alternative methodologies include the use of penalty methods \cite{alphonse2019b,bruggemann2023}, augmented Lagrangian methods \cite{kanzow2019}, or active set strategies \cite{alphonse2024}, which may be mesh-independent but either produce infeasible solutions or are confined to nodal low-order discretizations.

We now investigate the application of LVPP to QVIs. Since the obstacle appears explicitly in the equations to be solved in LVPP, it may, therefore, be treated in a straightforward manner as a function of the unknown solution. We restrict our discussion to a class of obstacle-type QVIs known as thermoforming problems \cite{alphonse2019,alphonse2019b,warby2003}. The objective of the thermoforming QVI is to find the equilibrium between the vertical displacement of a plastic membrane, denoted by $u \in H^1_0(\Omega)$, pushed upwards by a force $f \in L^2(\Omega)$ and heated to a temperature $T \in H^1(\Omega)$ into a metallic mold whose vertical displacement is denoted by $\Phi \in H^1(\Omega)$. In this case, $\Phi \coloneqq \Phi_0 + \xi T$, where $\Phi_0 \in H^1(\Omega)$ is the initial vertical displacement of the mold and $\xi \in C^2(\bar \Omega) \cap H^1_0(\Omega)$ is a given smoothing function. The smoothing function $\xi$ was introduced as a mathematical tool to ensure the existence of a thermoforming solution \cite{alphonse2019}. Physically, it causes the temperature to have a vanishingly small effect on the deformation of the mold close to the boundary of $\Omega$. Unlike in the obstacle problem of \Cref{sub:ObstacleProblem}, where the obstacle is fixed, the mold deforms due to the heat of the membrane, causing it to expand upwards. The deformation of the mold is modeled by a nonlinear screened Poisson equation that ensures a larger heat transfer from the membrane to the mold in regions where they are in close contact.  The behavior of the heat transfer is dictated by a positive conduction coefficient $\beta>0$ and a given globally Lipschitz non-increasing function $g \colon \mathbb{R} \to \mathbb{R}$ that defines an operator $g \colon H^1(\Omega) \to H^1(\Omega)$ (not relabeled). We define the feasible set as
\begin{align}
K(T) \coloneqq \{v \in H^1_0(\Omega) \mid v \leq \Phi = \Phi_0 + \xi T \}. 
\end{align}
The thermoforming QVI problem is to find $(u, T)\in K(T) \times H^1(\Omega)$ that satisfies, for all $(v, q) \in K(T) \times H^1(\Omega)$,
\begin{subequations} \label{eq:qvi}
\begin{align}
 (\nabla T, \nabla q) + \beta (T, q) &= (g(\Phi_0 + \xi T - u), q),\label{eq:qvi:b}\\
(\nabla u, \nabla(v-u)) &\geq (f, v-u).
\end{align}
\end{subequations} 

As in \Cref{sub:ObstacleProblem}, we choose the generalized Shannon entropy~\cref{eq:GeneralizedShannonEntropy} for the Legendre function in the LVPP subproblems (but with the relevant signs switched as we are dealing with a unilateral upper bound rather than a lower bound). The resulting saddle point formulation (in weak form) is:
for $\psi^0 = 0$, find $(u^{k}, \psi^{k}, T^{k}) \in H^1_0(\Omega) \times L^\infty(\Omega) \times H^1(\Omega)$ satisfying for all $(v,w,q) \in H^1_0(\Omega) \times L^\infty(\Omega) \times H^1(\Omega)$,
\begin{subequations} \label{eq:polyhedral:qvi:1}
\begin{align}
(\nabla T^k, \nabla q) + \beta(T^k, q) &=  (g(\exp(-\psi^{k})), q),
\label{eq:polyhedral:qvi:1c}\\
\alpha_{k} (\nabla u^k, \nabla v) + (\psi^k, v) &= \alpha_k (f, v) + (\psi^{k-1}, v), \\
(u^{k}, w) + (\exp(-\psi^{k}), w) &= (\Phi_0 + \xi T^k, w). \label{eq:polyhedral:qvi:1b}
\end{align}
\end{subequations}
The first equation \cref{eq:polyhedral:qvi:1c} arises thanks to the pointwise equality in \cref{eq:polyhedral:qvi:1b} and  may also be rewritten as
\begin{align}
\tag{\ref{eq:polyhedral:qvi:1c}'}
(\nabla T^k, \nabla q) + \beta(T^k, q) &= (g(\Phi_0 + \xi T^k - u^k), q),
\end{align}
which corresponds to \cref{eq:qvi:b} in the original QVI problem. We choose the former formulation as it achieves slightly better numerical convergence. 

We now fix the parameters: $\Omega = (0,1)^2$,  $\beta = 1$,  $\xi(x,y) = \sin(\pi x) \sin(\pi y)$,  $f \equiv 25$, $\Phi_0(x,y) = 1 - 2 \max(|x-1/2|, |y-1/2|)$, and
\begin{align*} 
g(s) = 
\begin{cases}
1 & s \leq 0,\\
1-100s & 0 < s < 10^{-2},\\
0 & s \geq 10^{-2}.
\end{cases}
\end{align*}
In \cref{tab:qvi}, we compare discretized LVPP with three alternative solvers that are experimentally observed to be mesh-independent: a direct Moreau--Yosida penalty scheme, a semismooth active set strategy \cite{alphonse2024}, and a fixed point approach. In all cases, we fix the same mesh with mesh size $h=0.01$ and devise a proximal Galerkin method by discretizing $u$, $T$, $\psi$ with continuous piecewise linear finite elements. We use the same discretization for $u$ and $T$ in the other solvers. The algorithms terminate once $\|u^{k} - u^{k-1}\|_{H^1(\Omega)} \leq 10^{-5}$. We note that proximal Galerkin could reach tighter tolerances, but the penalty schemes could not. 

We initialize proximal Galerkin with $(u^0, T^0) \equiv (0, 1)$, $\alpha_1 = 2^{-6}$ and adopt the update rule $\alpha_{k+1} = 4\alpha_k$. For robustness, a small Jacobian modification is added corresponding to $-10^{-10}\alpha_k^{-1}(\nabla \psi^k, \nabla w)$. Convergence was reached after solving 13 proximal Galerkin subproblems, reaching $\alpha_{13} = 2^{18}$, requiring a total of 20 Newton iterations, averaging 1.54 linear system solves per proximal Galerkin subproblem. Proximal Galerkin outperformed the other three solvers, resulting in the fastest run time\footnote{We acknowledge that parameter tuning may allow for faster run times in proximal Galerkin and the other solvers. Hence, the run times should be understood as a rough benchmark rather than a genuine ranking on solver performance.} and the fewest number of linear system solves.

The next best solver was the Moreau--Yosida penalty. This is a path following scheme where, at each Moreau--Yosida subproblem, one solves the following nonlinear and semismooth system:
\begin{subequations} \label{eq:polyhedral:qvi:2}
\begin{align}
(\nabla T^k, \nabla q) + \beta(T^k, q) -  (g(\Phi_0 + \xi T^k - u^k), q) &=0,
\label{eq:polyhedral:qvi:2b}\\
(\nabla u^k, \nabla v) + \gamma_k ((u^k - \Phi_0 - \xi T^k)_+, v)   &= (f, v),
\label{eq:polyhedral:qvi:2a}
\end{align}
\end{subequations}
where $(v)_+ = \max(0,v)$. The algorithm converges as $\gamma_k \to \infty$. The algorithm is initialized at $\gamma_1 = 1$ and we employ the $\gamma$-update rule described in \cite[\S4]{adam2019semismooth}.

The fixed point approach is a two-step iteration that first freezes the guess for $u^k$ and updates the obstacle via a nonlinear PDE solve, and subsequently freezes the obstacle, whereby $u$ is updated by a VI solve.
In other words: given $u^0$, for $k \in \mathbb{N}$ compute $T^k$ where, for all $q \in H^1_0(\Omega)$
\begin{subequations}
\label{eq:FixedPoint}
\begin{align}
(\nabla T^k, \nabla q) + \beta(T^k, q) -  (g(\Phi_0 + \xi T^k - u^{k-1}), q) &=0,
\end{align}
followed by the VI solve: find $u^k \in K(T^k)$ satisfying, for all $v \in K(T^k)$
\begin{align}
\label{eq:polyhedral:qvi:3b}
(\nabla u^k, \nabla(v-u^k)) &\geq (f, v-u^k).
\end{align}
\end{subequations}
Practically, \cref{eq:polyhedral:qvi:3b} is solved via a path-following Moreau--Yosida penalty scheme with the $\gamma$-update rule described in \cite[\S4]{adam2019semismooth}. The fixed point approach requires~\cref{eq:FixedPoint} to be a contraction. For this example, the contraction factor is close to one; as a result, we require 164 fixed point iterations to reach the termination criterion and a large wall clock solve time.

The semismooth active set strategy, derived in \cite[\S4.1]{alphonse2024}, is an acceleration of the fixed point approach that leads to local superlinear convergence. After each fixed point iteration in \cref{eq:FixedPoint}, one solves a single linear system akin to that found in an active set strategy \cite[Equations (62)--(64)]{alphonse2024}. This solver is well-posed on the infinite-dimensional level, which manifests as mesh-independent iteration counts.

From these preliminary results, it appears that LVPP is a promising choice for the numerical solution of QVIs.

\begin{figure}
  \begin{minipage}[c]{0.6\linewidth}
        \pgfplotsset{
    legend entry/.initial=,
    every axis plot post/.code={        \pgfkeysgetvalue{/pgfplots/legend entry}\tempValue
        \ifx\tempValue\empty
            \pgfkeysalso{/pgfplots/forget plot}        \else
            \expandafter\addlegendentry\expandafter{\tempValue}        \fi
    },
}
\newcommand{\equal}{=}


  \end{minipage}
    \begin{minipage}[c]{0.22\linewidth}
    \vspace*{-3mm}
  \includegraphics[clip, height=3.5cm]{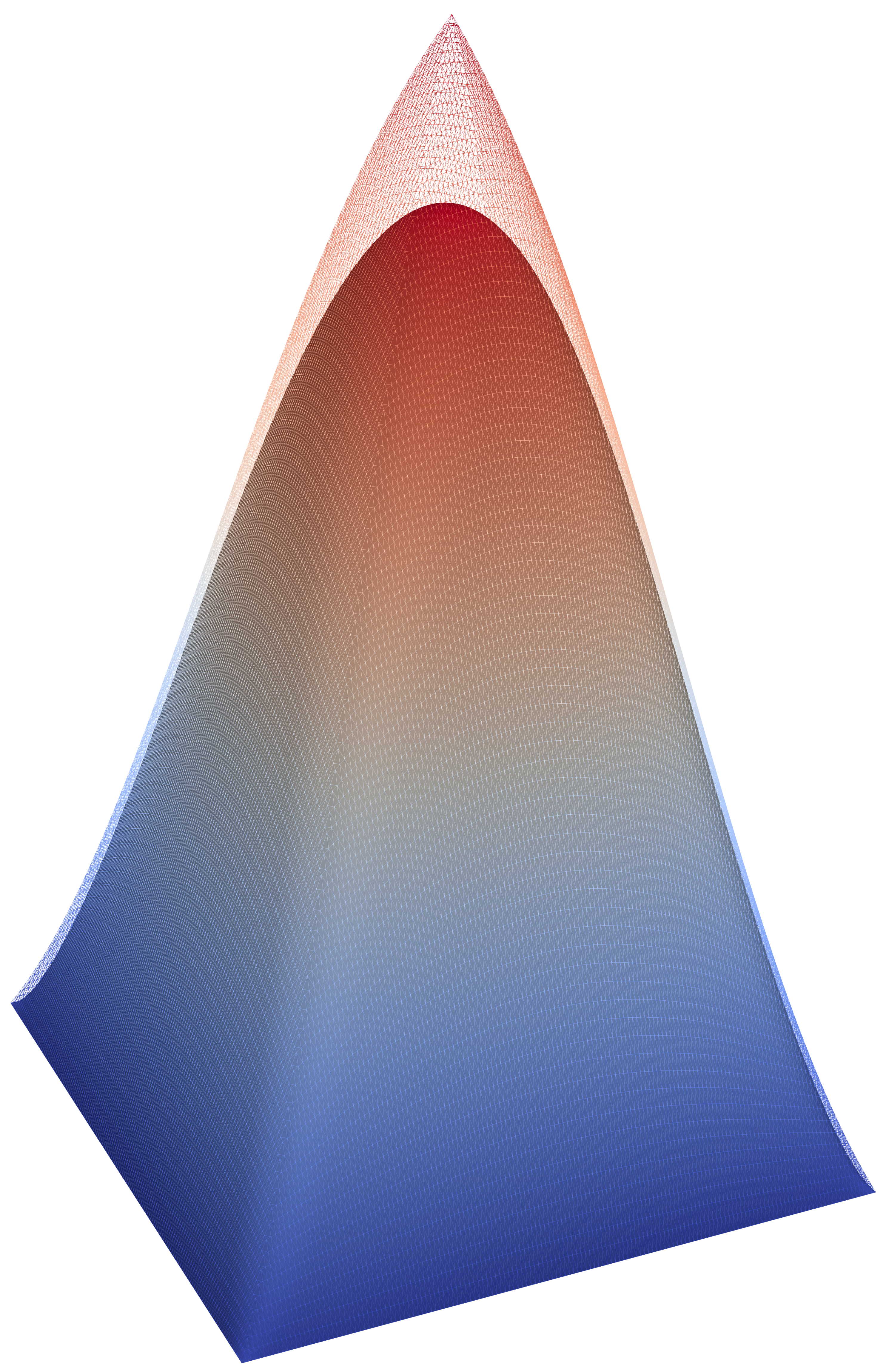}
  \end{minipage}
    \begin{minipage}[c]{0.06\linewidth}
  \includegraphics[height=2.5cm]{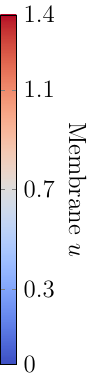}
  \end{minipage} 
\caption{Example 5. Left: One-dimensional slice of the membrane $u$ at $y=1/2$, mold $\Phi$, original mold $\Phi_0$, and temperature $T$ for the QVI problem described in \Cref{sec:qvi}. Right: Surface plot of the membrane $u$ (solid) and mold $\Phi$ (wireframe).}
\label{fig:qvi}
\end{figure}

\begin{table}[h!]
\renewcommand{\arraystretch}{1.3}
\centering
\small
\begin{tabular}{|l|c|c|c|}
\hhline{-|-|-|-|}
\rowcolor{lightgray!10} \multicolumn{1}{|c|}{Solver} &\multicolumn{1}{c|}{Outer loop} & \multicolumn{1}{c|}{Linear system solves}  & \multicolumn{1}{c|}{Run time (s)}  \\
\hline
Proximal Galerkin & 13 & 20 & 61.70 \\
Moreau--Yosida Penalty & 14 & 51 & 78.01 \\
Semismooth Active Set \cite{alphonse2024} & 7 & 236 & 112.60 \\
Fixed Point & 164 & 8493 & 3633.72\\ 
\hline
\end{tabular}
\caption{Example 5. The performance of four solvers, terminating when $\|u^{k} - u^{k-1}\|_{H^1(\Omega)} \leq 10^{-5}$. The Moreau--Yosida penalty solver and fixed point approach deliver an infeasible solution, whereas proximal Galerkin and the semismooth active set solver provide a feasible one.} 
\label{tab:qvi}
\end{table}

\section{Beyond bound constraints}
\label{sec:Non-Polyhedral}

We continue our study of LVPP by considering pointwise gradient and eigenvalue constraints.
These are important but less well-studied classes of problems, with fewer competing methods available.
The Lagrange multipliers for gradient constraints are much less regular than with bound constraints and, hence, much more subtle to discretize.
In particular, the higher regularity of the solution $u$ does not necessarily carry over to more regular Lagrange multipliers, as is often exploited to develop methods for pointwise bound constraints.
In addition, whereas one can, e.g.,~enforce a function's non-negativity using non-negative bases and coefficients, it is very difficult to enforce gradient or eigenvalue constraints by designing suitable bases.

We consider three representative problems.
The first is a simple plasticity model involving constraints on the norm of the gradient of a function $u \colon \Omega \to \mathbb{R}$.
The second considers enforcing natural bounds on the eigenvalues of an orientation variable in the Landau--de Gennes model for nematic liquid crystals.
The third demonstrates the treatment of intersections of inequality constraints using an obstacle problem with an additional gradient norm constraint.

\subsection{Example 6: Gradient norm constraints}
\label{sub:GradientNormConstraints}

Constraints on the gradient of the state are well-known from the optimal control literature; see, e.g., \cite{White1983,Mackenroth1986,Casas1993,Hintermller2010,Wollner2012}, where they are often used as a simplified model of a stress constraint. They also appear in the literature on (quasi-)variational inequalities; see the recent survey \cite{Rodrigues2019} and the references therein.
In particular, this type of constraint has been used to model the elastic-plastic torsion problem \cite{TWTing_1969,Brezis1971}, sandpile growth \cite{Prigozhin_1996,Antil_2022}, magnetization in type-II superconductors \cite{Prigozhin_1996a}, and in hydrology \cite{Prigozhin_1994}. Though the literature is decidedly more scarce for numerical methods for these types of variational inequalities compared to the obstacle problem, a quadratic penalty semismooth Newton approach has been successfully employed before, see e.g., \cite{Hintermller2012}.

Given a domain $\Omega \subset \mathbb{R}^n$ and a function $\phi \in L^{\infty}(\Omega)$ such that $\phi \ge \epsilon$ a.e.\ for some constant $\epsilon > 0$, a typical feasible set for gradient constraints is
\begin{equation}\label{eq:grad-constr}
K = \{ u \in H^1_0(\Omega) \mid |\nabla u| \le \phi \text{ a.e.\ in } \Omega \}\,.
\end{equation}
To obtain this feasible set $K$ from~\cref{eq:intro:feasible_set}, we set $V = H^1_0(\Omega)$, $B = \nabla$, $C(x)$ equal to the Euclidean ball of radius $\phi(x)$, and $\Omega_d = \Omega$.
For some $f \in L^2(\Omega)$, we consider minimizing the Dirichlet energy $J$, as defined in \cref{eq:dirichlet-energy},
in $K$. We use a modification of the Hellinger entropy~\cite{MTeboulle_2018}, which captures the geometry of the Euclidean ball, as the Legendre function:
\begin{equation}\label{eq:hellinger}
    R(a) := -\sqrt{\phi^2 - |a|^2}
    \,,
    \text{ with }
    \nabla R^\ast(a^\ast) = \frac{\phi}{\sqrt{1+|a^\ast|^2}} a^\ast
    \,.
\end{equation}
This leads to the following LVPP saddle point problem: for $\psi^0 = 0$, find $(u^k,\psi^k) \in H^1_0(\Omega) \times L^\infty(\Omega, \mathbb{R}^n)$ such that
\begin{subequations}
\label{eq:grad-lvpp}
\begin{align}
            \alpha_k(\nabla u^k, \nabla v) + (\psi^k, \nabla v)
            &= 
            \alpha_k(f, v) + (\psi^{k-1},\nabla v), \\
            (\nabla u^k, w) - \bigg(\frac{\phi\,\psi^k}{\sqrt{1+|\psi^k|^2}}, w\bigg)
            &= 0
\end{align}
\end{subequations}
for all $(v, w) \in H^1_0(\Omega) \times L^\infty(\Omega, \mathbb{R}^n)$.

\begin{figure}
\centering
    \begin{minipage}[c]{0.4\linewidth}
 {\includegraphics[clip, trim = 5cm 3cm 7cm 5cm, height=3cm]{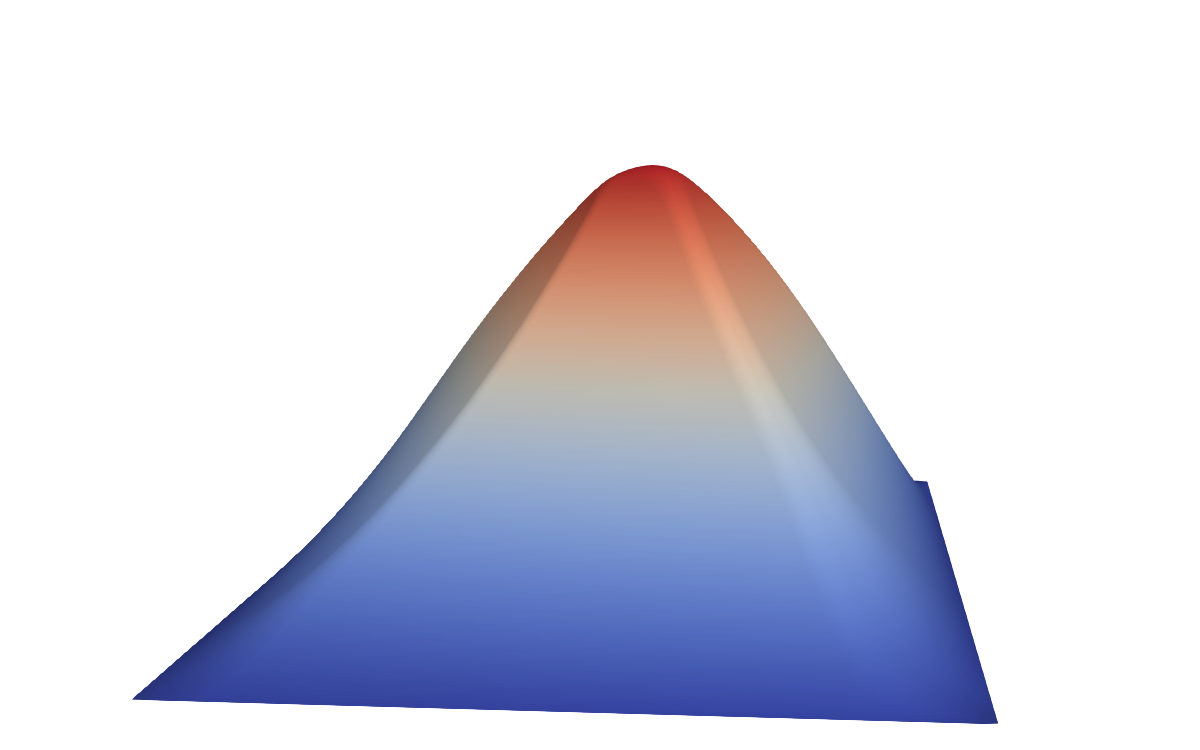}
 }
 \end{minipage}
 \begin{minipage}[c]{0.05\linewidth}
 \includegraphics[height=3cm]{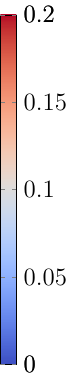}
 \end{minipage}
 \qquad
 \begin{minipage}[c]{0.28\linewidth}
 {\includegraphics[clip, trim = 6cm 6cm 6cm 6cm, height=3cm]{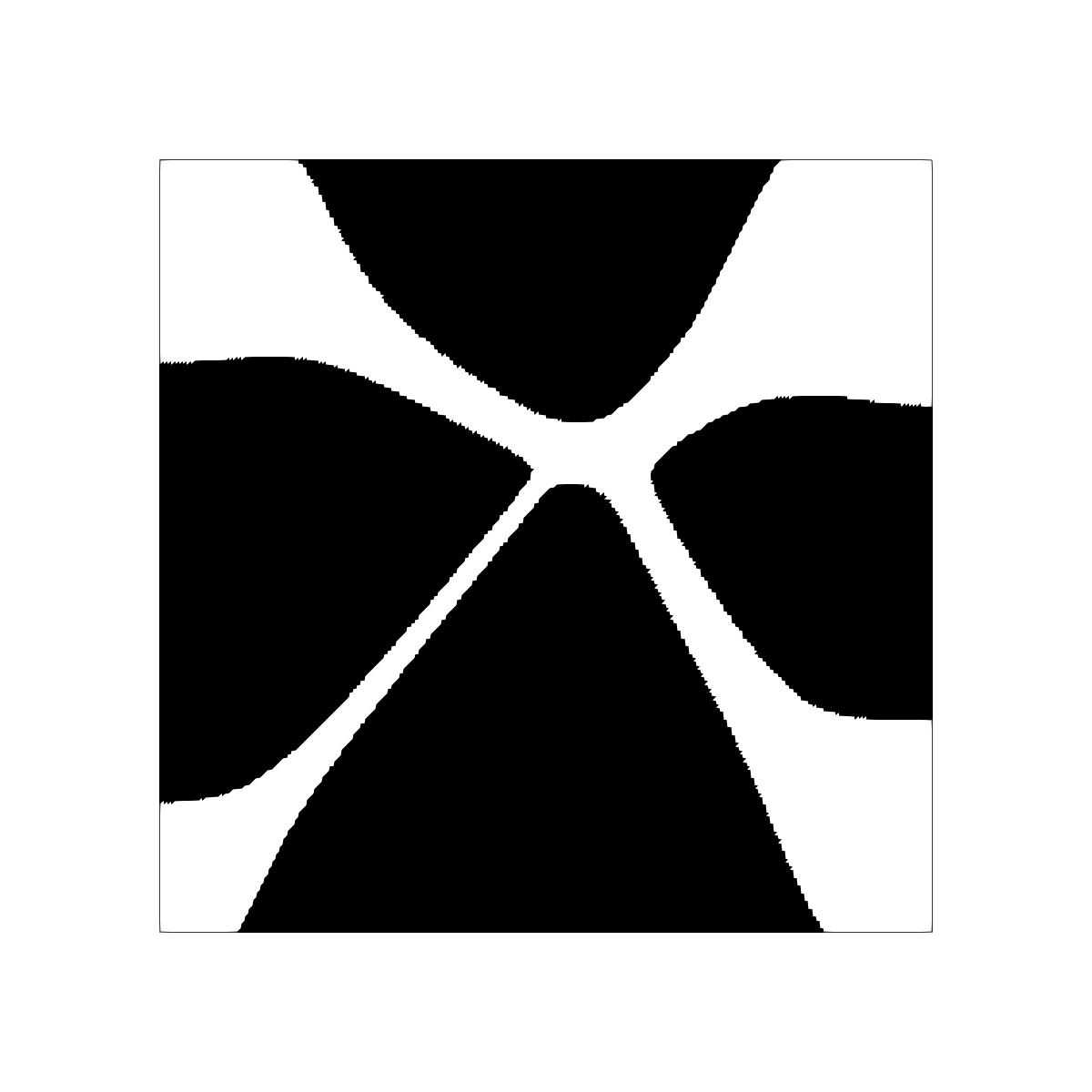}
 }
 \end{minipage}
    \caption{Example 6. Left: the discretized solution variable $u_h(x)$. Right: the estimated active set based on $\nabla R^*(\psi_h)$ (black). The algorithm was stopped once $\|u^k_h - u^{k-1}_h \|_{L^2(\Omega)} < 10^{-8}$.}
    \label{fig:grad-constr}
\end{figure}

For our numerical experiments, we solve~\cref{eq:grad-lvpp} on a unit square domain in $\mathbb R^2$ with homogeneous Dirichlet boundary conditions imposed on $u$. We set
\[
f(x_1,x_2) := 15 \sin^2(\pi x_1)
~~\text{ and }~~
\phi(x_1,x_2) := 0.1 + 0.2x_1 + 0.4x_2
\]
for all $(x_1,x_2) \in \Omega$.
We then discretize the underlying function spaces using continuous piecewise quadratic finite elements for $u$ and continuous piecewise linear finite elements for the latent variable $\psi$.
The stopping condition on~\cref{eq:grad-lvpp} is the $L^2(\Omega)$-norm of the increment in the primal variable from one proximal step to the next, $\|u^k_h - u^{k-1}_h \|_{L^2(\Omega)} < \mathtt{tol.}$, with the tolerance set to $10^{-8}$. The proximity parameter begins with $\alpha_1 = 1$ and doubles at each new proximal step.
The tolerance for Newton's method is set to $10^{-8}$. 
For a mesh of $200 \times 200$ cells, the method requires a total of 17 proximal iterations, amounting to 27 linear solves.
The maximum number of Newton steps was 4, and the minimum was 1, with the final 12 iterations requiring only one Newton step per proximal iteration. We plot the computed solution and approximated active set  \Cref{fig:grad-constr}. For the latter, we plot the characteristic function for the set where $|\nabla R^*(\psi_h)| \ge \phi - 10^{-8}$.

\subsection{Example 7: Eigenvalue constraints} \label{sec:ldg}

We illustrate how to apply LVPP to enforce eigenvalue bounds. This class of constraints is very difficult to enforce but also very important in engineering design.

As a didactic example, we enforce the eigenvalue constraints of the Landau--de Gennes model for nematic liquid crystals. In this model, the orientation of the liquid crystal is described by a tensor-valued order parameter $Q(x)$~\cite{degennes1995,ball2017,wang2021}. More specifically, on a domain $\Omega \subset \mathbb{R}^n$, at each point $Q$ is a symmetric and traceless $n \times n$ matrix, and must satisfy the eigenvalue constraints
\begin{equation} \label{eq:ldg:eigenvalue_constraints}
\lambda_i(Q) \in [-1/n, (n-1)/n], \quad i= 1, \dots, n.
\end{equation}
The equilibrium configuration of the liquid crystal is attained by solving a minimization problem with a suitable energy functional. In almost all cases, the eigenvalue constraints are not explicitly enforced in the numerical optimization (e.g.,~\cite{robinson2017,surowiec2023}), although they are usually satisfied through careful choice of the bulk energy term~\cite{ball2010,majumdar2010,ball2017}. A notable exception is the barrier approach of Wang \& Xu~\cite{wang2022}, which adds a quasi-entropy to the functional to be minimized that blows up as $\lambda(Q)$ approaches the bounds.

We consider the setup of Robinson et al.~\cite{robinson2017}, with a slightly modified nondimensionalization. The problem is to minimize the energy functional
\begin{equation} \label{eq:ldg:energy}
    J(u)
    =
    \frac{1}{2}\int_\Omega \nabla Q \fcolon \nabla Q \dd x
  + \frac{1}{2}\int_\Omega A \tr(Q^2) \dd x
  + \frac{1}{4}\int_\Omega C (\tr(Q^2)^2) \dd x
    \,,
\end{equation}
on the domain $\Omega = (0, 1)^2$, with parameters $A = 1, C = 4$. Here $\fcolon$ denotes the Frobenius inner product of third order tensors, analogous to $\cdot$ and $:$ for vectors and matrices.
In two dimensions, the feasible set simplifies to
\begin{equation} \label{eq:ldg:feasible_set}
    K = \{ Q \in H^1_g(\Omega, \mathbb{R}^{2 \times 2}_{\text{sym,tr}}) \mid -\frac{1}{2} I \preceq Q \preceq \frac{1}{2} I \text{ a.e.\ in } \Omega\}
    \,.
\end{equation}
Here $H^1_g(\Omega, \mathbb{R}^{2 \times 2}_{\text{sym,tr}})$ is the set of $H^1$ symmetric traceless matrices enforcing the boundary condition
\begin{equation}
g(x, y) =
\begin{pmatrix}
\frac{1}{2} s(x, y) \cos(2\theta) & \frac{1}{2} s(x, y) \sin(2\theta) \\
\frac{1}{2} s(x, y) \sin(2\theta) & -\frac{1}{2} s(x, y) \cos(2\theta) \\
\end{pmatrix},
\end{equation}
where $\theta = 0$ on $y \in \{0, 1\}$, $\theta = \pi/2$ on $x \in \{0, 1\}$, and the nematic ordering $s$ ramps to zero at the corners:
\begin{equation}
s(x, y) =
\begin{cases}
T(x) & y \in \{0, 1\}, \\
T(y) & x \in \{0, 1\},
\end{cases}
\end{equation}
with ramp function
\begin{equation}
T(z) = \begin{cases}
z/d & z \in [0, d), \\
1   & z \in [d, 1-d), \\
(1-z)/d & z \in [1-d, 1],
\end{cases}
\end{equation}
for ramp scale $d = 0.06$.

To obtain $K$ from~\cref{eq:intro:feasible_set}, we set $V = H^1_g(\Omega, \mathbb{R}^{2 \times 2}_{\text{sym,tr}})$, $B = I$, $C(x)$ to be the set of all symmetric $n \times n$ matrices with spectral radius less than or equal to $1/2$, and $\Omega_d = \Omega$.
We suggest selecting the entropy
\begin{align*}
    \begin{split}
    R(a) &= \operatorname{tr}\big( (a + I/2) \ln (a + I/2)
    +
    (I/2 - a) \ln (I/2 - a)
    \big)
    \,,\\
    \text{ with }
    \nabla R^\ast (a^\ast) &= \tanhm (a^\ast/2)/2
    \,,
    \end{split}
\end{align*}
where $\ln$ is the matrix logarithm \cite{petersen2012matrix} and $\tanhm$ is the matrix $\tanh$ function defined below.
Since $\tanh: \mathbb{R} \to (-1, 1)$, any output of the induced matrix function $\tanhm$ has spectral radius at most one. We implement $\tanhm$ using the formula
\begin{equation}
\tanhm(a) = 2 (\exp(2a) + I)^{-1} (\exp(2a) - I)
\end{equation}
and in turn implement the matrix exponential $\exp$ using the formulae of Bernstein \& So~\cite{bernstein1993} in two dimensions\footnote{There is a minor typographical error in~\cite{bernstein1993}. The prefactor in case \emph{i)} of Corollary 2.4 should be $\exp{((a + d)/2)}$, not $\exp{(a + d/2)}$.} and Cheng \& Yau~\cite{cheng1997} in three dimensions\footnote{There are several minor typographical errors in~\cite{cheng1997}. First, the indices $i$ and $j$ in (22) should be swapped. Second, the key formula (28') for the case with three distinct eigenvalues is missing several minus signs---the coefficients of $A^2$ and $I$ should be negated.}. In particular, with a careful implementation in the symbolic Unified Form Language~\cite{alnaes2012} employed by FEniCS and Firedrake, the derivatives of $\tanhm$ and $\exp$ needed in a Newton iteration are computed automatically.

In two dimensions, since $\tanh$ is odd, if $a$ is traceless, then $\tanhm(a)$ is also traceless. The natural space for the latent variable is therefore $L^\infty(\Omega, \mathbb{R}^{2 \times 2}_\mathrm{sym,tr}).$\footnote{For $n > 2$ the latent variable would not necessarily be traceless.} This leads to the following LVPP iteration: for $\psi^0 = 0$, find $(Q^k, \psi^k) \in H^1_g(\Omega, \mathbb{R}^{2 \times 2}_{\text{sym,tr}}) \times L^\infty(\Omega, \mathbb{R}^{2 \times 2}_\mathrm{sym,tr})$ such that
\begin{subequations}
\label{eq:ldg:lvpp}
\begin{align}
\alpha_k J'(Q; V) + (\psi^k, V) &= (\psi^{k-1}, V) \\
(Q, w) - (\frac{1}{2} \tanhm(\psi/2), w) &= 0 \label{eq:ldg:lvpp:b}
\end{align}
\end{subequations}
for all $(V, w) \in H^1_0(\Omega, \mathbb{R}^{2 \times 2}_{\text{sym,tr}}) \times L^\infty(\Omega, \mathbb{R}^{2 \times 2}_\mathrm{sym,tr})$.
For $n > 2$, \cref{eq:ldg:lvpp:b} would change to
\begin{equation}
(Q, w) - \left(\frac{1}{2} \tanhm(\psi/2) + \frac{n-2}{2n}I, w\right) = 0 \label{eq:ldg:lvpp:c}
\end{equation}
to scale and shift the eigenvalues into the physical range.

\begin{figure}
\centering
\begin{minipage}[c]{4cm}
\begin{tikzpicture}
    \node[anchor=south west,inner sep=0] (image) at (0,0) {\includegraphics[width=\textwidth]{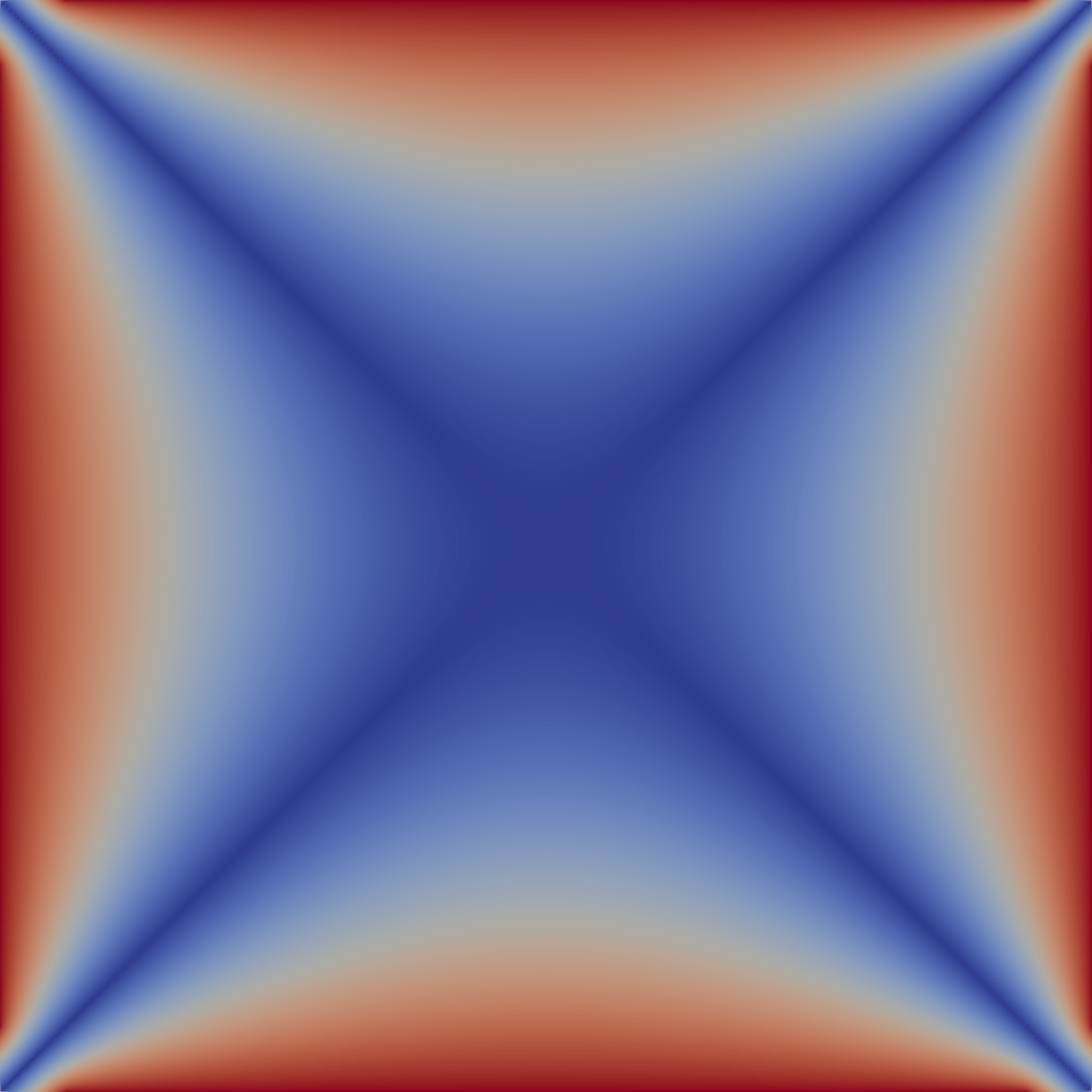}};
\end{tikzpicture}
\end{minipage}
\hspace*{0.2cm}
 \begin{minipage}[c]{0.8cm}
     \includegraphics[width=\linewidth]{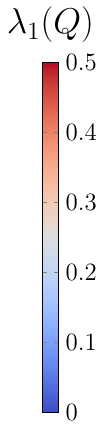}
 \end{minipage}
\caption{Example 7. The largest eigenvalue $\lambda_1(Q)$ of the minimizer of the Landau--de Gennes energy~\cref{eq:ldg:energy} over \cref{eq:ldg:feasible_set}. The other eigenvalue of the matrix is its negation.}
\label{fig:ldg}
\end{figure}

We discretize with piecewise cubic continuous finite elements on a $100 \times 100$ regular quadrilateral mesh for the components of both $Q$ and $\psi$. Since both variables are symmetric and traceless, they are described by two scalar fields, i.e.,~we employ the ans\"atze
\begin{equation}
Q = \begin{pmatrix}
Q_1 & Q_2 \\ Q_2 & -Q_1
\end{pmatrix},
\quad
\psi = \begin{pmatrix}
\psi_1 & \psi_2 \\ \psi_2 & -\psi_1
\end{pmatrix}.
\end{equation}
The solver took 6 proximal steps, requiring a total of 11 Newton steps. The converged solution, depicted in \Cref{fig:ldg}, satisfies the eigenvalue constraints \cref{eq:ldg:eigenvalue_constraints} pointwise. The solution is the well-known well order-reconstruction solution discovered by Kralj \& Majumdar~\cite{kralj2014}.

\subsection{Example 8: Intersections of constraints}
It is straightforward to apply LVPP when $K$ is defined as the intersection of several sets. We demonstrate this by minimizing the Dirichet energy $J$, as defined in \cref{eq:dirichlet-energy} with $f \equiv 0$
over the convex feasible set
\begin{equation}
    K = \{ u \in H^1_0(\Omega) \mid u \geq \phi_0 \text{ and } |\nabla u| \le \phi ~\text{ a.e.\ in } \Omega \}
    \,,
\end{equation}
which can be seen as the intersection of~\cref{eq:polyhedral:unilateral_K} and~\cref{eq:grad-constr}. We assume $\phi$ and $\phi_0$ satisfy suitable consistency conditions so that solutions to this problem exist.
To obtain $K$ from~\cref{eq:intro:feasible_set}, we set $V = H^1_0(\Omega)$, $B = (\operatorname{id}, \nabla)^\top$, $C(x)$ to be the Cartesian product the closed interval $[\phi_0(x), \infty)$ and the Euclidean ball of radius $\phi(x)$, and $\Omega_d = \Omega$. The geometry of this Cartesian product is captured by adding the entropies associated with the independent components:
\begin{align}
\begin{split}
    R((a_0,a))
    &=
    (a_0 - \phi_0) \ln (a_0 - \phi_0) - (a_0 - \phi_0) - \sqrt{\phi^2 - |a|^2}
    \,,\\
    \text{ with }
    \nabla R^\ast((a_0^\ast,a^\ast))
    &=
    \begin{pmatrix}
        \phi_0 + \exp a_0^\ast \\[3pt]
        \dfrac{\phi\, a^\ast}{\sqrt{1+|a^\ast|^2}}
    \end{pmatrix}
    \,.
\end{split}
\end{align}
For fixed $a_0$, the geometry of the ball of radius $\phi(x)$ is captured by the Hellinger entropy; for fixed $a$, the geometry of $[\phi_0(x), \infty)$ by the Fermi--Dirac entropy.

Since the gradient isomorphism $\nabla R^*$ has two components, we introduce two latent variables, just as a problem with two inequality constraints requires two Lagrange multipliers.
This leads to the following saddle point problem: for $(\psi_0^0,\psi^0) = (0,0)$, find $(u^k,\psi_0^k,\psi^k) \in H^1(\Omega) \times L^\infty(\Omega) \times L^\infty(\Omega, \mathbb{R}^n)$ such that
\begin{subequations}
\label{eq:grad-intersection-lvpp}
\begin{align}
            \alpha_k(\nabla u^k, \nabla v) + (\psi_0^k, v) + (\psi^k, \nabla v)
            &=
            (\psi_0^{k-1}, v) + (\psi^{k-1},\nabla v) \\
            (u^k, w_0) - (\exp \psi_0^k, w_0)
            &= (\phi_0, w_0)
            \\
            (\nabla u^k, w) - \bigg(\frac{\phi\,\psi^k}{\sqrt{1+|\psi^k|^2}}, w\bigg)
            &= 0
\end{align}
\end{subequations}
for all $(v, w_0, w) \in H^1(\Omega) \times L^\infty(\Omega) \times L^\infty(\Omega, \mathbb{R}^n)$.

We solve this problem for $\Omega = (0, 1)$, and obstacle and gradient constraints given by
\begin{align}
\begin{split}
\phi_0(x) &= \begin{cases}
c \exp{\left(-\frac{1}{10(x - 0.2)(0.8 - x)}\right)} & x \in (0.2, 0.8), \\
0 & \text{otherwise},
\end{cases}\\
\phi(x) &= \begin{cases}
\phi_c & x \in [0, 0.2] \cup [0.8, 1], \\
100 & \text{otherwise},
\end{cases}
\end{split}
\end{align}
where the normalisation constant $c \in \mathbb{R}$ is chosen so that $\phi_0(0.5) = 1$, and $\phi_c > 0$ is to be varied.
We discretize \cref{eq:grad-intersection-lvpp} with continuous piecewise linear finite elements for all variables $(u, \psi_0, \psi)$.

\begin{figure}
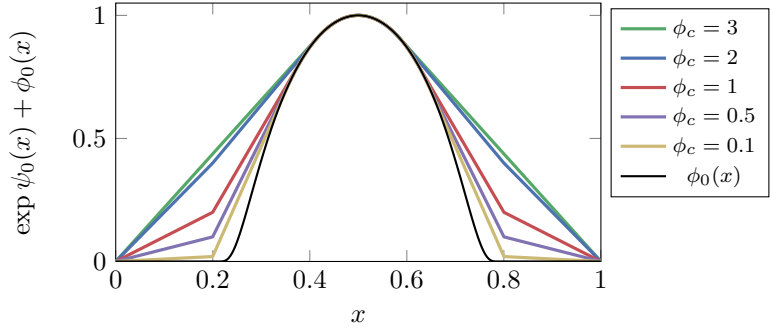

    \centering
        \pgfplotsset{
    legend entry/.initial=,
    every axis plot post/.code={        \pgfkeysgetvalue{/pgfplots/legend entry}\tempValue
        \ifx\tempValue\empty
            \pgfkeysalso{/pgfplots/forget plot}        \else
            \expandafter\addlegendentry\expandafter{\tempValue}        \fi
    },
}
\newcommand{\equal}{=}



    \caption{Example 8. Solutions to \cref{eq:grad-intersection-lvpp} with both obstacle and gradient constraints. The value of $\phi_c$ determines the gradient constraint imposed on $[0, 0.2] \cup [0.8, 1]$, while $\phi_0(x)$ defines the obstacle constraint.}
    \label{fig:grad-intersection-lvpp}
\end{figure}

The solutions for different values of $\phi_c$ are shown in \Cref{fig:grad-intersection-lvpp}. We plot the obstacle-conforming approximation $\phi_0 + \exp{\psi_0}$ (the primal solutions are similar). As $\phi_c$ is reduced, the gradient constraint becomes more and more restrictive, changing both the slope on $[0, 0.2] \cup [0.8, 1]$ and the active set in contact with the obstacle.

\section{Applications to fully-nonlinear PDEs}
\label{sec:NonlinearPDEs}
In addition to the applications demonstrated above, the LVPP methodology has a remarkable connection to certain fully-nonlinear PDEs.
In \Cref{ssec:eikonal}, we first illustrate how LVPP leads to new mixed forms for first-order nonlinear PDEs.
We provide computational evidence of the utility of such an approach by solving the eikonal equation.
We then close the section by deriving a new geometric formulation of the classical Monge--Amp\`ere equation by using matrix-valued latent variables and the matrix exponential function.

\subsection{Example 9: The eikonal equation}\label{ssec:eikonal}
Many nonlinear first-order 
PDEs have the following form \cite{Dacorogna1999-mj}:
\begin{equation}\label{eq:nonlinear-pde}
\text{find $u: \Omega \to \mathbb R$ such that  }
F(x,u,\nabla u) = 0 \text{ in } \Omega
\text{ and } u = g \text{ on }\partial \Omega\,.
\end{equation}
Generally, these equations lack a divergence structure, so one cannot integrate by parts to arrive at the correct notion of a weak solution.
Instead, for suitable $F$,
we opt for the concept of a viscosity solution introduced by Crandall and Lions \cite{Lions1982-jx,Crandall1983}.

Viscosity solutions are directly linked to constrained optimization; see, e.g., \cite[Proof of Thm.\ 1, Thm.\ 2]{SZagatti_2009}. 
Specifically, the viscosity solution of~\cref{eq:nonlinear-pde} minimizes
\begin{subequations}
\label{eqs:Viscosity}
\begin{equation}
\label{eq:ObjectiveFunctionViscosity}
J(u) = -\int_{\Omega} u \dd{x}
\end{equation}
over the feasible set
\begin{equation}
\label{eq:FeasibleSetViscosity}
K = \big\{
u \in W^{1,\infty}(\Omega) \mid F(x,u,\nabla u) \le 0 \text{ f.a.e.\ $x$ in } \Omega \text{ and } u|_{\partial \Omega} = g
\big\}
\,.
\end{equation}
\end{subequations}

Putting this feasible set $K$ into our framework for general mappings $F$ requires a generalization of~\cref{eq:intro:feasible_set} that would significantly complicate the paper.
However, there are many interesting examples that we can immediately consider. For instance,  upon defining $F(x,u,\nabla u) = |\nabla u| - 1$ and taking $g \equiv 0$, the minimizer of~\cref{eqs:Viscosity} is the viscosity solution of the classical eikonal equation \cite[1.1.2]{Dacorogna1999-mj},
\begin{equation}
\label{eq:EikonalStrongForm}
    |\nabla u| = 1 \text{ in } \Omega ~\text{ and }~ u = 0 \text{ on }\partial \Omega.
\end{equation}
To obtain~\cref{eq:FeasibleSetViscosity} from~\cref{eq:intro:feasible_set} for this problem, we need only set $V = H^1_0(\Omega)$, $B = \nabla$, $C$ to be the Euclidean ball of radius $1$, and $\Omega_d = \Omega$.

The two most well-known methods for solving the eikonal equation are the fast marching \cite{Osher1988,JASethian_1996} and the fast sweeping \cite{Zhao2000,Zhao2004} methods. They are both finite difference-based schemes that trace back to Dijkstra's algorithm \cite{Dijkstra1959}.
Both methods require specialized data structures that present a barrier to their integration into many of the high-order finite element codes used by the scientific computing community.
On the other hand, the LVPP-based approach derived below is immediately compatible with such codes as it involves a Galerkin discretization with widely available finite elements.

Selecting the Hellinger entropy~\cref{eq:hellinger} with $\phi \equiv 1$ and following the general framework for gradient-constrained problems described in \Cref{sub:GradientNormConstraints}, we deduce that minimizing~\cref{eq:ObjectiveFunctionViscosity} over the feasible set~\cref{eq:FeasibleSetViscosity} leads to the following sequence of LVPP subproblems:
for $\psi^0 = 0$, find $(u^{k}, \psi^{k}) \in H^1_0(\Omega) \times L^2(\Omega,\mathbb{R}^n)$ satisfying
\begin{subequations}
\label{eq:eikonal-lvpp-H1}
\begin{align}
    (\psi^k, \nabla v)
    &= 
    \alpha_k(1, v) + (\psi^{k-1},\nabla v), \\
    (\nabla u^k, w) - \bigg(\frac{\psi^k}{\sqrt{1+|\psi^k|^2}}, w\bigg)
    &= 0,
\end{align}
\end{subequations}
for all $(v,w) \in H^1_0(\Omega) \times L^2(\Omega,\mathbb{R}^n)$.
Similar systems of PDEs arise in strain-limiting elastic models \cite{beck2017existence,bulivcek2014,bulivcek2015}. These subproblems may also be seen as nonlinear Darcy flow models (cf.\ \cite[Chapter 51]{Ern2021}), which we highlight by integrating each of the differential terms in~\cref{eq:eikonal-lvpp-H1} by parts and swapping the order of the equations.
Formally, these operations suggest a change to the functional setting, revealing a new sequence of subproblems: for $\psi^0 = 0$, find $(u^{k}, \psi^{k}) \in L^2(\Omega) \times H(\operatorname{div},\Omega)$ satisfying
\begin{subequations}
\label{eq:eikonal-lvpp-Hdiv}
\begin{align}
    \bigg(\frac{\psi^k}{\sqrt{1+|\psi^k|^2}}, w\bigg) + (u^k, \nabla\cdot w)
    &= 0,
\label{eq:eikonal-lvpp-Hdiv-1}
    \\
    (\nabla\cdot \psi^k, v)
    &= 
    (\nabla\cdot\psi^{k-1},v) - \alpha_k(1,v),
\label{eq:eikonal-lvpp-Hdiv-2}
\end{align}
\end{subequations}
for all $(v,w) \in L^2(\Omega) \times H(\operatorname{div},\Omega)$.
As with most LVPP formulations in this paper, further work is required to rigorously establish the well-posedness of both~\cref{eq:eikonal-lvpp-H1,eq:eikonal-lvpp-Hdiv}.

\begin{figure}[t]
\centering
\begin{minipage}{0.34\textwidth}
    \centering
    \includegraphics[width=0.95\linewidth,keepaspectratio]{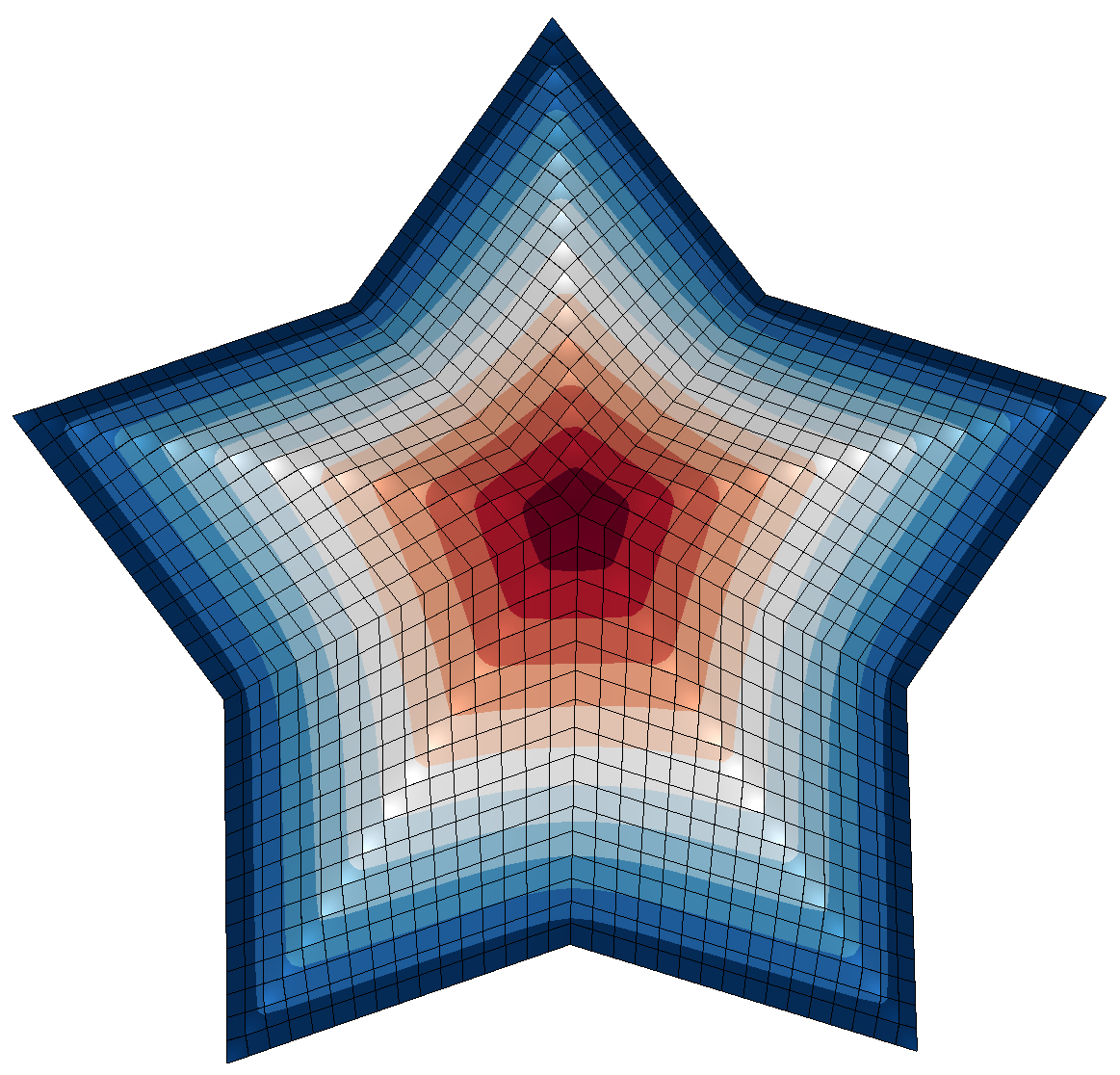}
\end{minipage}
\begin{minipage}{0.25\textwidth}
\centering
    \includegraphics[width=0.95\linewidth,keepaspectratio]{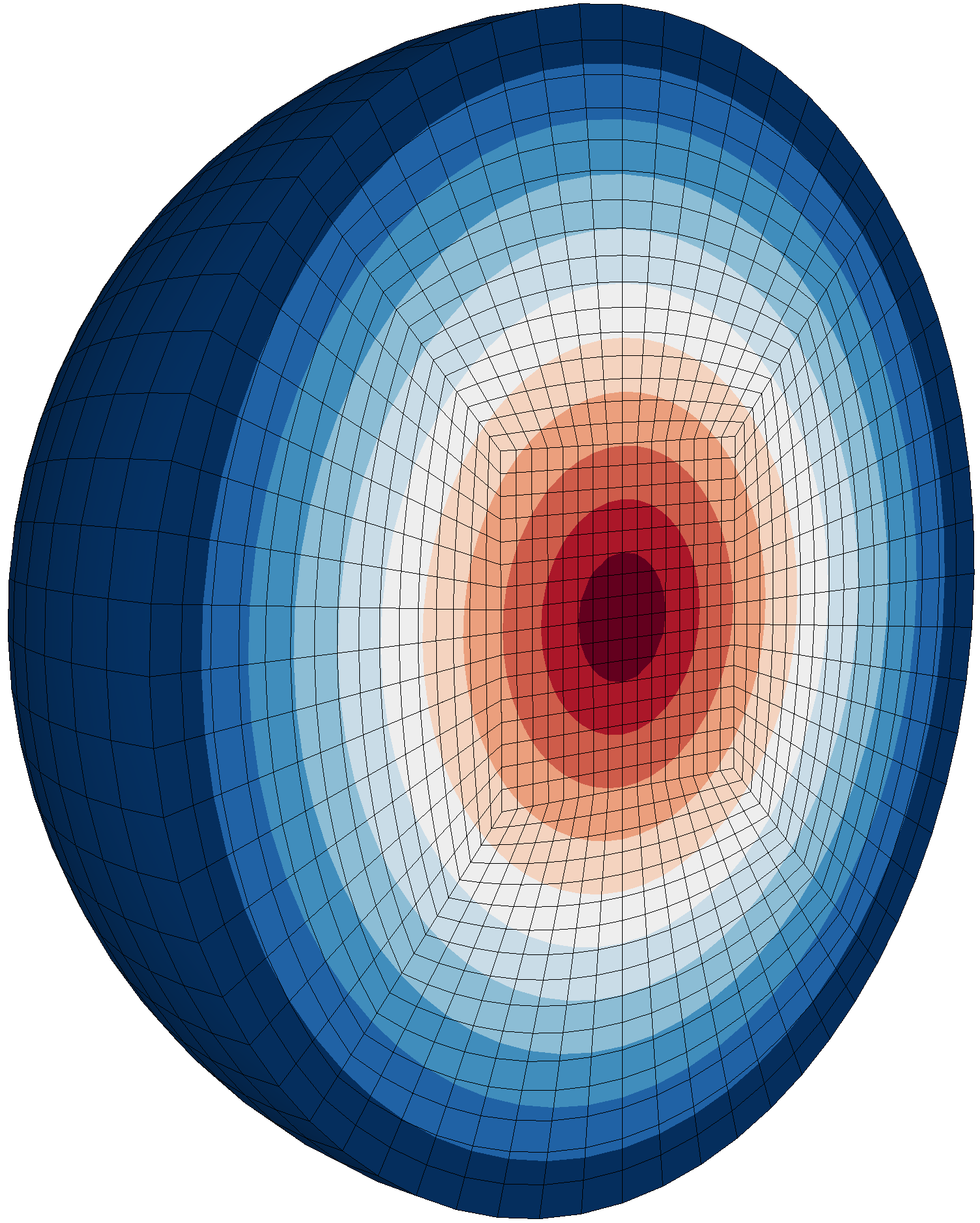}
\end{minipage}
\begin{minipage}{0.32\textwidth}
    \centering
    \includegraphics[width=0.95\linewidth,keepaspectratio]{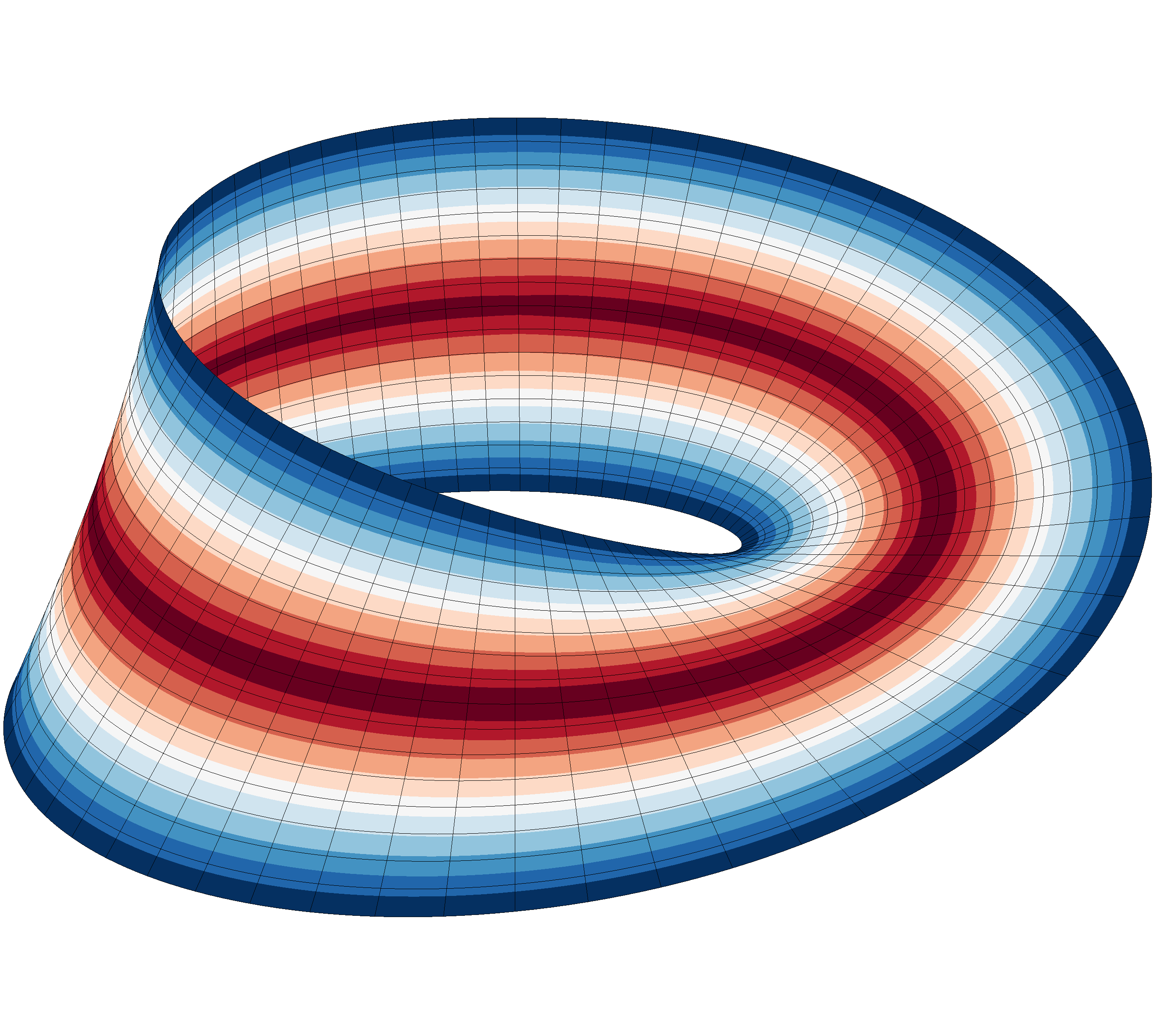}
\end{minipage}
\begin{minipage}{0.05\textwidth}
\centering
  \includegraphics[width=0.9\textwidth]{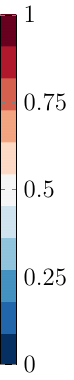}
\end{minipage}
\caption{Example 9. Viscosity solutions of the eikonal equation on a star-shaped domain (left), sphere (middle), and M\"obius strip (right).}
\label{fig:eikonal-lvpp}
\end{figure}

We solved~\cref{eq:eikonal-lvpp-Hdiv} on three different domains with the computed solutions presented in \Cref{fig:eikonal-lvpp}: a two-dimensional star-shaped domain, a solid (three-dimensional) sphere, and a M\"obius strip.
To compare the results, we rescaled each domain so that the maximum, taken over all points in the domain, of the shortest distance to the domain boundary would be one.
On the star-shaped domain and the sphere, we replaced $L^2(\Omega)$ in \cref{eq:eikonal-lvpp-Hdiv} with a first-order discontinuous piecewise polynomial space and $H(\operatorname{div},\Omega)$ with second-order Raviart--Thomas polynomial space, as this is the standard choice for Darcy flow problems.
Showcasing additional flexibility in the discretization choices, we used a continuous piecewise polynomial discretization of Taylor--Hood-type inspired by \cite{karper2009unified} for the M\"obius strip example.
In particular, we used second-order Lagrange elements to approximate $\psi^k$ and first-order Lagrange elements for $u^k$.
We then set $\alpha_k = 10\cdot\min\{2^k, 5\}$ and solved each discretized subproblem with a damped Newton method.
In each case, the solver took three Newton steps on the first subproblem ($k=1$) but only one Newton step for each subsequent subproblem ($k>1$).
We stopped the solver once $\|u^k - u^{k-1}\|_{L^2(\Omega)} < 10^{-4}$, which always occurred in under ten iterations.

The setting above easily extends to solving other related problems. 
Indeed, selecting the general Hellinger entropy given in~\cref{eq:hellinger} results in multiplying the nonlinear function in~\cref{eq:eikonal-lvpp-Hdiv-1} by $\phi$, and allows one to treat the general eikonal equation, $|\nabla u| = \phi$.
Alternatively, generalizing the objective function to $J(u) = -\int_{\Omega} u f \dd{x}$ and minimizing it over the feasible set $K = \{ u \in H^1_0(\Omega) \mid |\nabla u| \le 1 \text{ a.e.\ $x$ in } \Omega \}$
allows one to characterize the $p \to \infty$ limit of solutions to the $p$-Laplace equation with $f \in C(\overline{\Omega})$ \cite{ishii2005limits}:
\[
    -\Delta_p u = f \text{ in } \Omega ~\text{ and }~ u = 0 \text{ on }\partial \Omega.
    \,
\]
To arrive at an algorithm to solve this problem, we need only to replace the right-most term in~\cref{eq:eikonal-lvpp-Hdiv-2} by $-\alpha_k(f,v)$.

\subsection{Example 10: The Monge--Amp\`ere equation}\label{ssec:mongeamp}
One can also define viscosity solutions to problem~\cref{eq:nonlinear-pde} if $F$ depends on the Hessian $\nabla^2 u$ \cite{Dacorogna1999-mj,Evans2010,gilbarg1977elliptic}.
Thus, the innovations above raise the question of whether our techniques can also be used to solve fully nonlinear \textit{second-order} PDEs.
We choose to provide a partial answer to this question by considering the celebrated Monge--Amp\`ere equation, in which
$F(x,u,\nabla u, \nabla^2 u) = \operatorname{det}( \nabla^2 u ) - \rho$.
More specifically, given uniformly positive functions $\rho \in C(\overline{\Omega})$ and $g \in C^3(\overline{\Omega})$ defined over a smooth, bounded, and uniformly convex domain $\Omega \subset \mathbb R^n$, we seek the unique convex function $u \in H^2(\Omega) \cap H^1_g(\Omega)$ satisfying
\begin{equation}\label{eq:monge-ampere}
    \operatorname{det}( \nabla^2 u ) = \rho
    ~\text{ in } \Omega
    \,,
    \qquad
    u = g
    ~\text{ on } \partial\Omega
    \,.
\end{equation}
We refer to \cite{Caffarelli1990,Savin2013} for existence and regularity theory in this setting. 

Although the Monge--Amp\`ere equation is related to optimal transport \cite{santambrogio2015}, that connection to the optimization literature does not yield an immediate insight into our approach.
Instead, we focus on the fact that the set of admissible solutions to~\cref{eq:monge-ampere},
\[
    K
    =
    \{
        u \in H^2(\Omega) \cap H^1_g(\Omega) 
        \mid
        \nabla^2 u \succeq 0 \text{ a.e.\ in } \Omega
    \}
    ,
\]
fits into the setting of \cref{eq:intro:feasible_set} by taking $V = H^2(\Omega) \cap H^1_g(\Omega)$, $B = \nabla^2$, and $C$ equal to the set of all symmetric positive semidefinite $n \times n$ matrices.
We can now invoke a Legendre function to construct a structure-preserving change of variables reformulation of~\cref{eq:monge-ampere}.

Consider the (unnormalized negative) von Neumann entropy \cite{Neumann1927,Araki1970} for positive semidefinite symmetric matrices, defined by 
\begin{gather*}
    R(a) = \operatorname{tr}( a \ln a - a )
    \,,
    ~\text{ with }
    \nabla R^\ast (a^\ast) = \exp a^\ast
    .
\end{gather*}
The convex geometry of $K$ can be encoded by introducing a symmetric matrix-valued latent variable $\psi$ via the equation $B u = \nabla R^\ast (\psi)$; i.e., we define
\begin{subequations}
\label{eq:ig-monge}
\begin{equation}\label{eq:ig-monge-a}
    \psi = \ln \nabla^2 u
    \quad\iff\quad
    \exp \psi = \nabla^2 u.
\end{equation}
Finally, utilizing the identity $\mathrm{det}(\exp \psi) = \exp(\mathrm{tr}\, \psi)$ \cite{haber2018notes} and the original PDE~\cref{eq:monge-ampere}, we see that $\psi$ must satisfy the algebraic relation
\begin{equation}\label{eq:ig-monge-b}
    \operatorname{tr}\psi = \ln \rho
    \,.
\end{equation}
\end{subequations}
Together, \cref{eq:ig-monge-a} and \cref{eq:ig-monge-b} provide a reformulation of the Monge--Amp\`ere equation that does not appear to have been explored in the literature; see \cite{Gutirrez2016,Figalli2017,Neilan2020} and the many references therein.

The following variational formulation of~\cref{eq:ig-monge} is found by multiplying by smooth test functions and integrating over the domain $\Omega$: find $u \in H^2(\Omega) \cap H^1_g(\Omega)$ and $\psi \in L^\infty(\Omega; \mathbb{R}^{n\times n}_\text{sym})$ such that
\begin{subequations}
\label{eq:MA1}
\begin{align}
\label{eq:MA1-1}
    (\nabla^2 u, w) - (\exp \psi, w)
    &= 
    0, \\
\label{eq:MA1-2}
    (\operatorname{tr} \psi, v)
    &= (\ln \rho, v),
\end{align}
\end{subequations}
for all $v \in H^2(\Omega) \cap H^1_0(\Omega)$ and $w \in L^\infty(\Omega; \mathbb{R}^{n\times n}_\text{sym})$.
This formulation permits finite element discretizations with Argyris elements for $u$ and $v$.
However, due to the limited support for Argyris elements in most software, we suggest a first-order system discretization employing the most widely available types of elements.
To this end, we introduce the optimal transport map $T = \nabla u$ as a slack variable into the system of equations~\cref{eq:MA1}, leading to the following reformulation: find $u \in H^1_g(\Omega)$, $T \in H^1(\Omega,\mathbb{R}^n)$, and $\psi \in L^\infty(\Omega; \mathbb{R}^{n\times n}_\text{sym})$ such that
\begin{subequations}
\label{eq:MA2}
\begin{align}
\label{eq:MA2-0}
    (T, S) - (\nabla u, S)
    &= 
    0, \\
\label{eq:MA2-1}
    (\nabla T, w) - (\exp \psi, w)
    &= 
    0, \\
\label{eq:MA2-2}
    (\operatorname{tr} \psi, v)
    &= (\ln \rho, v),
\end{align}
\end{subequations}
for all $v \in H^1_0(\Omega)$, $S \in H^1(\Omega,\mathbb{R}^n)$, and $w \in L^\infty(\Omega; \mathbb{R}^{n\times n}_\text{sym})$.
We note that this is technically not an LVPP formulation, as there is no proximal loop.
Likewise, its Galerkin finite element discretization should not be considered a proximal Galerkin method.

\begin{figure}[t]
\centering
  \begin{minipage}[c]{0.4\linewidth}
    \includegraphics[clip, trim = 8cm 6cm 5cm 11cm, height=3.5cm]{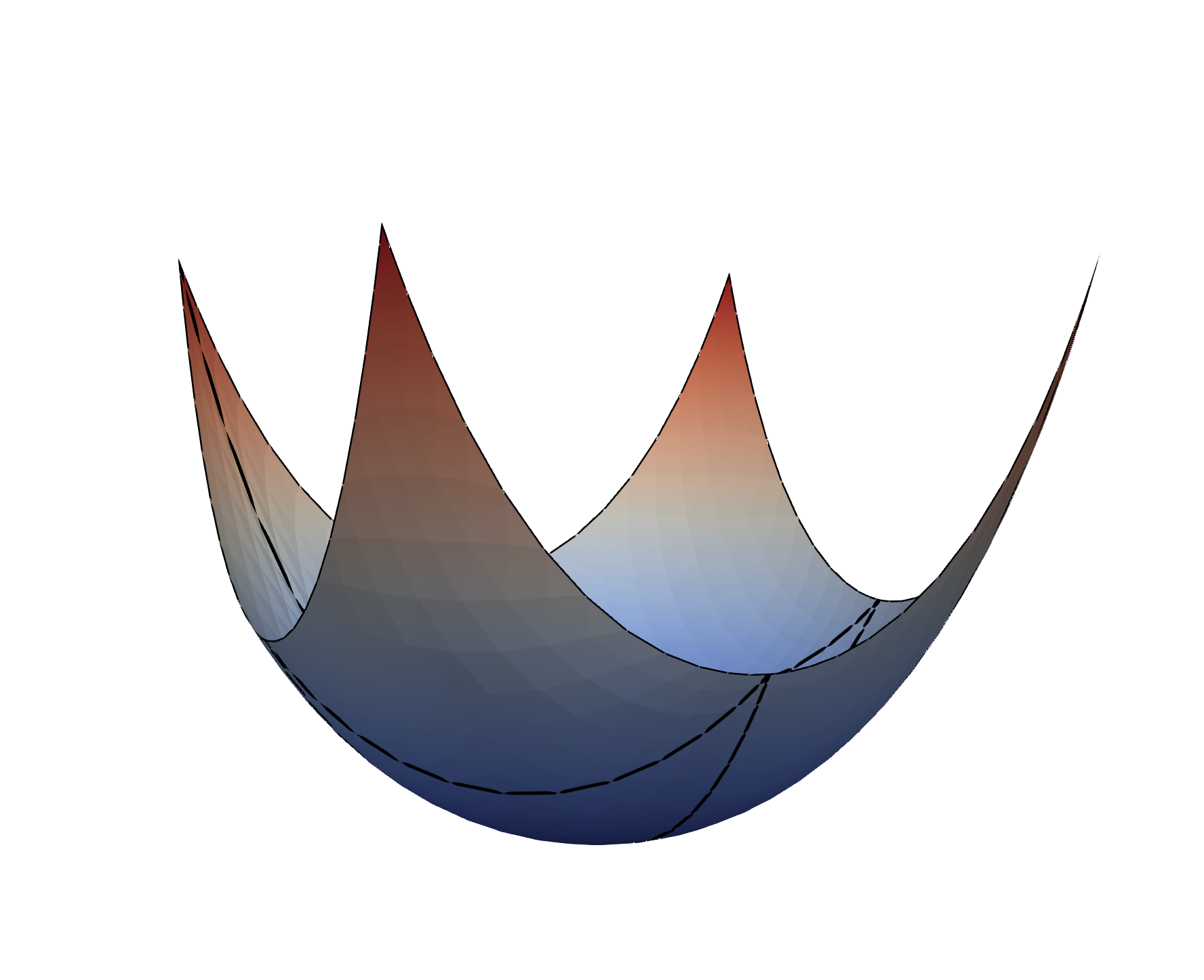}
  \end{minipage}
  \begin{minipage}[c]{0.04\linewidth}
  \includegraphics[height=3cm]{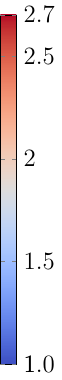}
  \end{minipage}
  ~~~
  \begin{minipage}[c]{0.4\linewidth}
        \pgfplotsset{
    legend entry/.initial=,
    every axis plot post/.code={        \pgfkeysgetvalue{/pgfplots/legend entry}\tempValue
        \ifx\tempValue\empty
            \pgfkeysalso{/pgfplots/forget plot}        \else
            \expandafter\addlegendentry\expandafter{\tempValue}        \fi
    },
}
\newcommand{\equal}{=}

\begin{tikzpicture}

\definecolor{seabornblue}{rgb}{0.2980392156862745, 0.4470588235294118, 0.6901960784313725}
\definecolor{seaborngreen}{rgb}{0.3333333333333333, 0.6588235294117647, 0.40784313725490196}
\definecolor{seabornred}{rgb}{0.7686274509803922, 0.3058823529411765, 0.3215686274509804}
\definecolor{seabornpurple}{rgb}{0.5058823529411764, 0.4470588235294118, 0.6980392156862745}
\definecolor{seabornsand}{rgb}{0.8, 0.7254901960784313, 0.4549019607843137}
\definecolor{seaborncyan}{rgb}{0.39215686274509803, 0.7098039215686275, 0.803921568627451}
\definecolor{seabornorange}{rgb}{0.958, 0.476, 0.206}

\begin{loglogaxis}[
xlabel={Degrees of Freedom},
ylabel={$\|u - u_h\|_{L^2(\Omega)}$},
xmin=100, xmax=10000,
ytick={1e-1, 1e-5, 1e-9, 1e-13},
width=4.5cm,
height=4cm,
]
\addplot [mark=*, very thick, seaborngreen]
table {198 0.10451552404823607
358 0.012027412025493565
566 0.0012939596098955697
822 0.00018992077659682886
1126 2.3025993957722756e-05
1478 3.0786892425448895e-06
1878 3.5627520395916875e-07
2326 4.229649002774819e-08
2822 4.5387319168003806e-09
3366 4.891370832841541e-10
3958 4.923387950569194e-11
4598 4.954193796369616e-12
5286 5.255031324605701e-13
};

\end{loglogaxis}

\end{tikzpicture}

  \end{minipage}
\caption{Example 10. Exponential convergence of a manufactured solution to the Monge--Amp\`ere equation on an eight-element mesh. Left: Computed solution with polynomial degree $p = 14$. Right: Degrees of freedom vs.\ $L^2(\Omega)$-error for polynomial degrees $p=2,3,\ldots,14$.}
\label{fig:MA}
\end{figure}
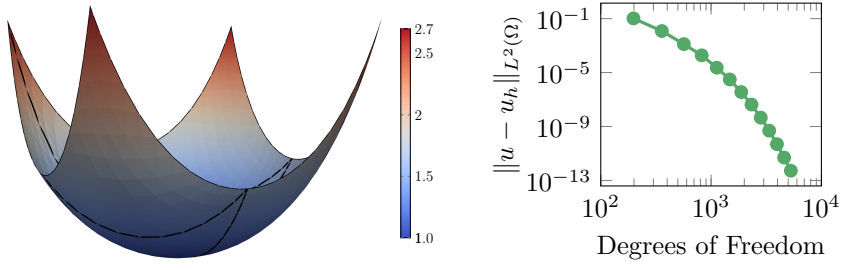

Our experiments indicate that formulation~\cref{eq:MA2} permits high-order accurate, stable discretizations using degree-$p$ continuous elements for $u$, vector-valued degree-$(p+1)$ continuous elements for $T$, and matrix-valued degree-$p$ continuous elements for $\psi$, if $p \geq 2$.
\Cref{fig:MA} showcases exponential convergence of the proposed method to a smooth manufactured solution $u_\text{man}(x) = \exp(|x|^2/2)$ as the polynomial degree $p$ is increased from $2$ to $14$.
In this experiment, we set $\Omega = (-1,1)^2$, $\rho = \operatorname{det} \nabla^2 u_\text{man}$ and $g = u_\text{man}|_{\partial\Omega}$.
Starting from a linearization around $\psi = 0$, the resulting discretizations converged in no more than six Newton iterations on a uniform eight-element triangular mesh.

\section{Conclusion}
\label{sec:Conclusion}

The LVPP algorithm provides a unifying mathematical framework for solving variational problems with pointwise inequality constraints.
This framework is entirely separate from the quadratic penalty, interior point, active set, trust region, and augmented Lagrangian methods, which have been studied extensively over many decades.
LVPP is derived at the level of the underlying function spaces, leaving it agnostic to the discretization strategy preferred for each application.
In turn, the LVPP framework can be used to construct novel numerical methods---including, but not limited to, novel finite element, finite difference, and spectral methods---for a wide class of challenging problems.

We have selected ten example problems and used the LVPP framework to derive state-of-the-art numerical methods for each of them.
The first five examples include free-boundary and free-discontinuity problems, as well as quasi-variational inequalities with pointwise bound constraints.
We find that LVPP delivers high-order, bound-preserving numerical methods for all five of these examples.
The next three examples demonstrate how to apply LVPP to solve problems with gradient constraints, eigenvalue constraints, and multiple independent inequality constraints, respectively.
Lastly, the LVPP framework is used to derive new methods for computing viscosity solutions of first and second-order fully nonlinear PDEs.
In each of these example problems, as well as every other problem we have attempted, LVPP or certain variants, see, e.g., \cite{keith2024analysis,kim2024simple}, have led to structure-preserving numerical methods that exhibit mesh-independence.

Each of the numerical methods in this paper requires further investigation.
In most cases, our derivation of the particular LVPP saddle-point problem is formal and calls for a rigorous mathematical treatment in future work.
Likewise, error analysis needs to be conducted on almost all of the derived methods to guarantee their convergence as the resolution of their discretization increases.
We hope that other applied mathematicians will join us in this task and in applying LVPP to other problems with inequality constraints.

\section*{Code availability} 
In order to facilitate the broader adoption of this work by the community, we have implemented the derived methods across a range of popular open-source software libraries and prepared a suite of scripts to reproduce the numerical experiments.
Scripts to generate all the examples are available at \cite{github:proximalgalerkin} and are archived on Zenodo~\cite{zenodo:proximalgalerkin}. We have implemented LVPP solvers in the finite element packages MFEM \cite{anderson2021mfem,andrej2024high}, Firedrake \cite{FiredrakeUserManual}, FEniCSx \cite{baratta2023}, and Gridap \cite{Badia2020}, and the spectral method package MultivariateOrthogonalPolynomials \cite{MulPoly.jl2024}.

\section*{Acknowledgements}

PEF was funded by the Engineering and Physical Sciences Research Council [grant numbers EP/R029423/1 and EP/W026163/1] and by the Donatio Universitatis Carolinae Chair ``Mathematical modelling of multicomponent systems''.
PEF would like to thank U.~Zerbinati for help with \Cref{sec:ldg}.
BK was supported in part by the U.S.\ Department of Energy Office of Science, Early Career Research Program under Award Number DE-SC0024335. IPAP was funded by the Deutsche Forschungsgemeinschaft (DFG, German Research Foundation) under Germany's Excellence Strategy -- The Berlin Mathematics Research Center MATH+ (EXC-2046/1, project ID: 390685689).
For the purpose of open access, the authors have applied a CC BY public copyright license to any author-accepted manuscript arising from this submission.

\appendix

\section{Linear equality constraints}
\label{app:LinearEqualityConstraints}
This appendix extends the LVPP framework to treat linear equality constraints.
In particular, by replacing $R^\ast$ with $\epsilon R^\ast$ and taking the limit $\epsilon \downarrow 0$, we uncover the standard Lagrange multiplier formulation for constraints of the form $Bu = 0$.
These problems, including notable examples like the Stokes and Maxwell equations, are well-studied in the context of saddle-point theory for linear constraints \cite{boffi2013mixed}, and we therefore omit a numerical example. 
Instead,  we highlight the connection between this classical theory and the LVPP algorithm, demonstrating that the LVPP approach converges in a single iteration when minimizing a quadratic objective function, mirroring the efficiency of Newton's method \cite{Zeidler1986} on quadratic problems.
While the setting is abstract and general, we do not seek maximal generality out of an effort to ease understanding.

Let $V$ and $W$ be real separable Hilbert spaces, with $W$ continuously embedded in $(L^2(\Omega))^m$, and let $B \colon V \to W$ be a bounded, surjective ($\operatorname{im} B = W$) linear operator.
We define
\begin{equation}
\label{eq:LinearEqualityConstraint}
    K = \operatorname{ker} B
    =
    \{ v \in V \mid Bv = 0 \}
    \,,
\end{equation}
and consider optimizing a quadratic functional over $K$,
\begin{equation}
\label{eq:GeneralQuadratic}
    J(v) = \frac12 a(v,v) - F(v)
    \,,
\end{equation}
where $a \colon V \times V \to \mathbb{R}$ is a continuous, symmetric, coercive bilinear form, and $F \in V^\prime$.
In this case, the associated optimality conditions are well-known \cite[Theorem~4.2.1]{boffi2013mixed}, and reduce to the following variational problem: find $u \in V$ and $\lambda \in W$ such that
\begin{subequations}
\label{eq:mixedmethod}
\begin{align}
    a(u,v) + (\lambda,B v) &= F(v), \\
    (B u,w) &= 0,
\end{align}
\end{subequations}
for all $v \in V$ and $w \in W$.
For example, taking $V = (H^1_0(\Omega))^d$, $W = L^2_0(\Omega) := \{ v \in L^2(\Omega) \mid \int_\Omega v \dd x = 0 \}$, $a(u,v) = (\nabla u, \nabla v)$, $Bv = -\nabla \cdot v$ (the negative divergence operator), and $F(v) = (f,v)$ for a forcing function $f \in (L^2(\Omega))^d$, we get the Stokes equations \cite{bernardi2024mathematics} for unknown velocity $u \in (H^1_0(\Omega))^d$ and pressure $p \in L^2_0(\Omega)$:
\begin{subequations}
\begin{align}
    (\nabla u, \nabla v) - (p, \nabla\cdot v) &= (f,v), \\
    (\nabla\cdot u,q) &= 0,
\end{align}
\end{subequations}
for all $v \in (H^1_0(\Omega))^d$ and $q \in L^2_0(\Omega)$.

Upon setting $C = \{ 0 \}$ and $\Omega_d = \Omega$, the linear subspace $K$ in~\cref{eq:LinearEqualityConstraint} seems to fit into the abstract setting~\cref{eq:intro:feasible_set}.
However, there is no Legendre function for $C$ since this convex set has only a single point; i.e., it has an empty interior.
As a remedy, we proceed by considering the closed unit ball centered at $0 \in \mathbb{R}^m$, denoted $S(0,1) = \{ a \in \mathbb{R}^m \mid |x| \leq 1 \}$, with the Hellinger entropy $R\colon S(0,1) \to [-1,0]$, $a \mapsto -\sqrt{1 - |a|^2}$. Outside of $S(0,1)$, we may define $R \equiv +\infty$.
Given $0 < \epsilon < 1$, we define $R_\epsilon = \epsilon R \circ (\epsilon^{-1}\operatorname{id})$, noting that $\operatorname{dom} R_\epsilon = S(0,\epsilon)$.

The first LVPP subproblem that arises from minimizing $J(v)$ in~\cref{eq:GeneralQuadratic} over the feasible set
\[
    K_\epsilon
    =
    \{ v \in V \mid Bv(x) \in S(0,\epsilon) \text{ for almost all } x \in \Omega \}
    \,,
\]
is as follows: for fixed $\psi^0 \in W$, find $u_\epsilon \in V$ and $\psi_\epsilon \in W$ such that
\begin{subequations}
\label{eq:LVPP_epsilon}
\begin{align}
\label{eq:LVPP_epsilon:a}
    \alpha_1 a(u_\epsilon,v) + (\psi_\epsilon,B v) &= \alpha_1 F(v) + (\psi^0,B v), \\
\label{eq:LVPP_epsilon:b}
    (B u_\epsilon,w) - \epsilon(\nabla R^\ast(\psi_\epsilon),w) &= 0,
\end{align}
\end{subequations}
for all $v \in V$ and $q \in W$.
Letting $\epsilon$ pass to zero from above, we derive a saddle-point problem similar to~\cref{eq:mixedmethod}.
In particular, we show below that $u_\epsilon  \rightharpoonup \overline{u}$ weakly in  $V$ and $\psi_\epsilon \rightharpoonup \overline{\psi}$ weakly in $W$, where $\overline{u} \in V$ and $\overline{\psi} \in W$ solve
\begin{subequations}
\label{eq:LVPP_epsilon_limit}
\begin{align}
    \alpha_1 a(\overline{u},v) + (\overline{\psi},B v) &= \alpha_1 F(v) + (\psi^{0},B v), \\
    (B \overline{u},w) &= 0,
\end{align}
\end{subequations}
for all $v \in V$ and $q \in W$. Taking $u = \overline{u}$ and $\lambda = \alpha_1^{-1}(\overline{\psi}-\psi^0)$, we arrive back at \cref{eq:mixedmethod}, which also ensures well-posedness of $\eqref{eq:LVPP_epsilon_limit}$. 
We can then conclude that the latent variable approach always converges in one iteration and is equivalent to the standard Lagrange multiplier mixed method~\cref{eq:mixedmethod} for this class of problems.

It remains to complete the analysis showing that the solution of~\cref{eq:LVPP_epsilon} weakly converges to the solution of~\cref{eq:LVPP_epsilon_limit} as $\epsilon \downarrow 0$.
As we explicitly require using the coercivity and continuity constants to complete the argument, we note that $a\colon V \times V \to \mathbb{R}$ being coercive means that there exists a positive number $\beta > 0$ such that
\(
    \beta \|v\|_V^2
    \leq
    a(v,v)
\)
for all $v \in V$.
Likewise, continuity means that there exists a number $C$ such that
\(
    a(u,v) \leq C \|u\|_V \|v\|_V
\)
for all $u,v \in V$.

Using $u_\epsilon$ as a test function in~\cref{eq:LVPP_epsilon:a}, we arrive at
\begin{equation}\label{eq:apriori}
    \beta \| u_\epsilon \|^2_V + \alpha_1^{-1}(\psi_\epsilon,B u_\epsilon)
    \le
    \| F \|_{V^\prime} \|u_\epsilon\|_{V} + \alpha_1^{-1} \| B^* \psi^0 \|_{V^\prime} \|u_\epsilon\|_{V}
    ,
\end{equation}
where $B^\ast \colon W  \to V^\prime$ denotes the $L^2$-adjoint of $B$, where we have identified $W$ with its dual $W'$.
Next, we recall the identity $\nabla R^\ast(0) = 0$ for the Hellinger entropy~\cref{eq:hellinger}.
Furthermore, note that $R^{\ast}$ is a smooth convex function, and therefore, its gradient is a monotone operator, i.e., $(\nabla R^{\ast}(a^*)-\nabla R^{\ast}(b^*),a^* - b^*)_{\mathbb R^m} \ge 0$ for all pairs $a^*, b^* \in \operatorname{dom} \nabla R^*$.
Consequently, the following inequality arises from~\cref{eq:LVPP_epsilon:a}:
\begin{equation}
    (\psi_\epsilon,B u_\epsilon)
    =
    (\psi_\epsilon, \epsilon\nabla R^\ast(\psi_\epsilon))
    =
    \epsilon(\psi_\epsilon - 0, \nabla R^\ast(\psi_\epsilon) - \nabla R^\ast(0))
    \geq
    0.
\end{equation}
In turn, we find that $u_\epsilon$ is bounded uniformly across $\epsilon>0$,
\begin{equation}
\label{eq:UniformBound_ueps}
    \beta \| u_\epsilon \|_V
    \leq
    \| F \|_{V^\prime} + \alpha_1^{-1} \| B^* \psi^0 \|_{V^\prime}
    \,.
\end{equation}
The uniform bound given above allows us to pass to the limit for a (positive) null sequence $\epsilon_k \downarrow 0$. Even without knowledge of the convergence properties of $\psi_{\epsilon_k}$, we can deduce that the sequence 
\[
\epsilon_{k} \nabla R^\ast(\psi_{\epsilon_k}) \to 0 \text{ in } L^\infty(\Omega)
\]
since the nonlinearity remains bounded in the bounded set $S(0,1)$.
By reflexivity, there exists a subsequence $\{\epsilon_{k_l}\} \subset \{\epsilon_{k}\}$ along which $u_{\epsilon_{k_l}}$ converges weakly in $V$ to some $\overline{u}$.
In particular, we may conclude that 
\(
B \overline{u} = 0
\).
Next, we invoke the closed range theorem and the surjectivity of $B \colon V \to W$ \cite[Theorem~2.20]{brezis2011functional} to establish the existence of a constant $c > 0$, independent of $\epsilon$, such that
\[
    c\|\psi_{\epsilon}\|_W
    \leq
    \|B^\ast\psi_{\epsilon}\|_{V^\prime}
    =
    \sup_{v \in V}
    \frac{| (\psi_{\epsilon}, Bv) |}{\|v\|_V}
        \le
    \alpha_1 \|F\|_{V^\prime}
    + \| B^* \psi^0 \|_{V^\prime} + \alpha_1 C \| u_{\epsilon} \|_{V}
    \,.
\]
This and~\cref{eq:UniformBound_ueps} shows that $\psi_\epsilon$ is bounded uniformly across $\epsilon>0$.
By reflexivity of $W$, there exists a weak accumulation point $\overline{\psi} \in W$ of $\psi_{\epsilon_k}$ for any sequence $\epsilon_k \downarrow 0$.
This implies that $(\psi_{\epsilon} - \psi^0)/\alpha_1$ will (weakly) converge along an appropriate subsequence to $(\overline{\psi} - \psi^0)/\alpha_1$. As observed above, we can take $\overline{u} = u$ and $\lambda = (\overline{\psi} - \psi^0)/\alpha_1$ in~\cref{eq:mixedmethod}. Since $u$ and $\lambda$ are unique it follows from the Fr\'echet--Urysohn property that
the full sequence converges weakly, i.e.,
$u_{\epsilon} \rightharpoonup u$ in $V$ and $(\psi_{\epsilon} - \psi^0)/\alpha_1 \rightharpoonup \lambda$ in $W$ as $\epsilon \downarrow 0$, as necessary.

\bibliographystyle{siamplain}
\bibliography{ref}

\begin{thebibliography}{100}

\bibitem{Adam2019}
{\sc L.~Adam, M.~Hinterm\"{u}ller, D.~Peschka, and T.~M. Surowiec}, {\em
  Optimization of a multiphysics problem in semiconductor laser design}, SIAM
  Journal on Applied Mathematics, 79 (2019), p.~257–283,
  \url{https://doi.org/10.1137/18m1179183}.

\bibitem{adam2019semismooth}
{\sc L.~Adam, M.~Hintermüller, and T.~M. Surowiec}, {\em A semismooth {Newton}
  method with analytical path-following for the {$H^1$}-projection onto the
  {Gibbs} simplex}, IMA Journal of Numerical Analysis, 39 (2018),
  pp.~1276--1295, \url{https://doi.org/10.1093/imanum/dry034}.

\bibitem{agarwal2003}
{\sc R.~Agarwal and D.~O'Regan}, {\em Nonlinear generalized quasi-variational
  inequalities: {A} fixed point approach}, Mathematical Inequalties \&
  Applications, 6 (2003), pp.~133--143,
  \url{https://doi.org/10.7153/mia-06-13}.

\bibitem{allgower1986mesh}
{\sc E.~L. Allgower, K.~B{\"o}hmer, F.~Potra, and W.~Rheinboldt}, {\em A
  mesh-independence principle for operator equations and their
  discretizations}, SIAM Journal on Numerical Analysis, 23 (1986),
  pp.~160--169, \url{https://doi.org/10.1137/0723011}.

\bibitem{alnaes2012}
{\sc M.~S. Aln\ae{}s, A.~Logg, K.~B. {\O}lgaard, M.~E. Rognes, and G.~N.
  Wells}, {\em {U}nified {F}orm {L}anguage: a domain-specific language for weak
  formulations of partial differential equations}, ACM Transactions on
  Mathematical Software, 40 (2014), pp.~9:1--9:37,
  \url{https://doi.org/10.1145/2566630}.

\bibitem{alphonse2024}
{\sc A.~Alphonse, C.~Christof, M.~Hinterm{\"u}ller, and I.~P.~A. Papadopoulos},
  {\em A globalized inexact semismooth {N}ewton method for nonsmooth
  fixed-point equations involving variational inequalities}, arXiv preprint
  arXiv:2409.19637,  (2024), \url{https://doi.org/10.48550/arXiv.2409.19637}.

\bibitem{alphonse2019}
{\sc A.~Alphonse, M.~Hinterm{\"u}ller, and C.~N. Rautenberg}, {\em Directional
  differentiability for elliptic quasi-variational inequalities of obstacle
  type}, Calculus of Variations and Partial Differential Equations, 58 (2019),
  p.~39, \url{https://doi.org/10.1007/s00526-018-1473-0}.

\bibitem{alphonse2019b}
{\sc A.~Alphonse, M.~Hinterm{\"u}ller, and C.~N. Rautenberg}, {\em Recent
  trends and views on elliptic quasi-variational inequalities}, Springer, 2019,
  \url{https://doi.org/10.1007/978-3-030-33116-0_1}.

\bibitem{alvarez2004hessian}
{\sc F.~Alvarez, J.~Bolte, and O.~Brahic}, {\em Hessian {R}iemannian gradient
  flows in convex programming}, SIAM journal on control and optimization, 43
  (2004), pp.~477--501.

\bibitem{Amari2016}
{\sc S.-i. Amari}, {\em Information Geometry and Its Applications}, Springer
  Japan, 2016, \url{https://doi.org/10.1007/978-4-431-55978-8}.

\bibitem{Ambrosio1990}
{\sc L.~Ambrosio and V.~M. Tortorelli}, {\em Approximation of functional
  depending on jumps by elliptic functional via t‐convergence},
  Communications on Pure and Applied Mathematics, 43 (1990), p.~999–1036,
  \url{https://doi.org/10.1002/cpa.3160430805}.

\bibitem{anderson2021mfem}
{\sc R.~Anderson, J.~Andrej, A.~Barker, J.~Bramwell, J.-S. Camier, J.~Cerveny,
  V.~Dobrev, Y.~Dudouit, A.~Fisher, T.~Kolev, W.~Pazner, M.~Stowell, V.~Tomov,
  I.~Akkerman, J.~Dahm, D.~Medina, and S.~Zampini}, {\em {MFEM}: A modular
  finite element methods library}, Computers \& Mathematics with Applications,
  81 (2021), pp.~42--74, \url{https://doi.org/10.1016/j.camwa.2020.06.009}.

\bibitem{andrej2024high}
{\sc J.~Andrej, N.~Atallah, J.-P. B{\"a}cker, J.-S. Camier, D.~Copeland,
  V.~Dobrev, Y.~Dudouit, T.~Duswald, B.~Keith, D.~Kim, et~al.}, {\em
  High-performance finite elements with {MFEM}}, The International Journal of
  High Performance Computing Applications, 38 (2024), pp.~447--467.

\bibitem{Antil_2022}
{\sc H.~Antil, R.~Arndt, C.~N. Rautenberg, and D.~Verma}, {\em Nondiffusive
  variational problems with distributional and weak gradient constraints},
  Advances in Nonlinear Analysis, 11 (2022), pp.~1466--1495,
  \url{https://doi.org/doi:10.1515/anona-2022-0227}.

\bibitem{Araki1970}
{\sc H.~Araki and E.~H. Lieb}, {\em Entropy inequalities}, Communications in
  Mathematical Physics, 18 (1970), p.~160–170,
  \url{https://doi.org/10.1007/bf01646092}.

\bibitem{arnold2018}
{\sc D.~N. Arnold}, {\em Finite element exterior calculus}, Society for
  Industrial and Applied Mathematics, 2018,
  \url{https://doi.org/10.1137/1.9781611975543}.

\bibitem{aubin2007}
{\sc J.~Aubin}, {\em Mathematical Methods of Game and Economic Theory}, Dover
  books on mathematics, Dover Publications, 2007.

\bibitem{Badia2020}
{\sc S.~Badia and F.~Verdugo}, {\em Gridap: An extensible finite element
  toolbox in {J}ulia}, Journal of Open Source Software, 5 (2020), p.~2520,
  \url{https://doi.org/10.21105/joss.02520}.

\bibitem{baiocchi1984}
{\sc C.~Baiocchi and A.~Capelo}, {\em Variational and Quasivariational
  Inequalities}, A Wiley-Interscience Publication, John Wiley \& Sons, Inc.,
  New York, 1984.

\bibitem{petsc-user-ref}
{\sc S.~Balay, S.~Abhyankar, M.~F. Adams, S.~Benson, J.~Brown, P.~Brune,
  et~al.}, {\em {PETSc/TAO} users manual}, Tech. Report ANL-21/39 - Revision
  3.21, Argonne National Laboratory, 2024,
  \url{https://doi.org/10.2172/2205494}.

\bibitem{ball2017}
{\sc J.~M. Ball}, {\em Mathematics and liquid crystals}, Molecular Crystals and
  Liquid Crystals, 647 (2017), pp.~1--27,
  \url{https://doi.org/10.1080/15421406.2017.1289425}.

\bibitem{ball2010}
{\sc J.~M. Ball and A.~Majumdar}, {\em Nematic liquid crystals: from
  {Maier--Saupe} to a continuum theory}, Molecular Crystals and Liquid
  Crystals, 525 (2010), pp.~1--11,
  \url{https://doi.org/10.1080/15421401003795555}.

\bibitem{banz2015}
{\sc L.~Banz and A.~Schr{\"o}der}, {\em Biorthogonal basis functions in
  $hp$-adaptive {FEM} for elliptic obstacle problems}, Computers \& Mathematics
  with Applications, 70 (2015), pp.~1721--1742,
  \url{https://doi.org/10.1016/j.camwa.2015.07.010}.

\bibitem{baratta2023}
{\sc I.~A. Baratta, J.~P. Dean, J.~S. Dokken, M.~Habera, J.~S. Hale, C.~N.
  Richardson, M.~E. Rognes, M.~W. Scroggs, N.~Sime, and G.~N. Wells}, {\em
  {DOLFINx}: {T}he next generation {FEniCS} problem solving environment}.
\newblock Unpublished, 2023, \url{https://doi.org/10.5281/zenodo.10447666}.

\bibitem{bauschke1997legendre}
{\sc H.~H. Bauschke and J.~M. Borwein}, {\em Legendre functions and the method
  of random {B}regman projections}, Journal of convex analysis, 4 (1997),
  pp.~27--67.

\bibitem{beck2017existence}
{\sc L.~Beck, M.~Bul{\'\i}{\v{c}}ek, J.~M{\'a}lek, and E.~S{\"u}li}, {\em On
  the existence of integrable solutions to nonlinear elliptic systems and
  variational problems with linear growth}, Archive for Rational Mechanics and
  Analysis, 225 (2017), pp.~717--769,
  \url{https://doi.org/10.1007/s00205-017-1113-4}.

\bibitem{benson2006}
{\sc S.~J. Benson and T.~S. Munson}, {\em Flexible complementarity solvers for
  large-scale applications}, Optimization Methods and Software, 21 (2006),
  pp.~155--168, \url{https://doi.org/10.1080/10556780500065382}.

\bibitem{bensoussan1975}
{\sc A.~Bensoussan and J.~L. Lions}, {\em Optimal impulse and continuous
  control: Method of nonlinear quasi-variational inequalities}, Trudy
  Matematicheskogo Instituta imeni V.A. Steklova, 134 (1975), pp.~5--22.

\bibitem{bensoussan1984}
{\sc A.~Bensoussan and J.-L. Lions}, {\em Impulse Control and Quasivariational
  Inequalities}, Gauthier-Villars, Montrouge; Heyden \& Son, Inc.,
  Philadelphia, PA, 1984.

\bibitem{bernardi2024mathematics}
{\sc C.~Bernardi, V.~Girault, F.~Hecht, P.-A. Raviart, and B.~Rivi{\`e}re},
  {\em Mathematics and Finite Element Discretizations of Incompressible
  Navier-Stokes Flows}, SIAM, 2024,
  \url{https://doi.org/10.1137/1.9781611978124}.

\bibitem{bernstein1993}
{\sc D.~Bernstein and W.~So}, {\em Some explicit formulas for the matrix
  exponential}, IEEE Transactions on Automatic Control, 38 (1993),
  pp.~1228--1232, \url{https://doi.org/10.1109/9.233156}.

\bibitem{biegler2010nonlinear}
{\sc L.~Biegler}, {\em Nonlinear Programming: Concepts, Algorithms, and
  Applications to Chemical Processes}, MOS-SIAM Series on Optimization, Society
  for Industrial and Applied Mathematics (SIAM, 3600 Market Street, Floor 6,
  Philadelphia, PA 19104), 2010, \url{https://doi.org/10.1137/1.9780898719383}.

\bibitem{boffi2013mixed}
{\sc D.~Boffi, F.~Brezzi, and M.~Fortin}, {\em Mixed finite element methods and
  applications}, vol.~44, Springer, 2013,
  \url{https://doi.org/10.1007/978-3-540-78319-0}.

\bibitem{Bourdin2000}
{\sc B.~Bourdin, G.~Francfort, and J.-J. Marigo}, {\em Numerical experiments in
  revisited brittle fracture}, Journal of the Mechanics and Physics of Solids,
  48 (2000), p.~797–826, \url{https://doi.org/10.1016/s0022-5096(99)00028-9}.

\bibitem{BREGMAN1967200}
{\sc L.~Bregman}, {\em The relaxation method of finding the common point of
  convex sets and its application to the solution of problems in convex
  programming}, USSR Computational Mathematics and Mathematical Physics, 7
  (1967), pp.~200--217, \url{https://doi.org/10.1016/0041-5553(67)90040-7}.

\bibitem{brezis2011functional}
{\sc H.~Brezis and H.~Br{\'e}zis}, {\em Functional analysis, {S}obolev spaces
  and partial differential equations}, vol.~2, Springer, 2011.

\bibitem{Brezis1971}
{\sc H.~Brezis and M.~Sibony}, {\em Equivalence de deux inéquations
  variationnelles et applications}, Archive for Rational Mechanics and
  Analysis, 41 (1971), p.~254–265, \url{https://doi.org/10.1007/bf00250529}.

\bibitem{bruggemann2023}
{\sc J.~A. Br{\"u}ggemann}, {\em A class of elliptic obstacle-type
  quasi-variational inequalities: {T}heory and solution methods},
  Humboldt-Universit{\"a}t zu Berlin,  (2023),
  \url{https://doi.org/10.18452/27750}.

\bibitem{bueler2023}
{\sc E.~Bueler and P.~E. Farrell}, {\em A full approximation scheme multilevel
  method for nonlinear variational inequalities}, SIAM Journal on Scientific
  Computing, 46 (2024), pp.~A2421--A2444,
  \url{https://doi.org/10.1137/23M1594200}.

\bibitem{bulivcek2014}
{\sc M.~Bul{\'\i}{\v{c}}ek, J.~M{\'a}lek, K.~R. Rajagopal, and E.~S{\"u}li},
  {\em On elastic solids with limiting small strain: modelling and analysis},
  EMS Surveys in Mathematical Sciences, 1 (2014), pp.~283--332,
  \url{https://doi.org/10.4171/emss/7}.

\bibitem{bulivcek2015}
{\sc M.~Bul{\'\i}{\v{c}}ek, J.~M{\'a}lek, and E.~S{\"u}li}, {\em Analysis and
  approximation of a strain-limiting nonlinear elastic model}, Mathematics and
  Mechanics of Solids, 20 (2015), pp.~92--118,
  \url{https://doi.org/10.1177/1081286514543601}.

\bibitem{Burke2010}
{\sc S.~Burke, C.~Ortner, and E.~S\"uli}, {\em An adaptive finite element
  approximation of a variational model of brittle fracture}, SIAM Journal on
  Numerical Analysis, 48 (2010), pp.~980--1012,
  \url{https://doi.org/10.1137/080741033}.

\bibitem{burns2020}
{\sc K.~J. Burns, G.~M. Vasil, J.~S. Oishi, D.~Lecoanet, and B.~P. Brown}, {\em
  Dedalus: {A} flexible framework for numerical simulations with spectral
  methods}, Physical Review Research, 2 (2020), p.~023068,
  \url{https://doi.org/10.1103/PhysRevResearch.2.023068}.

\bibitem{Caffarelli1990}
{\sc L.~A. Caffarelli}, {\em Interior {$W^{2,p}$} estimates for solutions of
  the {M}onge--{A}mp\`ere equation}, The Annals of Mathematics, 131 (1990),
  p.~135, \url{https://doi.org/10.2307/1971510}.

\bibitem{Casas1993}
{\sc E.~Casas and L.~A. Fernández}, {\em Optimal control of semilinear
  elliptic equations with pointwise constraints on the gradient of the state},
  Applied Mathematics \& Optimization, 27 (1993), p.~35–56,
  \url{https://doi.org/10.1007/bf01182597}.

\bibitem{chen1993convergence}
{\sc G.~Chen and M.~Teboulle}, {\em Convergence analysis of a proximal-like
  minimization algorithm using {B}regman functions}, SIAM Journal on
  Optimization, 3 (1993), pp.~538--543, \url{https://doi.org/10.1137/0803026}.

\bibitem{cheng1997}
{\sc H.-W. Cheng and S.~S.-T. Yau}, {\em More explicit formulas for the matrix
  exponential}, Linear Algebra and its Applications, 262 (1997), pp.~131--163,
  \url{https://doi.org/10.1016/s0024-3795(97)80028-6}.

\bibitem{Clarke1990}
{\sc F.~H. Clarke}, {\em Optimization and Nonsmooth Analysis}, Society for
  Industrial and Applied Mathematics, Jan. 1990,
  \url{https://doi.org/10.1137/1.9781611971309}.

\bibitem{Crandall1983}
{\sc M.~G. Crandall and P.-L. Lions}, {\em Viscosity solutions of
  {Hamilton--Jacobi} equations}, Transactions of the American Mathematical
  Society, 277 (1983), p.~1–42,
  \url{https://doi.org/10.1090/s0002-9947-1983-0690039-8}.

\bibitem{Dacorogna1999-mj}
{\sc B.~Dacorogna and P.~Marcellini}, {\em Implicit partial differential
  equations}, Progress in Nonlinear Differential Equations and Their
  Applications, Birkh{\"a}user, Cambridge, MA, 1999~ed., Aug. 1999,
  \url{https://doi.org/10.1007/978-1-4612-1562-2}.

\bibitem{DalMaso1993}
{\sc G.~Dal~Maso}, {\em An Introduction to $\Gamma$-Convergence},
  Birkh\"{a}user Boston, 1993, \url{https://doi.org/10.1007/978-1-4612-0327-8}.

\bibitem{degennes1995}
{\sc P.~G. de~Gennes and J.~Prost}, {\em The physics of liquid crystals},
  International Series of Monographs on Physics, Oxford University Press,
  2~ed., 1995.

\bibitem{dean2024}
{\sc J.~P. Dean}, {\em Mathematical and computational aspects of solving
  mixed-domain problems using the finite element method}, PhD thesis,
  University of Cambridge, Cambridge, UK, 2024,
  \url{https://doi.org/10.17863/CAM.108292}.

\bibitem{Dijkstra1959}
{\sc E.~W. Dijkstra}, {\em A note on two problems in connexion with graphs},
  Numerische Mathematik, 1 (1959), p.~269–271,
  \url{https://doi.org/10.1007/bf01386390}.

\bibitem{zenodo:proximalgalerkin}
{\sc J.~S. Dokken, P.~E. Farrell, B.~Keith, I.~P. Papadopoulos, and T.~M.
  Surowiec}, {\em The latent variable proximal point algorithm for variational
  problems with inequality constraints}, May 2025,
  \url{https://doi.org/10.5281/zenodo.15554771}.

\bibitem{github:proximalgalerkin}
{\sc J.~S. Dokken, P.~E. Farrell, B.~Keith, I.~P.~A. Papadopoulos, and T.~M.
  Surowiec}, {\em {ProximalGalerkin}}, 2025,
  \url{https://github.com/METHODS-Group/ProximalGalerkin}.

\bibitem{Elliott1989}
{\sc C.~Elliott, D.~French, and F.~Milner}, {\em A second order splitting
  method for the {C}ahn--{H}illiard equation.}, Numerische Mathematik, 54
  (1989), pp.~575--590, \url{https://doi.org/10.1007/BF01396363}.

\bibitem{Ern2021}
{\sc A.~Ern and J.-L. Guermond}, {\em Finite Elements II: Galerkin
  Approximation, Elliptic and Mixed PDEs}, Springer International Publishing,
  Cham, 2021, \url{https://doi.org/10.1007/978-3-030-56923-5_51}.

\bibitem{Evans2010}
{\sc L.~Evans}, {\em Partial Differential Equations}, American Mathematical
  Society, Mar. 2010, \url{https://doi.org/10.1090/gsm/019}.

\bibitem{Fanchi2005-lp}
{\sc J.~R. Fanchi}, {\em Principles of applied reservoir simulation}, Gulf
  Publishing, Oxford, England, 3~ed., Dec. 2005,
  \url{https://doi.org/10.1016/B978-0-7506-7933-6.X5000-4}.

\bibitem{farrell2018a}
{\sc P.~E. Farrell, M.~Croci, and T.~M. Surowiec}, {\em Deflation for
  semismooth equations}, Optimization Methods and Software, 35 (2019),
  pp.~1248--1271, \url{https://doi.org/10.1080/10556788.2019.1613655}.

\bibitem{fichera1963}
{\sc G.~Fichera}, {\em Sul problema elastostatico di {Signorini} con ambigue
  condizioni al contorno}, Atti della Accademia Nazionale dei Lincei, Classe di
  Scienze Fisiche, Matematiche e Naturali, 34 (1963), pp.~138--142.

\bibitem{fife2000models}
{\sc P.~C. Fife}, {\em Models for phase separation and their mathematics},
  Electronic Journal of Differential Equations, 2000 (2000), pp.~1--26.

\bibitem{Figalli2017}
{\sc A.~Figalli}, {\em The Monge–Ampère Equation and Its Applications}, EMS
  Press, Jan. 2017, \url{https://doi.org/10.4171/170},
  \url{http://dx.doi.org/10.4171/170}.

\bibitem{Francfort1998}
{\sc G.~Francfort and J.-J. Marigo}, {\em Revisiting brittle fracture as an
  energy minimization problem}, Journal of the Mechanics and Physics of Solids,
  46 (1998), p.~1319–1342,
  \url{https://doi.org/10.1016/s0022-5096(98)00034-9}.

\bibitem{fu2024locallyconservativeproximalgalerkinmethod}
{\sc G.~Fu, B.~Keith, and R.~Masri}, {\em A locally-conservative proximal
  {G}alerkin method for pointwise bound constraints}, arXiv preprint
  arXiv:2412.21039,  (2024).

\bibitem{Garcke1998}
{\sc H.~Garcke, B.~Nestler, and B.~Stoth}, {\em On anisotropic order parameter
  models for multi-phase systems and their sharp interface limits}, Physica D:
  Nonlinear Phenomena, 115 (1998), p.~87–108,
  \url{https://doi.org/10.1016/s0167-2789(97)00227-3}.

\bibitem{garcke1999multiphase}
{\sc H.~Garcke, B.~Nestler, and B.~Stoth}, {\em A multiphase field concept:
  {N}umerical simulations of moving phase boundaries and multiple junctions},
  SIAM Journal on Applied Mathematics, 60 (1999), pp.~295--315,
  \url{https://doi.org/10.1137/S0036139998334895}.

\bibitem{Gerasimov2020}
{\sc T.~Gerasimov, U.~Römer, J.~Vondřejc, H.~G. Matthies, and L.~{de
  Lorenzis}}, {\em Stochastic phase-field modeling of brittle fracture:
  computing multiple crack patterns and their probabilities}, Computer Methods
  in Applied Mechanics and Engineering, 372 (2020), p.~113353,
  \url{https://doi.org/10.1016/j.cma.2020.113353}.

\bibitem{gilbarg1977elliptic}
{\sc D.~Gilbarg and N.~S. Trudinger}, {\em Elliptic partial differential
  equations of second order}, vol.~224, Springer, 1977.

\bibitem{Glowinski2014-ii}
{\sc R.~Glowinski, J.-L. Lions, and R.~Tr\'emoli\`eres}, {\em Numerical
  analysis of variational inequalities}, Studies in Mathematics and Its
  Applications, North-Holland, May 2014.

\bibitem{gould2003}
{\sc N.~I.~M. Gould, D.~Orban, and P.~L. Toint}, {\em {GALAHAD}, a library of
  thread-safe {F}ortran 90 packages for large-scale nonlinear optimization},
  ACM Transactions on Mathematical Software, 29 (2003), p.~353–372,
  \url{https://doi.org/10.1145/962437.962438}.

\bibitem{gustafsson2017on}
{\sc T.~Gustafsson, R.~Stenberg, and J.~Videman}, {\em On finite element
  formulations for the obstacle problem--mixed and stabilised methods},
  Computational Methods in Applied Mathematics, 17 (2017), pp.~413--429,
  \url{https://doi.org/10.1515/cmam-2017-0011}.

\bibitem{Gutirrez2016}
{\sc C.~E. Gutiérrez}, {\em The Monge-Ampère Equation}, Springer
  International Publishing, 2016,
  \url{https://doi.org/10.1007/978-3-319-43374-5}.

\bibitem{haber2018notes}
{\sc H.~E. Haber}, {\em Notes on the matrix exponential and logarithm}, Santa
  Cruz Institute for Particle Physics, University of California: Santa Cruz,
  CA, USA,  (2018).

\bibitem{hackbusch1983}
{\sc W.~Hackbusch and H.~D. Mittelmann}, {\em On multi-grid methods for
  variational inequalities}, Numerische Mathematik, 42 (1983), pp.~65--76,
  \url{https://doi.org/10.1007/BF01400918}.

\bibitem{FiredrakeUserManual}
{\sc D.~A. Ham, P.~H.~J. Kelly, L.~Mitchell, C.~J. Cotter, R.~C. Kirby,
  K.~Sagiyama, N.~Bouziani, et~al.}, {\em Firedrake User Manual}, Imperial
  College London and University of Oxford and Baylor University and University
  of Washington, first edition~ed., 5 2023,
  \url{https://doi.org/10.25561/104839}.

\bibitem{hintermueller2002}
{\sc M.~Hinterm\"uller, K.~Ito, and K.~Kunisch}, {\em The primal-dual active
  set strategy as a semismooth {N}ewton method}, SIAM Journal on Optimization,
  13 (2002), pp.~865--888, \url{https://doi.org/10.1137/S1052623401383558}.

\bibitem{hintermuller2006feasible}
{\sc M.~Hinterm{\"u}ller and K.~Kunisch}, {\em Feasible and noninterior
  path-following in constrained minimization with low multiplier regularity},
  SIAM Journal on Control and Optimization, 45 (2006), pp.~1198--1221,
  \url{https://doi.org/10.1137/050637480}.

\bibitem{Hintermller2010}
{\sc M.~Hinterm\"{u}ller and K.~Kunisch}, {\em {PDE}-constrained optimization
  subject to pointwise constraints on the control, the state, and its
  derivative}, SIAM Journal on Optimization, 20 (2010), p.~1133–1156,
  \url{https://doi.org/10.1137/080737265}.

\bibitem{Hintermller2012}
{\sc M.~Hinterm\"{u}ller and C.~N. Rautenberg}, {\em A sequential minimization
  technique for elliptic quasi-variational inequalities with gradient
  constraints}, SIAM Journal on Optimization, 22 (2012), p.~1224–1257,
  \url{https://doi.org/10.1137/110837048}.

\bibitem{hintermueller2006b}
{\sc M.~Hintermüller and K.~Kunisch}, {\em Path-following methods for a class
  of constrained minimization problems in function space}, SIAM Journal on
  Optimization, 17 (2006), pp.~159--187,
  \url{https://doi.org/10.1137/040611598}.

\bibitem{hoppe1987}
{\sc R.~H. Hoppe}, {\em Multigrid algorithms for variational inequalities},
  SIAM Journal on Numerical Analysis, 24 (1987), pp.~1046--1065,
  \url{https://doi.org/10.1137/0724069}.

\bibitem{ishii2005limits}
{\sc H.~Ishii and P.~Loreti}, {\em Limits of solutions of $p$-{L}aplace
  equations as $p$ goes to infinity and related variational problems}, SIAM
  Journal on Mathematical Analysis, 37 (2005), pp.~411--437,
  \url{https://doi.org/10.1137/S0036141004432827}.

\bibitem{jaillet1990}
{\sc P.~Jaillet, D.~Lamberton, and B.~Lapeyre}, {\em Variational inequalities
  and the pricing of {American} options}, Acta Applicandae Mathematicae, 21
  (1990), pp.~263--289, \url{https://doi.org/10.1007/bf00047211}.

\bibitem{kanzow2019}
{\sc C.~Kanzow and D.~Steck}, {\em Quasi-variational inequalities in {B}anach
  spaces: {T}heory and augmented {L}agrangian methods}, SIAM Journal on
  Optimization, 29 (2019), pp.~3174--3200,
  \url{https://doi.org/10.1137/18M1230475}.

\bibitem{karper2009unified}
{\sc T.~Karper, K.-A. Mardal, and R.~Winther}, {\em Unified finite element
  discretizations of coupled {D}arcy--{S}tokes flow}, Numerical Methods for
  Partial Differential Equations: An International Journal, 25 (2009),
  pp.~311--326.

\bibitem{keith2024analysis}
{\sc B.~Keith, D.~Kim, B.~S. Lazarov, and T.~M. Surowiec}, {\em Analysis of the
  {SiMPL} method for density-based topology optimization}, arXiv preprint
  arXiv:2409.19341,  (2024).
\newblock To appear in SIAM Journal on Optimization.

\bibitem{keith2023proximal}
{\sc B.~Keith and T.~M. Surowiec}, {\em Proximal {G}alerkin: {A}
  structure-preserving finite element method for pointwise bound constraints},
  Foundations of Computational Mathematics,  (2024),
  \url{https://doi.org/10.1007/s10208-024-09681-8}.

\bibitem{NKikuchi_JTOden_1988}
{\sc N.~Kikuchi and J.~T. Oden}, {\em Contact Problems in Elasticity}, Society
  for Industrial and Applied Mathematics, 1988,
  \url{https://doi.org/10.1137/1.9781611970845}.

\bibitem{kim2024simple}
{\sc D.~Kim, B.~S. Lazarov, T.~M. Surowiec, and B.~Keith}, {\em A simple
  introduction to the {SiMPL} method for density-based topology optimization},
  arXiv preprint arXiv:2411.19421,  (2024).

\bibitem{DKinderlehrer_GStampacchia_2000}
{\sc D.~Kinderlehrer and G.~Stampacchia}, {\em An Introduction to Variational
  Inequalities and Their Applications}, Society for Industrial and Applied
  Mathematics, 2000, \url{https://doi.org/10.1137/1.9780898719451}.

\bibitem{Kirby2010}
{\sc R.~C. Kirby}, {\em Fast simplicial finite element algorithms using
  {Bernstein} polynomials}, Numerische Mathematik, 117 (2010), p.~631–652,
  \url{https://doi.org/10.1007/s00211-010-0327-2}.

\bibitem{kornhuber1994}
{\sc R.~Kornhuber}, {\em Monotone multigrid methods for elliptic variational
  inequalities {I}}, Numerische Mathematik, 69 (1994), pp.~167--184,
  \url{https://doi.org/10.1007/BF03325426}.

\bibitem{kralj2014}
{\sc S.~Kralj and A.~Majumdar}, {\em Order reconstruction patterns in nematic
  liquid crystal wells}, Proceedings of the Royal Society A, 470 (2014),
  p.~20140276, \url{https://doi.org/10.1098/rspa.2014.0276}.

\bibitem{Lions1982-jx}
{\sc P.-L. Lions}, {\em Generalized solutions of {Hamilton-Jacobi} equations},
  Pitman Publishing, Harlow, England, 1982.

\bibitem{lions1986}
{\sc P.-L. Lions and B.~Perthame}, {\em Quasi-variational inequalities and
  ergodic impulse control}, SIAM Journal on Control and Optimization, 24
  (1986), pp.~604--615, \url{https://doi.org/10.1137/0324036}.

\bibitem{Mackenroth1986}
{\sc U.~Mackenroth}, {\em On some elliptic optimal control problems with state
  constraints}, Optimization, 17 (1986), p.~595–607,
  \url{https://doi.org/10.1080/02331938608843174}.

\bibitem{majumdar2010}
{\sc A.~Majumdar}, {\em Equilibrium order parameters of nematic liquid crystals
  in the {Landau--de Gennes} theory}, European Journal of Applied Mathematics,
  21 (2010), pp.~181--203, \url{https://doi.org/10.1017/s0956792509990210}.

\bibitem{Modica1987a}
{\sc L.~Modica}, {\em The gradient theory of phase transitions and the minimal
  interface criterion}, Archive for Rational Mechanics and Analysis, 98 (1987),
  p.~123–142, \url{https://doi.org/10.1007/bf00251230}.

\bibitem{Modica1987b}
{\sc L.~Modica}, {\em Gradient theory of phase transitions with boundary
  contact energy}, Annales de l’Institut Henri Poincaré C, Analyse non
  linéaire, 4 (1987), p.~487–512,
  \url{https://doi.org/10.1016/s0294-1449(16)30360-2}.

\bibitem{mosco1976}
{\sc U.~Mosco}, {\em Implicit variational problems and quasi variational
  inequalities}, in Nonlinear operators and the calculus of variations
  ({S}ummer {S}chool, {U}niv. {L}ibre {B}ruxelles, {B}russels, 1975), Lecture
  Notes in Math., Vol. 543, Springer, Berlin, 1976, pp.~83--156,
  \url{https://doi.org/10.1007/BFb0079943}.

\bibitem{Mumford1989}
{\sc D.~Mumford and J.~Shah}, {\em Optimal approximations by piecewise smooth
  functions and associated variational problems}, Communications on Pure and
  Applied Mathematics, 42 (1989), p.~577–685,
  \url{https://doi.org/10.1002/cpa.3160420503}.

\bibitem{munson2001}
{\sc T.~S. Munson, F.~Facchinei, M.~C. Ferris, A.~Fischer, and C.~Kanzow}, {\em
  The semismooth algorithm for large scale complementarity problems}, INFORMS
  Journal on Computing, 13 (2001), pp.~294--311,
  \url{https://doi.org/10.1287/ijoc.13.4.294.9734}.

\bibitem{Neilan2020}
{\sc M.~Neilan, A.~J. Salgado, and W.~Zhang}, {\em The Monge–Ampère
  equation}, Elsevier, 2020, p.~105–219,
  \url{https://doi.org/10.1016/bs.hna.2019.05.003}.

\bibitem{Neumann1927}
{\sc J.~v. Neumann}, {\em Thermodynamik quantenmechanischer {G}esamtheiten},
  Nachrichten von der Gesellschaft der Wissenschaften zu Göttingen,
  Mathematisch-Physikalische Klasse, 1927 (1927), pp.~273--291,
  \url{http://eudml.org/doc/59231}.

\bibitem{MulPoly.jl2024}
{\sc S.~Olver et~al.}, {\em {MultivariateOrthogonalPolynomials.jl}}, 2024,
  \url{https://github.com/JuliaApproximation/MultivariateOrthogonalPolynomials.jl}.

\bibitem{Olver2013}
{\sc S.~Olver and A.~Townsend}, {\em A fast and well-conditioned spectral
  method}, SIAM Review, 55 (2013), pp.~462--489,
  \url{https://doi.org/10.1137/120865458}.

\bibitem{Osher1988}
{\sc S.~Osher and J.~A. Sethian}, {\em Fronts propagating with
  curvature-dependent speed: {A}lgorithms based on {H}amilton--{J}acobi
  formulations}, Journal of Computational Physics, 79 (1988), p.~12–49,
  \url{https://doi.org/10.1016/0021-9991(88)90002-2}.

\bibitem{papadopoulos2024hierarchicalproximalgalerkinfast}
{\sc I.~P.~A. Papadopoulos}, {\em Hierarchical proximal {G}alerkin: a fast
  $hp$-{FEM} solver for variational problems with pointwise inequality
  constraints}, 2024, \url{https://doi.org/10.48550/arXiv.2412.13733}.

\bibitem{Papadopoulos2024}
{\sc I.~P.~A. Papadopoulos, T.~S. Gutleb, R.~M. Slevinsky, and S.~Olver}, {\em
  Building hierarchies of semiclassical {J}acobi polynomials for spectral
  methods in annuli}, SIAM Journal on Scientific Computing, 46 (2024),
  pp.~A3448--A3476, \url{https://doi.org/10.1137/23M160846X}.

\bibitem{petersen2012matrix}
{\sc K.~B. Petersen and M.~S. Pedersen}, {\em The matrix cookbook},  (2012),
  \url{https://www2.compute.dtu.dk/pubdb/pubs/3274-full.html}.

\bibitem{Prigozhin_1994}
{\sc L.~Prigozhin}, {\em Sandpiles and river networks: {E}xtended systems with
  nonlocal interactions}, Phys. Rev. E, 49 (1994), pp.~1161--1167,
  \url{https://doi.org/10.1103/PhysRevE.49.1161}.

\bibitem{Prigozhin_1996a}
{\sc L.~Prigozhin}, {\em On the {B}ean critical-state model in
  superconductivity}, European Journal of Applied Mathematics, 7 (1996),
  p.~237–247, \url{https://doi.org/10.1017/S0956792500002333}.

\bibitem{Prigozhin_1996}
{\sc L.~Prigozhin}, {\em Variational model of sandpile growth}, European
  Journal of Applied Mathematics, 7 (1996), p.~225–235,
  \url{https://doi.org/10.1017/S0956792500002321}.

\bibitem{Risken1996-wo}
{\sc H.~Risken}, {\em The Fokker-Planck Equation: Methods of Solution and
  Applications}, Springer Berlin Heidelberg, 1996,
  \url{https://doi.org/10.1007/978-3-642-61544-3}.

\bibitem{robinson2017}
{\sc M.~Robinson, C.~Luo, P.~E. Farrell, R.~Erban, and A.~Majumdar}, {\em From
  molecular to continuum modelling of bistable liquid crystal devices.}, Liquid
  Crystals, 44 (2017), pp.~2267--2284,
  \url{https://doi.org/10.1080/02678292.2017.1290284}.

\bibitem{Rockafellar1967Conjugates}
{\sc R.~T. Rockafellar}, {\em Conjugates and {L}egendre transforms of convex
  functions}, Canadian Journal of Mathematics, 19 (1967), pp.~200--205,
  \url{https://doi.org/10.4153/CJM-1967-012-4}.

\bibitem{RTRockafellar_1970}
{\sc R.~T. Rockafellar}, {\em Convex Analysis}, Princeton University Press,
  Princeton, 1970.

\bibitem{Rodrigues2019}
{\sc J.~F. Rodrigues and L.~Santos}, {\em Variational and Quasi-Variational
  Inequalities with Gradient Type Constraints}, Springer International
  Publishing, 2019, p.~319–361,
  \url{https://doi.org/10.1007/978-3-030-33116-0_13}.

\bibitem{santambrogio2015}
{\sc F.~Santambrogio}, {\em Optimal Transport for Applied Mathematicians:
  Calculus of Variations, PDEs, and Modeling}, no.~87 in Progress in Nonlinear
  Differential Equations and Their Application, Birkhäuser, 2015,
  \url{https://doi.org/10.1007/978-3-319-20828-2}.

\bibitem{Savin2013}
{\sc O.~Savin}, {\em Global {$W^{2,p}$} estimates for the {M}onge--{A}mpère
  equation}, Proceedings of the American Mathematical Society, 141 (2013),
  p.~3573–3578, \url{https://doi.org/10.1090/s0002-9939-2013-11748-x}.

\bibitem{schwedes2017mesh}
{\sc T.~Schwedes, D.~A. Ham, S.~W. Funke, and M.~D. Piggott}, {\em Mesh
  dependence in {PDE}-constrained optimisation}, Springer, 2017,
  \url{https://doi.org/10.1007/978-3-319-59483-5}.

\bibitem{JASethian_1996}
{\sc J.~A. Sethian}, {\em A fast marching level set method for monotonically
  advancing fronts}, Proceedings of the National Academy of Sciences of the
  United States of America, 93 (1996), pp.~1591--1595,
  \url{https://doi.org/10.1073/pnas.93.4.1591}.

\bibitem{signorini1959}
{\sc A.~Signorini}, {\em Questioni di elasticit{\`a} non linearizzata e
  semilinearizzata}, Rendiconti Di Matematica e Delle Sue Applicazioni, 18
  (1959), pp.~95--139.

\bibitem{surowiec2023}
{\sc T.~M. Surowiec and S.~W. Walker}, {\em Optimal control of the {Landau--de
  Gennes} model of nematic liquid crystals}, SIAM Journal on Control and
  Optimization, 61 (2023), pp.~2546--2570,
  \url{https://doi.org/10.1137/22m1506158}.

\bibitem{MTeboulle_2018}
{\sc M.~Teboulle}, {\em A simplified view of first order methods for
  optimization}, Mathematical Programming, 170 (2018), p.~67–96,
  \url{https://doi.org/10.1007/s10107-018-1284-2}.

\bibitem{TWTing_1969}
{\sc T.~W. Ting}, {\em Elastic-plastic torsion of convex cylindrical bars},
  Journal of Mathematics and Mechanics, 19 (1969), pp.~531--551,
  \url{http://www.jstor.org/stable/24901809} (accessed 2024-03-06).

\bibitem{ulbrich2003}
{\sc M.~Ulbrich}, {\em Semismooth {N}ewton methods for operator equations in
  function spaces}, SIAM Journal on Optimization, 13 (2003), pp.~805--841,
  \url{https://doi.org/10.1137/S1052623400371569}.

\bibitem{ulbrich2011}
{\sc M.~Ulbrich}, {\em Semismooth {N}ewton Methods for Variational Inequalities
  and Constrained Optimization Problems in Function Spaces}, vol.~11 of
  MOS-SIAM Series on Optimization, SIAM, 2011,
  \url{https://doi.org/10.1137/1.9781611970692}.

\bibitem{VasilDisk}
{\sc G.~M. Vasil, K.~J. Burns, D.~Lecoanet, S.~Olver, B.~P. Brown, and J.~S.
  Oishi}, {\em Tensor calculus in polar coordinates using {J}acobi
  polynomials}, Journal of Computational Physics, 325 (2016), pp.~53--73,
  \url{https://doi.org/10.1016/j.jcp.2016.08.013}.

\bibitem{wachter2006}
{\sc A.~W{\"a}chter and L.~T. Biegler}, {\em On the implementation of an
  interior-point filter line-search algorithm for large-scale nonlinear
  programming}, Mathematical Programming, 106 (2006), pp.~25--57,
  \url{https://doi.org/10.1007/s10107-004-0559-y}.

\bibitem{wang2021}
{\sc W.~Wang, L.~Zhang, and P.~Zhang}, {\em Modelling and computation of liquid
  crystals}, Acta Numerica, 30 (2021), pp.~765--851,
  \url{https://doi.org/10.1017/s0962492921000088}.

\bibitem{wang2022}
{\sc Y.~Wang and J.~Xu}, {\em Q-tensor gradient flow with quasi-entropy and
  discretizations preserving physical constraints}, Journal of Scientific
  Computing, 94 (2022), \url{https://doi.org/10.1007/s10915-022-02060-x}.

\bibitem{warby2003}
{\sc M.~Warby, J.~R. Whiteman, W.-G. Jiang, P.~Warwick, and T.~Wright}, {\em
  Finite element simulation of thermoforming processes for polymer sheets},
  Mathematics and computers in simulation, 61 (2003), pp.~209--218,
  \url{https://doi.org/10.1016/S0378-4754(02)00077-0}.

\bibitem{weiser2005asymptotic}
{\sc M.~Weiser, A.~Schiela, and P.~Deuflhard}, {\em Asymptotic mesh
  independence of {N}ewton's method revisited}, SIAM journal on numerical
  analysis, 42 (2005), pp.~1830--1845,
  \url{https://doi.org/10.1137/S0036142903434047}.

\bibitem{White1983}
{\sc L.~W. White}, {\em Control of a hyperbolic problem with pointwise stress
  constraints}, Journal of Optimization Theory and Applications, 41 (1983),
  p.~359–369, \url{https://doi.org/10.1007/bf00935231}.

\bibitem{Wollner2012}
{\sc W.~Wollner}, {\em Optimal control of elliptic equations with pointwise
  constraints on the gradient of the state in nonsmooth polygonal domains},
  SIAM Journal on Control and Optimization, 50 (2012), p.~2117–2129,
  \url{https://doi.org/10.1137/110836419}.

\bibitem{wu2023geometric}
{\sc K.~Wu and C.-W. Shu}, {\em Geometric quasilinearization framework for
  analysis and design of bound-preserving schemes}, SIAM Review, 65 (2023),
  pp.~1031--1073, \url{https://doi.org/10.1137/21M1458247}.

\bibitem{wu2017multiphase}
{\sc S.~Wu and J.~Xu}, {\em Multiphase {A}llen--{C}ahn and {C}ahn--{H}illiard
  models and their discretizations with the effect of pairwise surface
  tensions}, Journal of Computational Physics, 343 (2017), pp.~10--32,
  \url{https://doi.org/10.1016/j.jcp.2017.04.039}.

\bibitem{SZagatti_2009}
{\sc S.~Zagatti}, {\em On viscosity solutions of {H}amilton--{J}acobi
  equations}, Transactions of the American Mathematical Society, 361 (2009),
  pp.~41--59, \url{http://www.jstor.org/stable/40302759}.

\bibitem{Zeidler1986}
{\sc E.~Zeidler}, {\em Nonlinear Functional Analysis and its Applications: Part
  1: Fixed-Point Theorems}, vol.~1 of Nonlinear Functional Analysis and its
  Applications, Springer, 1986.

\bibitem{Zhao2004}
{\sc H.~Zhao}, {\em A fast sweeping method for eikonal equations}, Mathematics
  of Computation, 74 (2004), p.~603–627,
  \url{https://doi.org/10.1090/s0025-5718-04-01678-3}.

\bibitem{Zhao2000}
{\sc H.-K. Zhao, S.~Osher, B.~Merriman, and M.~Kang}, {\em Implicit and
  nonparametric shape reconstruction from unorganized data using a variational
  level set method}, Computer Vision and Image Understanding, 80 (2000),
  p.~295–314, \url{https://doi.org/10.1006/cviu.2000.0875}.

\end{thebibliography}

\end{document}